\PassOptionsToPackage{czech,english}{babel}
\documentclass[a4paper, 12pt]{article} %define title
\usepackage{natbib}
\usepackage[linesnumbered, ruled,vlined]{algorithm2e}
\usepackage{textcomp}
\usepackage{amsmath}
\usepackage{authblk}
\usepackage{mathtools}
\usepackage{enumitem} 
\usepackage{amssymb}
\usepackage{amsthm}
\usepackage{setspace}
\usepackage{mleftright}
\usepackage{multirow}
\usepackage[figuresright]{rotating}
\usepackage[utf8]{inputenc}
\usepackage[english]{babel}
\usepackage{filecontents}
\usepackage{graphicx}
\usepackage{titletoc}
\usepackage[toc,page]{appendix}
\graphicspath{ {figures/} }
\usepackage{array}
\usepackage{float}
\usepackage{scalerel,stackengine}
\usepackage{accents}
\usepackage{color, graphics, graphicx}
\usepackage{xcolor,colortbl}
\usepackage{undertilde}
\usepackage[T1]{fontenc}
\usepackage{lmodern}
\usepackage[babel=true]{microtype}
\usepackage[czech]{babel}
\usepackage{hyperref}

\usepackage{etoolbox}
\preto\tabular{\shorthandoff{-}}

\usepackage{caption}
\usepackage{subcaption}

\newtheorem{assumption}{Assumption}
\newtheorem{theorem}{Theorem}[section]
\newtheorem{proposition}{Proposition}[section]
\newtheorem{corollary}{Corollary}[section]
\newtheorem{lemma}[theorem]{Lemma}
\numberwithin{equation}{section}
\numberwithin{theorem}{section}

\newcommand{\mbf}[1]{\mbox{\boldmath $#1$}}
\newcommand{\pr}{^{\prime}}
\newcommand{\n}{^{(n)}}
\newcommand{\bth}{{\mbf \theta}}

\newcommand{\bTH}{{\mbf \Theta}}

\newcommand{\al}{{\mbf \alpha}}

\newcommand{\btau}{{\mbf \tau}}

\DeclareMathAlphabet\mathbfcal{OMS}{cmsy}{b}{n}

\usepackage{stackengine}
\stackMath
\newcommand\tenq[2][1]{%
\def\useanchorwidth{T}%
\ifnum#1>1%
\stackunder[0pt]{\tenq[\numexpr#1-1\relax]{#2}}{\scriptscriptstyle\thicksim}%
\else%
\stackunder[1pt]{#2}{\scriptscriptstyle\thicksim}%
\fi%
}

\def\pms{\mspace{-1mu}{\scriptscriptstyle\pm}}

\def\s{\sum_{t=1}^n}

\def\bth{\mbox{\boldmath$\theta$}}

\def\bLam{\mbox{\boldmath$\Lambda$}}

\def\bDelta{\mbox{\boldmath$\Delta$}}
\def\bSigma{\mbox{\boldmath$\Sigma$}}
\def\bXi{\mbox{\boldmath$\Xi$}}
\def\bmu{\mbox{\boldmath$\mu$}}
\def\bxi{\mbox{\boldmath$\xi$}}
\def\bpsi{\mbox{\boldmath$\psi$}}
\def\bGamma{\mbox{\boldmath$\Gamma$}}

\def\bUpsilon{\mbox{\boldmath$\Upsilon$}}
\def\bvp{\mbox{\boldmath$\varphi$}}
\def\bepsilon{\mbox{\boldmath$\epsilon$}}

\def\vp{\varphi}

\def\vp{\varphi}

\def\R{{\mathbb R}}
\def\Z{{\mathbb Z}}
\def\D{{\mbf D}}
\def\N{{\mathbb N}}
\def\a{{\mbf a}}
\def\b{{\mbf b}}

\def\vbf{{\mbf v}}
\def\K{{\mbf K}}

\def\J{{\mbf J}}
\def\L{{\mbf L}}
\def\u{{\mbf u}}

\def\M{{\mbf M}}
\def\m{{\mbf m}}

\def\0{{\mbf 0}}

\def\P{{\mbf P}}
\def\Q{{\mbf Q}}
\def\S{{\mbf S}}

\def\T{{\mbf T}}
\def\Y{{\mbf Y}}
\def\A{{\mbf A}}

\def\CC{{\mbf C}}
\def\C{{\mbf C}}
\def\DD{{\mbf D}}
\def\F{{\mbf F}}
\def\H{{\mbf H}}
\def\hbf{{\mbf h}}
\def\G{{\mbf G}}
\def\I{{\mbf I}}
\def\X{{\mbf X}}
\def\x{{\mbf x}}
\def\y{{\mbf y}}

\def\ZZ{{\mbf Z}}

\def\0{{\mbf 0}}

\def\z{{\mbf z}}

\def\1{{\mbf 1}}
\newcommand{\Acal}{{\mathcal A} }
\newcommand{\Bcal}{{\mathcal B} }

\newcommand{\betab}{{\pmb \beta} }

\newtheorem{innercustomgeneric}{\customgenericname}
\providecommand{\customgenericname}{}
\newcommand{\newcustomtheorem}[2]{%
  \newenvironment{#1}[1]
  {%
   \renewcommand\customgenericname{#2}%
   \renewcommand\theinnercustomgeneric{##1}%
   \innercustomgeneric
  }
  {\endinnercustomgeneric}
}

\newcustomtheorem{customlemma}{Lemma}

\title{ \Large \bf \sc  $\,$\vspace{-25mm}\\%Center-Outward Rank  Tests
Rank-Based Testing\vspace{-1mm} \\ for %of the Order of \\ a
 Semiparametric VAR Models: \vspace{-1mm} \\ a measure transportation approach\vspace{-3mm} }

\author{M. Hallin,$^{(a)}$ D. La Vecchia,$^{(b)}$  and H. Liu$^{(c)}$\vspace{-1mm}   \\
$^{(a)}${\small{ECARES, Université libre de Bruxelles CP 114/4 \\Avenue F.D. Roosevelt 50 - 
B-1050 Bruxelles, Belgium \\
 Email: mhallin@ulb.ac.be \\ \vspace{0.1cm}
 
 $^{(b)}$Research Center for Statistics, University of Geneva \\ Boulevard du Pont d'Arve 40 - CH-1211 Geneva, Switzerland \\
 Email: davide.lavecchia@unige.ch \\ \vspace{0.1cm}

 $^{(c)}$Department of Mathematics and Statistics, Lancaster University \\ LA1 4YF Lancaster, UK\vspace{-2mm} \\
 Email: h.liu11@lancaster.ac.uk 
    } }}
    \date{}\vspace{-4mm}

\begin{document}

\maketitle

\vspace{-10mm}

\begin{abstract}
{\small{We develop a class of tests for semiparametric vector autoregressive (VAR) models with unspecified innovation densities,   based on the recent  measure-transportation-based concepts of multivariate {\it center-outward ranks} and {\it signs}.  We show that these concepts, combined with  Le Cam's asymptotic  theory of statistical experiments, yield novel testing procedures, which (a)~are valid under a broad class of    innovation densities (possibly non-elliptical, skewed,  and/or  with infinite moments), (b)~are optimal (locally asymptotically maximin or most stringent) at selected %spherical 
ones, and (c) are robust against additive outliers. In order to do so, we establish a  H\' ajek asymptotic representation result, of independent interest,  for a general class of center-outward rank-based serial statistics. 
 As an illustration, we   consider the problems of testing the absence of serial correlation in multiple-output and possibly non-linear regression (an extension of the classical Durbin-Watson problem) and the sequential identification of  the order $p$ of a vector autoregressive (VAR($p$)) model. A Monte Carlo comparative study of our tests and their routinely-applied Gaussian competitors  demonstrates the benefits (in terms of size, power, and robustness) of our methodology; these benefits are particularly significant in the presence of asymmetric and leptokurtic innovation densities. A real data application concludes the paper.}\vspace{-1mm}}
\end{abstract}

\section{Introduction}  Despite the overwhelming empirical evidence,  in most real-life datasets,  of non-Gaussian and non-symmetric (non-elliptical) distributions, multivariate  time series analysis  remains very deeply marked by explicit or implicit  Gaussian assumptions: correlogram-based and spectral methods, pseudo-Gaussian tests,  Gaussian quasi-likelihoods, and their possible robustifications, are ubiquitous in methodological developments  as well as in daily practice.  

The theoretical justification for this   prevalence of Gaussian methods  looks quite solid at first sight and is seldom questioned: Gaussian quasi-likelihood methods 
%---such as pseudo-Gaussian tests and Gaussian quasi-maximum-likelihood estimators (QMLEs)---
are {\it asymptotically valid} (viz., pseudo-Gaussian tests have correct asymptotic nominal size and Gaussian quasi-maximum-likelihood estimators (QMLEs)  are root-$n$ consistent) under a broad range $\mathcal F$ of non-Gaussian innovation densities $f$ (typically, under  finite fourth-order moments).  

A closer look at this asymptotic validity argument, however, reveals a fundamental weakness:  while it holds,  pointwise, for any given $f\in{\mathcal F}$, it fails to hold {\it uniformly} over~$\mathcal F$. Consider, for instance, the sequence $\phi^{(n)}$ of Gaussian level-$\alpha$ tests rejecting a null hypothesis of the form~$\bth=\bth_0$ whenever some test statistic $S\n$ exceeds the standard normal quantile~$z_\alpha$ of order $(1-\alpha)$. Denoting by ${\rm P}^{(n)}_{\bth ,  f}$ the distribution of the observation under parameter value~$\bth$ and  innovation density~$f\in{\mathcal F}$,  assume that $\phi^{(n)}$ is a valid pseudo-Gaussian test at asymptotic level $\alpha$, that is, satisfies 
\begin{equation}\label{pseudoG}
%$
\lim_{n\to\infty}{\rm P}^{(n)}_{\bth_0 , f}\left[S\n > z_\alpha%\phi^{(n)}=1
\right] \leq\alpha%$ 
\end{equation} 
 for any $f\in\mathcal F$. This, however,  is not sufficient if $f$ remains unspecified. The null, then,  is the    semiparametric hypothesis~$\big\{
{\rm P}^{(n)}_{\bth ,  f}\vert\, \bth=\bth_0,\,  f\in\mathcal F\big\}$ and, in order for $\phi^{(n)}$ to qualify   as an asymptotically valid  level-$\alpha$ test,    the stronger condition 
\begin{equation}\label{UpseudoG}
%$
\lim_{n\to\infty}\sup_{f\in\mathcal F}{\rm P}^{(n)}_{\bth_0 , f}\left[S\n > z_\alpha%\phi^{(n)}=1
\right]  \leq\alpha 
%$
 \end{equation}
which, in general, does not follow from \eqref{pseudoG},  is required. 

This lack of uniformity in the asymptotics of pseudo-Gaussian procedures is not just a mathematical detail: in the context of VAR estimation, \cite{HLL2019} have shown that, depending on the actual innovation density $f$,  the finite-sample performance of a valid QMLE, for given sample size $n$, can be extremely far from its asymptotic performance. That discrepancy between finite-sample and asymptotic performance is particularly significant (see the motivating example in the introduction of   the same reference) in the presence of  non-elliptical innovation densities such as   mixture or skewed densities---a type of innovation distribution which is likely to occur in the case of  omitted variables or covariates---and only gets worse as the dimension of the observation space  increases. The same phenomenon is bound to take place in testing problems.

A remedy to this lack of uniformity of pseudo-Gaussian asymptotics is distribution-freeness. If indeed the test statistic $S\n$ is distribution-free under $\big\{
{\rm P}^{(n)}_{\bth ,  f}\vert\, \bth\in\bTH_0,\,  f\in\mathcal F\big\}$,  then \eqref{pseudoG} and \eqref{UpseudoG} are equivalent. Rank-based inference then naturally enters the picture: it follows indeed from Basu's classical theorems on the relation between minimal sufficiency and distribution-freeness (see Appendix~E in~\cite{Hallinetal2020} or Chapters~5 and~10 in~\cite{LehRom}) that, in the univariate case, the $\sigma$-field generated by residual ranks  is {\it essentially}\footnote{That is, maximal up to $\mathcal F$-null sets.} {\it maximal ancillary}---intuitively, {\it ``maximal distribution-free"}. 

Moreover, it follows from \cite{HW03} that the semiparametric efficiency bounds at given density$f\in{\mathcal F}$---which, for ARMA models, coincide with the parametric ones---can be reached by rank-based tests: normal-score (van der Waerden) rank tests in ARMA models, in particular, are achieving, under Gaussian $f$, the same asymptotic performance as the  (pseudo-)Gaussian ones while doing strictly and uniformly better under non-Gaussian $f$  \citep{HJTSA}. Pseudo-Gaussian tests, thus, are not admissible and uniformly dominated by their normal-score rank-based counterparts. 

All this, which strongly suggests abandoning pseudo-Gaussian methods in favor of rank-based ones, unfortunately, until recently was limited to univariate ARMA models due to the lack of   an adequate multivariate extension of the univariate concept of ranks. Many proposals have been made in the literature---see \cite{Hallinetal2020} for a commented bibliography. None of them, however, is enjoying distribution-freeness---let alone maximal ancillarity---except for the so-called {\it Mahalanobis ranks} and {\it signs} \citep{HP02, HP04}  under the quite restrictive assumption, however,  of elliptical symmetry.  This lack of a satisfactory concept of multivariate ranks   has been dealt with recently  with the introduction of the measure-transportation-based  {\it Monge-Kantorovich} \citep{Chernetal17} and {\it center-outward}  \citep{Hallin2017} {\it ranks and signs},  the maximal ancillarity of which is established in~\cite{Hallinetal2020}. Building on that concept, \cite{HHH20} are proposing distribution-free center-outward rank-based tests for multiple-output regression and MANOVA models, \cite{GhS19} for goodness-of-fit, \cite{DebSen19} and \cite{Shietal20}   for independence between vectors, while  
 \cite{HLL2019}  are constructing    center-outward  R-estimators for VARMA models. The finite-sample performance of the new rank-based methods in all cases appears to be quite remarkable.  In this paper, we similarly study a class of center-outward rank-based testing procedures for VAR models and illustrate their role in the problem of sequential VAR order selection (see P\" otscher~1983) in the presence of unspecified innovation density.

The paper is organized as follows. Section~\ref{Sec: meth} introduces notation and recalls some methodological facts about local asymptotic normality for VARs (Section~\ref{secLAN}) and the measure-transportation approach to center-outward ranks and signs (Section~\ref{secranks}).  Section~\ref{SecAsRepr} contains the main theoretical result of this paper: a H\' ajek representation theorem for serial statistics based on center-outward ranks and signs. This key result, of independent interest, extends to general scores the spherical-score result of \cite{HLL2019}, as  well as the univariate results of \cite{HIP85} and \cite{HV96}. Based on LAN and the H\' ajek representation of Section~\ref{SecAsRepr}, Section~\ref{SecSR} introduces the notion of rank-based central sequence from which the rank-based test statistics are  constructed in Section~\ref{Sec:Rtests}. Focusing on two particular problems, we develop, in Section~\ref{Sec.theta0}, a rank-based  extension  of the classical Durbin-Watson procedure to multiple-output and possibly non-linear regression and, in Section~\ref{Sec.testorder} a sequential rank-based procedure for VAR order identification. Section~\ref{GaussSec} derives the Gaussian counterparts of the tests developed in Section~\ref{Sec:Rtests} while  Section~\ref{Sec:NumStud} presents a Monte Carlo investigation of the finite-sample performance of the rank-based methods we are proposing (Sections~\ref{sec.NumWhite} and~\ref{VARorder}) and concludes with a real-data application (Section~\ref{Sec:emp}).
\color{black}

\section{Methodological background\label{Sec: meth}}

Our approach is combining two basic methodological ingredients:
\begin{enumerate}
\item[(a)] the   local asymptotic normality of VAR models, as established in  \cite{GH95} and \cite{HP04}; the explicit statement (Section~\ref{secLAN} and Proposition~\ref{Prop.LAN1}) of this property (with closed-form expressions\footnote{The difficulty in the derivation of these closed forms mainly stems from the fact that the VAR order $p_0$ under the null hypothesis can be strictly less than the order $p$ under the alternative.} for central sequences and information matrices)  
%. Explicit expressions of the corresponding central sequences
 requires some  algebraic preparation,  which  is  provided Appendix~A and can be skipped at first reading;  %  borrowed from \cite{GH95} and \cite{HP04}.
\item[(b)]  the recently proposed (\cite{Chernetal17, Hallin2017, Hallinetal2020}) measure-transportation-based concepts of multivariate {\it center-outward ranks} and {\it signs}: see Section~\ref{secranks} for precise definitions . 
% and Appendix~B for an intuitive presentation.
\end{enumerate}
In this section, we first settle the notation (Section~\ref{sec2}) before turning  (Sections~\ref{secLAN}  and~\ref{secranks}) to~(a) and~(b). 

\subsection{Notation and general setting}\label{sec2}

Consider the $d$-dimensional VAR model of order $p\leq p_1$% ($p_1 \geq 1$)  
\begin{equation}
\A_{\bth}(L)\X_t := (\I_d - \sum_{i=1}^{p_1} \A_i L^i)  \X_t = \bepsilon_t, \quad t \in \mathbb{Z}, \label{VAR_mod}
\end{equation}
where $\I_d$ is the $d\times d$ identity matrix, $\A_1, \ldots , \A_{p_1}$ are $d \times d$ matrix coefficients, $L$ is the lag operator,~$\{\bepsilon_t; t \in \Z\}$ is an i.i.d.\  innovation  process with mean $\0$ and probability density~$f$, and~$\bth := ((\text{vec}{\A_1})^\prime, \ldots , (\text{vec}{\A_{p_1}})^\prime)^\prime$   
denotes the $({p_1}d^2)$-dimensional vector of autoregressive parameters. 
%Writing  
%$$\bth := ((\text{vec}{\A_1})^\prime, \ldots , (\text{vec}{\A_{p_1}})^\prime)^\prime$$  
%for the autoregressive parameter under $({p_1}d^2)$-dimensional vector form, let
For $\bth$ of the form 
 \begin{equation}\label{bth0}\bth_0 := ((\text{vec}{\A_1})^\prime, \ldots , (\text{vec}{\A_{p_0}})^\prime, \0_{d^2(p_1-p_0)\times 1}^\prime)^\prime 
 \end{equation}
  (with $0< p_0<p_1$),  $\A_{\bth}(L)$  is the VAR operator
\begin{equation}% \A_0(L) := %\A_0
\A_{\bth_0}(L)=\I_d - \sum_{i = 1}^{p_0} \A_i L^i\label{VAR_mod0}
\end{equation}
 of order $p_0<p_1$.  For   $p_0=0$,   \eqref{bth0}, with $\bth_0=(0,\ldots,0)$ yields the white noise model~$ \X_t = \bepsilon_t$.

 Denote by $\bTH$ the set of parameter values $\bth\in\mathbb{R}^{{p_1}d^2}$ such that all solutions of the determinantal equation
$\text{\rm det} \left(\A_{\bth} (z) 
\right)
%=\text{\rm det} \left(\I_d - \sum_{i = 1}^{p_0} \A_i z^i \right)
 = 0$, $z \in \mathbb{C}$ 
 lie outside the unit ball (the standard VAR stationarity  condition) and by $\bTH_0$ the set of parameter values $\bth_0\in\bTH$ of the form~\eqref{bth0} with~$|\A_{p_0}| \neq 0$. The notation $\A_0(L) $  will be used for  a null VAR$(p_0)$ operator of the form~\eqref{VAR_mod0} with unspecified $\bth_0\in\bTH_0$, the notation  $\A_1(L) $ for an alternative VAR$(p)$ operator~$\A_{\bth_1}(L) $ of order $p$, that is, with unspecified   $\bth_1\in\bTH$ and $\A_p\neq {\bf 0}$ for some $p_0~\!<~\!p\leq~\!p_1$. Denoting by~$\X^{(n)}:=\{\X^{(n)}_1, \ldots , \X^{(n)}_n\}$  (superscript~$^{(n)}$  omitted whenever possible) an  observed triangular array of solutions of \eqref{VAR_mod}, our objective is to construct   rank-based tests of (a) the null hypothesis ${\mathcal H}\n_{\bth_0}$ under which $\X^{(n)}$ was generated by a   VAR operator of the form~$\A_{\bth_0}(L)$ (specified ${\bth_0}\in\bTH_0$, unspecified innovation density) and (b)  the null hypothesis ${\mathcal H}\n_{0}$ under which $\X^{(n)}$ was generated by a   VAR operator of the form $\A_{0}(L)$ (unspecified ${\bth_0}\in\bTH_0$, unspecified innovation density) against the alternative of a VAR operator of the form $\A_1(L) $ (unspecified $\bth_1\in\bTH_1$, unspecified innovation density). The sequential  order identification procedure then consists of performing such tests for $p_0= 0, 1, ...$ and~$p_1=p_0+1$ until the null hypothesis ${\mathcal H}\n_{0}$ no longer gets rejected.

\subsection{Local asymptotic normality}\label{secLAN}

In this section, we state the local asymptotic normality (LAN) result on which    the construction of our tests heavily relies.  LAN, of course, requires some regularity assumptions on $f$; we throughout assume the following (essentially borrowed from  \cite{GH95}).
% \medskip 

\begin{assumption}\label{ass.den}{\rm 
\begin{enumerate}
\item[{\it (i)}] The innovation density $f$ belongs to the class $\mathcal{F}_d$ of {\it non-vanishing} (with respect to the Lebesgue measure~$\mu_d$ on~$\R^d$) densities\footnote{The requirement that $f$ has support $\mathbb{R}^d$ can be relaxed to a requirement of a convex support (see \cite{BSH20}) at the expense, however, of a less direct definition of center-outward distribution and quantile functions. For the sake of simplicity, we are sticking to the assumption made here. } i.e., for all $c \in \R^+\!$, there exist~$b_{c; f}$ and~$a_{c; f}$ in $\mathbb{R}$ such that 
$0<b_{c; f}  \leq a_{c; f}<\infty$ and $b_{c; f} \leq f(\x) \leq a_{c; f}$  for any~$\Vert \x\Vert  \leq c$;
\item[{\it (ii)}]  $\displaystyle{\int \x f(\x) \mathrm{d}\mu = \0}$ and $\displaystyle{ \int \x \x^\prime f(\x) \mathrm{d}\mu = \bXi}$ where~$\bXi$ is positive definite; 
\item[{\it (iii)}] $f^{1/2}$ is mean-square differentiable with mean-square gradient~$\DD f^{1/2}$, that is, 
there exists a square-integrable   vector $\DD f^{1/2}$ such that, for all sequence~$\hbf\in\mathbb{R}^d$ such\linebreak that~$\0 \neq \hbf \rightarrow \0$, 
$\displaystyle{(\hbf^\prime \hbf)^{-1} \int \left[ f^{1/2}(\x + \hbf) - f^{1/2}(\x) - \hbf^\prime \DD  f^{1/2} (\x)\right]^2 \mathrm{d}\mu \rightarrow 0;}$
\item[{\it (iv)}] letting (the location score function) 
%\begin{equation}\label{bvpf}
$\bvp_{f}(\x) := (\vp_1(\x), \ldots , \vp_d(\x))^\prime := -2 (\DD f^{1/2})/f^{1/2}$, %\end{equation}
 $$\int [\vp_i (\x)]^4 f(\x) \mathrm{d}\mu < \infty,\ i = 1, \ldots , d;$$
\item[{\it (v)}]   the   function $\x\mapsto\bvp_f(\x)$ is piecewise Lipschitz, i.e., there exists a finite measurable partition of $\R^d$ into $J$ non-overlapping subsets $I_1,  \ldots , I_J$ such that,  for all $\x, \y$\linebreak in~$I_j$, $j = 1, \ldots , J$,~$\Vert \bvp_f(\x) - \bvp_f(\y)\Vert  \leq K \Vert \x - \y\Vert $.
\end{enumerate}}
\medskip
\end{assumption}
%These assumptions guarantee {\it (i)} the finiteness of Fisher information for $\bth$
Denote by ${\mathcal F}_d^*$   the family of innovation densities satisfying Assumption~\ref{ass.den}.\medskip 

Let $\ZZ_1^{(n)}(\bth_0),\ldots,\ZZ_n^{(n)}(\bth_0)$ with
 $\ZZ_t^{(n)}(\bth_0):= \A_{\bth_0}(L)\X^{(n)}_t
$ 
\color{black} denote the residuals computed 
from the initial values%$\bepsilon_{-q+1}, \ldots , \bepsilon_0$ and
~$\X_{-p + 1}, \ldots , \X_{0}$, 
the parameter value $\bth_0$, and the observations~$\X^{(n)}$.  Clearly, $\X^{(n)}$ is the finite realization of a solution of~\eqref{VAR_mod0} with parameter value~$\bth_0$  iff the residuals~$\ZZ^{(n)}_1(\bth_0), \ldots , \ZZ^{(n)}_n(\bth_0)$ 
and~$\bepsilon_1, \ldots , \bepsilon _n$ coincide. Denoting by ${\rm P}^{(n)}_{\bth_0 ;f}$ the distribution of~$\X^{(n)}$ under parameter value~$\bth_0$ and innovation density~$f$, the residuals $\ZZ_1^{(n)}(\bth_0),\ldots,\ZZ_n^{(n)}(\bth_0)$ 
under~${\rm P}^{(n)}_{\bth_0 ;f}$ thus are i.i.d.\ with density $f$.

Write~$L^{(n)}_{\bth_0 + n^{-1/2}\btau^{(n)}/\bth_0; f}$ for the log-likelihood ratio $\log {\rm d}{\rm P}^{(n)}_{\bth_0 + n^{-1/2}\btau^{(n)};f} / {\rm d}{\rm P}^{(n)}_{\bth_0;f}$, 
%of ${\rm P}^{(n)}_{\bth_0 + n^{-1/2}\btau^{(n)};f}$
% % {\mathcal H}^{(n)}(\bth_0 + n^{-1/2}\btau^{(n)};f)$
%  with respect to~${\rm P}^{(n)}_{\bth_0;f}$, 
  where~$\btau^{(n)} $ is a bounded sequence of $\mathbb{R}^{p_1 d^2}$.   
  %The following LAN result  will be used to motivate the definition of our R-estimator and to establish the asymptotic normality of our R-estimators. Note that it does not require $f$ to be elliptic. 
%As remarked in Garel and Hallin (1995), the above assumptions (i) and (iii) guarantee that the innovation term has a finite second moment and existence of the Fisher information matrix $\mathbfcal{I}(f) = \int \bvp_f(\x) \bvp_f^\prime(\x) f(\x) d\mu$, which appears in the covariance matrix of the central sequence in following Proposition \ref{Prop.LAN1}.
Define  
%I assume these (or at least the product $\M_{\bth_0}^\prime  \P_{\bth_0}^\prime \Q_{\bth_0}^{(n)\prime}$) do not depend on the $\psi$'s ... correct?
\begin{equation}\label{Delta}
\bDelta^{(n)}_{f} (\bth_0) :=   n^{1/2}   \M_{\bth_0}^\prime  \P_{\bth_0}^\prime \Q_{\bth_0}^{(n)\prime} \bGamma_{f}^{(n)}(\bth_0),
\end{equation}
where $\M_{\bth_0}$, $ \P_{\bth_0}$, and $\Q^{(n)}_{\bth_0}$ are given in Appendix~A (\eqref{defM} and \eqref{defPQ}) and 
\begin{equation}\label{Sf1}
\bGamma_{f}^{(n)}(\bth_0) :=   n^{-1/2}  \left((n-1)^{1/2} (\text{vec}{\bGamma_{1, f}^{(n)}(\bth_0)})^\prime, \ldots , (n-i)^{1/2} (\text{vec}{\bGamma_{i, f}^{(n)}(\bth_0)})^\prime, \ldots , (\text{vec}{\bGamma_{n-1, f}^{(n)}(\bth_0)})^\prime\right)^\prime,
\end{equation}
with the so-called  $f$-{\it cross-covariance matrices}
\begin{equation}\label{Gamma1}
\bGamma_{i, f}^{(n)}(\bth_0) := (n-i)^{-1} \sum_{t=i+1}^n \bvp_{f}(\ZZ_t^{(n)}(\bth_0)) \ZZ_{t-i}^{(n)\prime} (\bth_0).
\end{equation}
%($\bvp_{f}$ as in \eqref{bvpf}). 
%Recall that under Assumption (A1), the Green matrices $\G_u$ and $\H_u$ (see Appendix~A) decrease exponentially in $u$: we can thus safely 
Finally, let 
%\begin{equation}\label{bLamf.cov}
$$\bLam_{f} (\bth_0) := \M_{\bth_0}^\prime  \P_{\bth_0}^\prime \underset{n\rightarrow \infty}{\lim} \left\lbrace \Q_{\bth_0}^{(n)\prime} [\I_{n-1} \otimes (\bXi \otimes \mathbfcal{I}(f))] \Q_{\bth_0}^{(n)} \right\rbrace  \P_{\bth_0} \M_{\bth_0}.$$
%\end{equation}
%due to the convergence of the limit.
We then have the following LAN result \citep{HLL2019}.
%The following LAN result then is essentially the same as in Garel and Hallin (1995, (LAN~2) in their Proposition 3.1)%---exploiting, however, the proof of Proposition~1 in Hallin and Paindaveine~(2004) in order to obtain the form \eqref{Delta} of $\bDelta^{(n)}_{f} (\bth_0)$.

\begin{proposition}\label{Prop.LAN1}
Let Assumption~\ref{ass.den} hold. Then, for any $\bth_0\in\bTH_0$ and any bounded sequence~$\btau^{(n)}$ in~$\R^{p_1 d^2}\!$, under ${\rm P}^{(n)}_{\bth_0;f}$, as $n \rightarrow \infty$, 
\begin{equation}\label{lik1}
L^{(n)}_{\bth_0 + n^{-1/2}\btau^{(n)}/\bth_0; f} = \btau^{{(n)}\prime} \bDelta^{(n)}_{f}(\bth_0) - \frac{1}{2} \btau^{(n)\prime} \bLam_{f} (\bth_0) \btau^{(n)} + o_{\rm P}(1),
\end{equation}
and
$\bDelta^{(n)}_{f} (\bth_0) \rightarrow {\mathcal N}(\0, \bLam_{f} (\bth_0)).$
\end{proposition}

%%%%%%%%%%%%%%%%%%%%%%%%%%%%%%%%%%%%%%%%%%%%
%%%%%%%%%%%%%%%%%%%%%%%%%%%%%%%%%%%%%%%%%%%%
%\end{document}%%%%%%%%%%%%%%%%%%%%%%%%%%%%%%%%%%%%%%%%%%%%
%%%%%%%%%%%%%%%%%%%%%%%%%%%%%%%%%%%%%%%%%%%%
%%%%%%%%%%%%%%%%%%%%%%%%%%%%%%%%%%%%%%%%%%%%

Parametrically efficient (in the H\' ajek-Le Cam asymptotic sense) rank- and sign-based (hence, distribution-free) inference procedures 
in LAN families (with given $f$) are possible when the LAN central sequence~\eqref{Delta} can be expressed in terms of signs and ranks. More precisely (a  central sequence~$\bDelta^{(n)}_{f} (\bth_0)$ is only defined up to $o_{{\rm P}}(1)$-under-${\rm P}^{(n)}_{\bth_0;f}$  terms), when there exists some~$\utilde{\bDelta}^{(n)}_{f} (\bth_0)$ measurable with respects to the ranks and signs of the residuals~$\ZZ^{(n)}_t(\bth_0)$ such that $\utilde{\bDelta}^{(n)}_{f} (\bth_0)-\bDelta^{(n)}_{f} (\bth_0)$ is~$o_{{\rm P}}(1)$ under ${\rm P}^{(n)}_{\bth_0;f}$ as $n\to\infty$. 

In the univariate setting, this idea was exploited  in \citet{HPu88,HPu91,HPu94}, \cite{HMel88}, and \cite{GH99}, leading to a fairly complete toolbox of distribution-free procedures based on traditional ranks or signed ranks for ARMA models, in  \cite{BenH92, BenH96} and \cite{HL2017, HL2020} for nonlinear time series models.  In the multivariate context,
under the assumption of ellipticity, \cite{HP02, HP05, HP06, HP08}  for location, regression, and shape, \cite{HPV10, HPV13, HPV14}  for principal and common principal components, and\linebreak  \cite{HP04, HP04b, HP05} for VAR and VARMA models, proposed test procedures   based on the so-called {\it Mahalanobis ranks} and {\it signs}. However,  the assumption of elliptic innovation density (which reduces to symmetry in the univariate case)  severely restricts  the applicability of these tests, specially in the analysis of economic data, where  actual innovations, typically, are skewed and/or leptokurtic. Thus, with the aim of accommodating a more general and realistic setting of non-elliptical distribution,  we develop here novel test procedures based on  the notions of multivariate ranks and signs proposed by \cite{Chernetal17} (under the name of {\it Monge-Kantorovich ranks} and {\it signs}) and developed in \cite{Hallin2017} and \cite{Hallinetal2020} under the name of {\it center-outward ranks} and {\it signs}. These new notions  hinge on measure transportation theory; in their empirical version, they are based on the idea of an optimal  %(in the sense of a quadratic loss function) 
 coupling of the residuals $\{\ZZ^{(n)}_t; 1\leq t \leq n\}$ with a regular grid $\mathfrak{G}\n$  over the open unit ball~${\mathbb S}_d$  in~$\mathbb{R}^d$, which we now describe. %; see Section~\ref{secranks} for mathematical definitions and  graphical representations. %We  briefly review them in section \ref{secranks}.} %We refer to Villani (2009) for a book-length mathematical discussion  and to Hallin (2017) for a discussion in a statistical context.}  

\subsection{Center-outward ranks and signs}\label{secranks}

% \subsubsection{Mathematical aspects} \label{Sec: mathdef}
 
Let $\mathcal{P}_d$ denote the family of all distributions~$\rm P$ with densities in ${\mathcal F}_d$. For $\rm P$ in  this family,  the center-outward distribution functions defined below  are continuous: see \cite{Hallinetal2020}. More general cases are studied in \cite{BSH20}, but require  more cautious and less intuitive definitions which, for the sake of simplicity, we do not consider here. Denote by~${\rm U}_d$ the spherical uniform distribution over  ${\mathbb S}_d$, that is, the product of a uniform measure over the hypersphere ${\mathcal S}_{d-1}$ and a uniform over the unit interval of distances to the origin. The {\it center-outward distribution function}~$\F_{{\pms}}$ of $\rm P$ is defined as the a.e.~unique gradient of convex function mapping $\R^d$ to $\mathbb{S}_{d}$  and  {\it pushing~$\rm P$ forward} to    ${\rm U}_d$.\footnote{That is, such that $\F_{{\pms}}({\bf X})\sim{\rm U}_d$ if ${\bf X}\sim{\rm P}.$} For ${\rm P}\in{\mathcal P}_d$, such mapping is  a homeomorphism between ${\mathbb S}_d\setminus\{{\bf 0}\}$ and $\mathbb{R}^d\setminus \F_{{\pms}}^{-1}(\{{\bf 0}\})$ (\citealt{Figalli2018}) and   the corresponding {\it center-outward quantile function} is defined  (letting, with a small abuse of notation,~$\Q_{\pms} ({\bf 0}) := \F_{\pms}^{-1}(\{{\bf 0}\})$) as~$\Q_{\pms} := \F_{\pms}^{-1}$. For any given distribution $\rm P$, $\Q_{{\pms}}$ induces a collection of continuous, connected, and nested quantile  contours and regions;   the {\it center-outward median}  $\Q_{\pms}(\0)$ is a uniquely defined  compact set of Lebesgue measure zero. We refer to \cite{Hallinetal2020} for details.

 Turning to the sample,  
% denote by ${\rm P}^{(n)}_{\bth_0; f}$ the distribution of $\X^{(n)}$ under parameter value $\bth_0$ and innovation density~$f$. T
 the residuals $\ZZ_1^{(n)}(\bth_0),\ldots,\ZZ_n^{(n)}(\bth_0)$ under ${\rm P}^{(n)}_{\bth_0 ;f}$ are i.i.d.\ with density $f\in{\cal F}_d$ and center-outward 
distribution function $\F_\pm$. 
% (when no confusion arises, we drop the dependence on $\bth_0$). 
 For the empirical counterpart~$\F^{(n)}_\pm$ of $\F_\pm$,   let $n$ factorize into $n=n_R n_S +n_0,$ for $n_R, n_S, n_0 \in \mathbb{N}$ and~$0\leq n_0 < \min \{ n_R, n_S \}$, where $n_R \rightarrow \infty$ and $n_S \rightarrow \infty$ as $n \rightarrow \infty$, and consider a sequence $\mathfrak{G}\n$ of grids, where each grid
%contains  $n_R n_S$ points in the unit ball $\mathbb{S}_d$ and
consists of the intersection between an $n_S$-tuple $(\boldsymbol{u}_1,\ldots \boldsymbol{u}_{n_S})$ of unit vectors, and the~$n_R$ hyperspheres  with radii $1/(n_R+1),\ldots ,n_R/(n_R+1)$ centered at the origin, along with~$n_0$ copies of the origin. The only requirement\footnote{One exception will be made in the sequel when the so-called {\it sign test scores}  are considered (see Section~\ref{Sec: examples}). Those scores being entirely based on directions, the grid can be constructed over the unit hypersphere rather than the unit ball, with an empirical distribution converging weakly to the uniform over the  unit hypersphere. This is obtained by letting $n_S=n$, $n_R=1$, and $n_0=0$.} is that the sequence  $\mathfrak{G}\n$ of grids is such that the discrete distribution with probability masses~$1/n$ at each gridpoint and probability mass~$n_0/n$ at the origin converges weakly to the uniform~${\rm U}_d$ over the ball $\mathbb{S}_d$. 
%\textcolor{red}{Once we define  $\F^{(n)}_\pm$, we set $\Q^{(n)}_\pm := (\F_{\pm}^{(n)})^{-1}$ which is the sample counterpart of $\Q_\pm$ (again, as in the case of the population quantities, we do a small abuse of notation and set~$\Q^{(n)}_\pm ({\bf 0}) := (\F_{\pm}^{(n)})^{-1}(\{{\bf 0}\})$).} 
Then, we  define~$\F_{\pm}^{(n)}(\ZZ^{(n)}_t) $, for~$t = 1,\ldots ,n$ as the solution (optimal mapping) of a coupling problem between the residuals %observations 
and the grid. Specifically, the empirical center-outward distribution function is the (random) discrete mapping\vspace{-2mm} 
$$
\F^{(n)}_{\pm}: \ZZ^{(n)}:= (\ZZ^{(n)}_1,\ldots ,\ZZ^{(n)}_n) \mapsto (\F^{(n)}_{\pm}(\ZZ^{(n)}_1),\ldots ,\F^{(n)}_{\pm}(\ZZ^{(n)}_n))
$$
satisfying 
%\begin{equation}\label{Fpm0}
%\sum_{t=1}^{n}  \Vert \ZZ^{(n)}_t - (\F^{(n)}_{\pm}(\ZZ^{(n)}_t)\Vert^2 = \min_{\pi}  \sum_{t=1}^{n}  \Vert \ZZ^{(n)}_{\pi t} - (\F^{(n)}_{\pm}(\ZZ^{(n)}_t)\Vert^2,
%\end{equation}
%(where $\ZZ^{(n)}_t=\ZZ^{(n)}_t(\bth)$ and~$\Vert  \cdot \Vert $ stands for the Euclidean norm) or, equivalently, 
\begin{equation}\label{Fpm0}
\s \Vert \ZZ^{(n)}_t - \F^{(n)}_\pm(\ZZ^{(n)}_t)\Vert ^2 = \underset{T \in \mathcal{T}}{\min} \s \Vert \ZZ^{(n)}_t - T(\ZZ^{(n)}_t)\Vert ^2,
\end{equation}
where $\ZZ^{(n)}_t=\ZZ^{(n)}_t(\bth)$, the set $\{\F^{(n)}_{\pm}(\ZZ^{(n)}_t)\vert t = 1,\ldots ,n\}$ coincides with the $n$ points of the grid, 
$\Vert  \cdot \Vert $ stands for the Euclidean norm,
%where  
%$\pi$ in (\ref{Fpm0})   ranges over the $n!$ possible permutations of $\{1,\ldots ,n\}$, and 
and $\mathcal{T}$ denotes the set of all possible bijective mappings between $\ZZ^{(n)}_1, \ldots , \ZZ^{(n)}_n$ and the $n$ gridpoints.

Based on this empirical center-outward distribution function,  %$\F^{(n)}_\pm$, for every $1 \leq t \leq n$,
  the {\it center-outward ranks} are defined as 
\begin{equation}
{R^{(n)}_{\pm, t}}:=  {R^{(n)}_{\pm, t}}(\bth):=({n_R + 1}) \Vert \F^{(n)}_\pm(\ZZ^{(n)}_t) \Vert,  \label{Ranks}
\end{equation}
 the {\it center-outward signs}  as 
\begin{equation}
\S^{(n)}_{\pm, t}:=  {\S^{(n)}_{\pm, t}}(\bth):=  \F^{(n)}_\pm(\ZZ^{(n)}_t)  I [\F^{(n)}_\pm(\ZZ^{(n)}_t) \neq \boldsymbol{0}]/ 
\Vert \F^{(n)}_\pm(\ZZ^{(n)}_t) \Vert. \label{Signs}
\end{equation}
%\color{red} and the \textit{center-outward quantile contours} as
%\begin{equation}
%\mathcal{C}^{(n)}_{\pm;\ZZ^{(n)}_t}\left(\frac{j}{n_R +1} \right):= \left\{  \ZZ^{(n)}_t \vert R^{(n)}_{\pm, t} = j  \right\},
%\label{Cont}
%\end{equation}
%where ${j}/{(n_R +1)}$, $j = 0, 1, . . . , n_R$ is an empirical probability content, to be interpreted as a quantile order. 
%We refer to \cite{Hallin2017} for details.  
It follows that $\F^{(n)}_\pm(\ZZ_t^{(n)})$ factorizes into\vspace{-2mm}
\begin{equation}\label{factFpm}
\F^{(n)}_\pm(\ZZ_t^{(n)})=   \frac{R^{(n)}_{\pm, t}}{n_R + 1}  \S^{(n)}_{\pm, t},\quad\text{ whence }\quad \ZZ_t^{(n)}= \Q^{(n)}_\pm\Big(
 \frac{R^{(n)}_{\pm, t}}{n_R + 1}  \S^{(n)}_{\pm, t}\Big).
\end{equation}

%In dimension $d=1$, 
Those ranks and signs are jointly distribution-free (for  $f\in{\cal F}_d$): more precisely,  under~${\rm P}^{(n)}_{\bth_0;f}$,  the $n$-tuple $\F^{(n)}_\pm(\ZZ_1^{(n)}),\ldots , \F^{(n)}_\pm(\ZZ_n^{(n)})$ is uniformly distributed over the $n!/n_0!$ permutations with repetition of the~$n$ underlying  gridpoints (the origin having multiplicity~$n_0$). 
%For $n_0=0$, the vector of ranks  ${R^{(n)}_{\pm, 1}},\ldots, {R^{(n)}_{\pm, n}}$  is uniformly distributed over the $n_R!$ permutations of the integers~$(1,\ldots, n_R)$, the $n_S$-tuple of signs $\S^{(n)}_{\pm, 1},\ldots,\S^{(n)}_{\pm, n}$ is uniformly distributed.
Moreover, the center-outward distribution functions, ranks, and signs inherit, from the invariance properties of Euclidean distances, elementary but remarkable invariance and equivariance properties:   \cite{HHH20} %Hallin et al. (2020)
 show that center-outward quantities enjoy invariance/equivariance with respect to shift, global scale, and orthogonal transformations.
 
 \section{A  H\' ajek asymptotic representation result for serial center-outward rank statistics }\label{SecAsRepr}
 
\subsection{H\' ajek asymptotic representation}\label{hajeksec}
Throughout this section, denote by $\F_{\pm, t}:= \F_{\pm}(\ZZ_t^{(n)}(\bth_0))$  the  value of the center-outward distribution function $\F_\pm$ associated with innovation density $f$  computed at   $\ZZ_t=\ZZ_t^{(n)}(\bth_0)$ and \linebreak by~$\F_{\pm, t}\n:= \F_{\pm}\n(\ZZ_t^{(n)}(\bth_0))$  its empirical counterpart computed from the $n$-tuple of resi\-duals~$\ZZ_1^{(n)}(\bth_0),\ldots,\ZZ_n^{(n)}(\bth_0)$.  
Considering  two {\it score functions} $\J_1$ and $\J_2$  from ${\mathbb S}_d$ to $\mathbb{R}^d$, define, for $1 \leq i \leq n - 1,$\vspace{-2mm}
 \begin{equation}\label{tildeGam}
%\tenq{\bGamma}_{i, \J_1, \J_2}^{(n)}(\bth_0) := 
\tenq{\bGamma}_{i, \J_1, \J_2}^{(n)}(\bth_0) := (n-i)^{-1} \sum_{t=i+1}^n \a (\F^{(n)}_{\pm, t}, \F^{(n)}_{\pm, t-i})
\vspace{-2mm}\end{equation}
and 
 \begin{equation}\label{barGam}
%\underline{\bGamma}_{i, \J_1, \J_2}^{(n)}(\bth_0) := 
\underline{\bGamma}_{i, \J_1, \J_2}^{(n)}(\bth_0) := (n-i)^{-1} \sum_{t=i+1}^n \J_1(\F_{\pm, t}) \J_2\pr(\F_{\pm, t-i}),
\vspace{-1mm}\end{equation}
where the function $\a : \mathbb{S}_d \times \mathbb{S}_d \rightarrow \mathbb{R}^{d^2}$ satisfies, under ${\rm P}^{(n)}_{\bth_0 ;f}$,\vspace{-2mm}
\begin{equation}\label{ass.a}
\lim_{n\rightarrow \infty}{\rm E} \left\Vert \text{vec} \left( \a (\F^{(n)}_{\pm, 2}, \F^{(n)}_{\pm, 1}) - \J_1(\F_{\pm, 2}) \J_2\pr(\F_{\pm, 1}) \right) \right\Vert^2 \rightarrow 0. \vspace{-3mm}
\end{equation}
Let \vspace{-2mm} 
\begin{align*}
\m^{(n)} :=& {\rm E}_{{\bth_0};f} [\a (\F^{(n)}_{\pm, 2}, \F^{(n)}_{\pm, 1})]  \\
=& \,[n(n-1)]^{-1}\Bigg[\sum_{\mathfrak{g}_1\neq\mathfrak{g}_2\in\mathfrak{G}\n\setminus{\bf 0} }\!\!\!\!\a (\mathfrak{g}_2,\mathfrak{g}_1) +n_0\!\!\!\!\sum_{\mathfrak{g}\in\mathfrak{G}\n\setminus{\bf 0} }\!\!\!\!\a (\mathfrak{g},{\bf 0}) 
 + n_0\!\!\!\!\sum_{\mathfrak{g}\in\mathfrak{G}\n\setminus{\bf 0} }\!\!\!\!\a ({\bf 0},\mathfrak{g}) + 
n_0(n_0-1)\a ({\bf 0},{\bf 0})
\Bigg] 
 \end{align*}\vspace{-6mm}
 
\noindent  and\vspace{-3mm} 
\begin{align*}
 \m :=& {\rm E}_{{\bth_0};f} [\J_1(\F_{\pm, 2}) \J_2\pr(\F_{\pm, 1})] 
 =  {\rm E}_{{\bth_0};f} [\J_1(\F_{\pm, 1})] {\rm E}_{{\bth_0};f} [\J_2\pr(\F_{\pm, 1})], 
 \end{align*}
where ${\rm E}_{{\bth_0};f}$ denotes expectation under ${\rm P}^{(n)}_{\bth_0 ;f}$. 
 Clearly, $\m^{(n)}\!$, which does not depend on~$f$, is centering $\tenq{\bGamma}_{i, \J_1, \J_2}^{(n)}$ under the (null) hypothesis that~$(\F\n_{\pm, 1},\ldots, \F\n_{\pm, n})$ is uniform over the permutations of the grid it was constructed from (that is, under any innovation density in~$\mathcal{F}_d$%{\rm P}^{(n)}_{\bth_0 ;g}$
 ), while $\m =\m_f$, which depends on $f$ through $\F_\pm$, is cen\-tering~$\underline{\bGamma}_{i, \J_1, \J_2}^{(n)}$ under ${\rm P}^{(n)}_{\bth_0 ;f}$ only; none of them 
%neither $ \m$ nor $ \m\n$
 depends on $i$. 

An essential step in H\' ajek's approach to the asymptotics of univariate ranks is the so-called H\' ajek asymptotic representation of linear rank statistics  which, contrary to earlier approaches based on empirical processes, allows for unbounded square-integrable score functions. Such results first were obtained (H\' ajek~1961; H\' ajek and \v Sid\'ak 1967) for the nonserial rank statistics used in the context of  linear models (single-output regression, ANOVA, etc.);  Hallin et al.~(1985)  extend them to serial rank statistics---with scores involving several ranks at a time---appearing in the analysis of time-series models. The matrices $ \tenq{\bGamma}_{i, \J_1, \J_2}^{(n)}$ are of that serial type, hence require an extension of Hallin et al.~(1985) and \cite{HV96} (which only deal with classical ``univariate ranks") to the present case of multivariate center-outward ranks and signs. 
Proposition~\ref{Prop.bar.til.Gam} establishes such asymptotic representation for 
%result  the asymptotic equivalence of
 $ \text{\rm vec} ( \tenq{\bGamma}_{i, \J_1, \J_2}^{(n)}(\bth_0) - \m^{(n)})$ under 
 % and $(n-i)^{1/2} \text{\rm vec}(\underline{\bGamma}_{i, \J_1, \J_2}^{(n)}(\bth_0) - \m)$.
%We make
 the following assumption on the score functions $\J_1$ and~$\J_2$. 

\setcounter{assumption}{1}

\begin{assumption}\label{ass.vp}
{\it (i)} {\rm $\J_1$ and $\J_2$ are continuous over $\mathbb{S}_d$;}\\
{\it (ii)} {\rm  $\J_1$ and $\J_2$ are  square-integrable, that is, 
$%\[
 \int_{\mathbb{S}_d} \Vert  \J_\ell (\u) \Vert^2 {\rm d U}_d < \infty$ for $ \ell =1, 2$, 
%$ %\]
 and, for any sequence $\mathfrak{s}^{(n)} := \{{\bf s}_1^{(n)}, \ldots, {\bf s}_n^{(n)}\}$ of $n$-tuples in  $\mathbb{S}_d$ such that the uniform discrete distribution over $\mathfrak{s}^{(n)}$ converges weakly to ${\rm U}_d$ as $n \rightarrow \infty$,}
\begin{equation}\label{ass2}
\underset{n\rightarrow\infty}{\lim} n^{-1} \sum_{t=1}^n 
 \Vert \J_\ell ({\bf s}_t^{(n)}) \Vert^2 =    \int_{\mathbb{S}_d} \Vert  \J_\ell (\u) \Vert^2 {\rm d U}_d  , \quad \ell =1, 2.
\end{equation}
\end{assumption}

%{\it (iii)} $\int_{\mathbb{S}_d}  \J_\ell (\u)  {\rm d} U_d = \0, \quad \ell =1, 2$.
%\end{assumption}
%
%
%\textbf{Remark}. Assumption~\ref{ass.vp}(iii) is quite mild as  it is satisfied by the optimal scores $\J_1 = \bvp_{f} \circ \Q_{\pm}$ and $\J_2 = \Q_{\pm}$ for any $f \in \mathcal{F}_d$. Specifically, for any $f$ with mean zero, we have $\int_{\mathbb{S}_d}  \Q_{\pm} (\u)  {\rm d} U_d = \0$.
%Moreover, for $\J_1 = \bvp_{f} \circ \Q_{\pm}$, we have 
%$\int_{\mathbb{S}_d}  \bvp_{f} \circ \Q_{\pm} (\u)  {\rm d} U_d = \int_{\mathbb{R}^d} \bvp_{f}(\z) f(\z) {\rm d} \z = \0$.

When $\J _\ell ({\u})$ is of the form ({\it spherical} score) $J_\ell(\Vert \u\Vert)\frac{\u}{\Vert \u\Vert}$ with $J_\ell:[0,1) \to\mathbb{R}$ continuous, a sufficient condition for \eqref{ass2} is the traditional assumption of bounded variation (that is,~$J_\ell$ is the continuous difference of two nondecreasing functions).

\begin{proposition}\label{Prop.bar.til.Gam}
Let Assumption~2 % \ref{ass.vp}
 hold.  Then, for any positive integer $i$, %\color{blue} 
 %(should we use ${\rm E}$ instead of ${\rm E}_{{\bth_0};f}$ here?) \color{black}
\begin{equation*}
 {\rm E}\left\Vert   \text{\rm vec} \left( \tenq{\bGamma}_{i, \J_1, \J_2}^{(n)}(\bth_0) - \m^{(n)}   - \underline{\bGamma}_{i, \J_1, \J_2}^{(n)}(\bth_0) + \m \right) \right\Vert^2 = o((n-i)^{-1/2}) 
\end{equation*}
under ${\rm P}^{(n)}_{\bth_0 ;f}$ (any $f\in {\cal F}_d$) and ${\rm P}^{(n)}_{{\bth_0} + n^{-1/2} \btau ;f}$ (any $f\in {\cal F}^*_d$),  as $n \rightarrow \infty$.
\end{proposition}

See Section~\ref{Sec3} for the proof. \medskip

Adopting H\' ajek's terminology, examples of functions $\a$ satisfying \eqref{ass.a} for given $\J_1$ and~$\J_2$ are the \textit{approximate scores} \vspace{-1mm}
$$\a_{\rm a} (\F^{(n)}_{\pm, t}, \F^{(n)}_{\pm, s}) := \J_1(\F^{(n)}_{\pm, t}) \J_2\pr(\F^{(n)}_{\pm, s})\vspace{-1mm}$$
and (still with   ${\rm E}_{{\bth_0};f}$ denoting expectation under ${\rm P}^{(n)}_{\bth_0 ;f}$) the \textit{exact scores}\vspace{-1mm}
$$\a_{\rm e} (\F^{(n)}_{\pm, t}, \F^{(n)}_{\pm, s}) :=  {\rm E}_{{\bth_0};f}\left[\J_1(\F_{\pm, t}) \J_2\pr (\F_{\pm, s}) \big| \F^{(n)}_{\pm, t}, \F^{(n)}_{\pm, s}\right].\vspace{-1mm}$$
% The following Lemma~\ref{lem.xi} %and \ref{lem.xi2}
%  states that both $\a_{\rm a}$ and $\a_{\rm e}$ satisfy \eqref{ass.a}.
We then have, for $\a_{\rm a}$ and $\a_{\rm e}$,  the following lemma.

\begin{lemma}\label{lem.xi}
Let Assumption~2  % \ref{ass.vp}
 hold. Put\vspace{-2mm}
    \[
  \bxi_{t, s; {\rm a}}^{(n)}\! := \text{\rm vec} (\a_{\rm a} (\F^{(n)}_{\pm, t}, \F^{(n)}_{\pm, s})), \ \ \bxi_{t, s; {\rm e}}^{(n)}\! := \text{\rm vec} (\a_{\rm e} (\F^{(n)}_{\pm, t}, \F^{(n)}_{\pm, s})),\vspace{-2mm}\]
and
     \[ \bxi_{t, s}\! := \text{\rm vec} (\J_1(\F_{\pm, t}) \J_2\pr(\F_{\pm, s})).\]
Then, under ${\rm P}^{(n)}_{\bth_0 ;f}$ (any $f\in{\cal F}_d$) and ${\rm P}^{(n)}_{\bth_0 + n^{-1/2} \btau ;f}$ (any $f\in{\cal F}_d^*$), 
  \begin{equation}\label{lem21}
{\rm E} \Vert \bxi_{2, 1; {\rm a}}^{(n)} - \bxi_{2, 1} \Vert^2 \rightarrow 0\quad\text{and}\quad 
{\rm E} \Vert \bxi_{2, 1; {\rm e}}^{(n)} - \bxi_{2, 1} \Vert^2 \rightarrow 0
\quad\text{as $n \rightarrow \infty$. }\end{equation}

\end{lemma}

Letting\vspace{-2mm} 
\[
\m_{\rm a}^{(n)} := {\rm E}_{{\bth_0};f}  [\a_{\rm a} (\F^{(n)}_{\pm, t}, \F^{(n)}_{\pm, s})],\quad\m_{\rm e}^{(n)} := {\rm E}_{{\bth_0};f}  [\a_{\rm e} (\F^{(n)}_{\pm, t}, \F^{(n)}_{\pm, s})]
\vspace{-1mm}\]
and, for $1 \leq i \leq n - 1$,\vspace{-2mm}
$$
\tenq{\bGamma}_{i; {\rm a}}^{(n)}(\bth_0)  := (n-i)^{-1} \sum_{t=i+1}^n \a_{\rm a} (\F^{(n)}_{\pm, t}, \F^{(n)}_{\pm, t-i})
\vspace{-2mm}$$
and 
\begin{equation*} 
\tenq{\bGamma}_{i; {\rm e}}^{(n)}(\bth_0) := (n-i)^{-1} \sum_{t=i+1}^n \a_{\rm e} (\F^{(n)}_{\pm, t}, \F^{(n)}_{\pm, t-i}), 
\end{equation*}
the following proposition follows as a corollary to Proposition~\ref{Prop.bar.til.Gam} and Lemma~\ref{lem.xi}.

\begin{proposition}\label{Prop.HajekRep} 
Let Assumption 2 hold.  Then, for any positive integer $i$, \vspace{-1mm}
\begin{equation*} 
{(i)} \  {\rm E}\left\Vert (n-i)^{1/2}  \text{\rm vec} ( \tenq{\bGamma}_{i; {\rm a}}^{(n)}(\bth_0) - \m_{\rm a}^{(n)}   - \underline{\bGamma}_{i, \J_1, \J_2}^{(n)}(\bth_0) + \m) \right\Vert ^2 = o(1),
\end{equation*}
\begin{equation*}
{(ii)} \  {\rm E}\left\Vert(n-i)^{1/2}  \text{\rm vec} ( \tenq{\bGamma}_{i; {\rm e}}^{(n)}(\bth_0) - \m_{\rm e}^{(n)}   - \underline{\bGamma}_{i, \J_1, \J_2}^{(n)}(\bth_0) + \m) \right\Vert^2 = o(1),
\vspace{-1mm}\end{equation*}
and, consequently,\vspace{-2mm}
\begin{equation*}
{(iii)} \  {\rm E}\left\Vert (n-i)^{1/2}  \text{\rm vec} ( \tenq{\bGamma}_{i; {\rm a}}^{(n)}(\bth_0) - \m_{\rm a}^{(n)}    - \tenq{\bGamma}_{i; {\rm e}}^{(n)}(\bth_0) + \m_{\rm e}^{(n)})  \right\Vert^2 = o(1)
\vspace{-1mm}\end{equation*}
under ${\rm P}^{(n)}_{\bth_0 ;f}$ (any $f\in{\cal F}_d$) and ${\rm P}^{(n)}_{\bth_0 + n^{-1/2} \btau ;f}$ (any $f\in{\cal F}_d^*$), as $n \rightarrow \infty$.
\end{proposition}

Assume, without loss of generality, that $\J_1$ and $\J_2$ are such that $\m ={\bf 0}$  (a sufficient condition is either $\int_{\mathbb{S}_d} \J_1(\u) {\rm d U}_d ={\bf 0}$ or $\int_{\mathbb{S}_d} \J_2(\u) {\rm d U}_d ={\bf 0}$). Then, for $d=~\!1$,~$\m_{\rm a}^{(n)}$ and~$\m_{\rm e}^{(n)}$ both are 
%$o_{\rm P}(n^{-1/2})$ \color{blue} (maybe
 $o(n^{-1/2})$  %?) \color{black}:
  (see Lemma~1 in \cite{HL2017}) and %As a consequence,
    $\m$, $\m_{\rm a}^{(n)}$, and $\m_{\rm e}^{(n)}$ can be omitted in Proposition~\ref{Prop.HajekRep}. That simplification, however, is no longer valid when~$d>1$.

The proofs  of Proposition~\ref{Prop.bar.til.Gam} and Lemma~\ref{lem.xi} are given in Appendix~B.

%\section{Center-outward rank-based tests}

%\end{document}

\section{Center-outward rank-based central sequences} \label{SecSR}
The basic idea underlying the construction of   our rank-based tests is the definition of  rank-based versions $\utilde\bDelta^{(n)}_{f} (\bth_0)$ of the central sequences~$\bDelta^{(n)}_{f} (\bth_0)$. Recall that $\utilde\bDelta^{(n)}_{f} (\bth_0)$ qualifies as a central sequence as soon as $\utilde\bDelta^{(n)}_{f} (\bth_0)-\bDelta^{(n)}_{f} (\bth_0)=o_{\rm P}(1) $ under ${\rm P}^{(n)}_{\bth_0 ;f}$.

%\color{red}STOPPED HERE September 25\color{black}
\subsection{Construction and asymptotic representation}\label{Sec41}

Writing $\F^{(n)}_{\pm, t}$, $R^{(n)}_{\pm, t}$ and $\S^{(n)}_{\pm, t}$   for~$\F^{(n)}_\pm(\ZZ_t^{(n)}(\bth_0))$, $ {R^{(n)}_{\pm, t}}(\bth_0)$ and~$ {\S^{(n)}_{\pm, t}}(\bth_0)$, respectively, consider the center-outward rank-based counterpart of ${\bGamma}^{(n)}_f (\bth_0)$. Specifically, define \vspace{-2mm}  
 \begin{align}
\tenq{\bGamma}_{\J_1 \J_2}^{(n)}(\bth_0) := &   n^{-1/2} \left( (n-1)^{1/2} (\text{vec}\,  \tenq{\bGamma}_{1, \J_1, \J_2}^{(n)}(\bth_0))^\prime, \ldots , \right. \nonumber \\
&\qquad \left.  (n-i)^{1/2} (\text{vec}\,  \tenq{\bGamma}_{i, \J_1, \J_2}^{(n)}(\bth_0))^\prime, \ldots ,
(\text{vec}\, \tenq{\bGamma}_{n-1, \J_1, \J_2}^{(n)}(\bth_0))^\prime\right)^\prime, \label{tilde.S.m}
\end{align}
\vspace{-10mm}

\noindent with  \vspace{-2mm} 
 \begin{equation}\label{tildeGam1}
\tenq{\bGamma}_{i, \J_1, \J_2}^{(n)}(\bth_0) := (n-i)^{-1} \sum_{t=i+1}^n \J_1\left(\frac{R^{(n)}_{\pm, t}}{n_R + 1}\S^{(n)}_{\pm, t}\right) \J_2\pr\left(\frac{R^{(n)}_{\pm, t-i}}{n_R + 1}\S^{(n)}_{\pm, t-i}\right)   , \quad i = 1, \ldots , n - 1,
\end{equation}
where  $\J_1$ and~$ \J_2: {\mathcal S}_d \rightarrow \R$
 are score functions satisfying Assumption~\ref{ass.vp}. Call $\tenq{\bGamma}_{i, \J_1, \J_2}^{(n)}(\bth_0) $  a (residual)  {\it lag-$i$ rank-based~cross-covariance matrix}.   As an approximate-score  rank-based statistic, $\tenq{\bGamma}_{i, \J_1, \J_2}^{(n)}(\bth_0)$  under ${\rm P}^{(n)}_{\bth_0 ;f}$   has expectation\vspace{-2mm}
 $$\m_{\rm a}^{(n)} = {\rm E}_{{\bth_0};f} \left[\tenq{\bGamma}_{i, \J_1, \J_2}^{(n)}(\bth_0)\right] = {\rm E}_{{\bth_0};f}  \left[\J_1(\F^{(n)}_{\pm, 2}) \J_2\pr(\F^{(n)}_{\pm, 1})\right]$$
 under ${\rm P}^{(n)}_{\bth_0 ;f}$, hence under any ${\rm P}^{(n)}_{\bth_0 ;g}$, $g\in{\cal F}_d$.  
 We thus have\vspace{-2mm} 
 \begin{equation}\label{EGam}
 {\rm E}_{{\bth_0};f}  \left[\tenq{\bGamma}_{\J_1 \J_2}^{(n)}(\bth_0)\right] =  n^{-1/2}   \left( {\rm diag} ((n-1)^{1/2}, \ldots , (n-i)^{1/2}, \ldots , 1) \otimes \I_{d^2} \right) \, \text{vec}\, (\m_{\rm a}^{(n)}).
 \end{equation}
 Similarly define\vspace{-2mm}
  \begin{align}
\underline{\bGamma}_{\J_1 \J_2}^{(n)}(\bth_0) := &   n^{-1/2} \left( (n-1)^{1/2} (\text{vec}\,  \underline{\bGamma}_{1, \J_1, \J_2}^{(n)}(\bth_0))^\prime, \ldots , \right. \nonumber \\
&\qquad \left.  (n-i)^{1/2} (\text{vec}\,  \underline{\bGamma}_{i, \J_1, \J_2}^{(n)}(\bth_0))^\prime, \ldots ,
(\text{vec}\, \underline{\bGamma}_{n-1, \J_1, \J_2}^{(n)}(\bth_0))^\prime\right)^\prime, \label{underline.S.m}
\end{align}
\vspace{-8mm}

\noindent with $\underline{\bGamma}_{i, \J_1, \J_2}^{(n)}(\bth_0)$   defined in \eqref{barGam}. 
Clearly, $\underline{\bGamma}_{\J_1 \J_2}^{(n)}(\bth_0)$ is a counterpart of $\tenq{\bGamma}_{\J_1 \J_2}^{(n)}(\bth_0)$ constructed from the center-outward distribution function and under ${\rm P}^{(n)}_{\bth_0 ;f}$  has mean\vspace{-2mm}
$${\rm E}_{{\bth_0};f}  \left[\underline{\bGamma}_{\J_1 \J_2}^{(n)}(\bth_0)\right] =  n^{-1/2}   \left( {\rm diag} ((n-1)^{1/2}, \ldots , (n-i)^{1/2}, \ldots , 1) \otimes \I_{d^2} \right) \, \text{vec}\, (\m)$$

 Denote by $\utilde\bDelta^{(n)}_{\J_1 \J_2} (\bth_0)$ and $\underline{\bDelta}^{(n)}_{\J_1 \J_2} (\bth_0)$   the    statistics resulting from substituting\vspace{-2mm} 
$$\tenq{\bGamma}_{\J_1 \J_2}^{(n)}(\bth_0)-  {\rm E}_{{\bth_0};f}  \left[\tenq{\bGamma}_{\J_1 \J_2}^{(n)}(\bth_0)\right]
\quad\text{ and }\quad   \underline{\bGamma}_{\J_1 \J_2}^{(n)}(\bth_0) -  {\rm E}_{{\bth_0};f}  \left[\underline{\bGamma}_{\J_1 \J_2}^{(n)}(\bth_0)\right],\vspace{-2mm}$$
respectively,   for ${\bGamma}_{f}^{(n)}(\bth_0)$  in the definition~\eqref{Delta} of the central sequence~$\bDelta^{(n)}_{f} (\bth_0)$.  For the scores~$\J_1=\bvp\circ {\F}_\pm^{-1}$ and $\J_2={\F}_\pm^{-1}$, the following corollary to Proposition
\ref{Prop.HajekRep} {\it (i)} implies that~$\utilde\bDelta^{(n)}_{\J_1 \J_2} (\bth_0)$ is asymptotically equivalent, under ${\rm P}^{(n)}_{\bth_0 ;f}$, to $\bDelta^{(n)}_{f} (\bth_0)$, hence constitutes a rank-based version of the same central sequence. These scores depend on   $f$: if $f$ stands for the actual innovation density, thus, they are {\it oracle scores}. They can be computed, however, at any chosen {\it reference density} such as the spherical or skew Gaussian; see Section~\ref{Sec: examples}.  With a slight abuse of terminology,  irrespective of $\J_1$ and $\J_2$, we therefore call $\utilde\bDelta^{(n)}_{\J_1 \J_2}(\bth_0)$ a {\it rank-based central sequence}. 
Proposition~\ref{Prop.HajekRep} {\it (iii)} moreover implies that, asymptoti\-cally,~$\utilde\bDelta^{(n)}_{\J_1 \J_2} (\bth_0)$    is equivalent to the expectation of $\bDelta^{(n)}_{f} (\bth_0)$ conditional on the $\sigma$-field generated by the  center-outward ranks and signs---that is, the projection of $\bDelta^{(n)}_{f} (\bth_0)$ onto the space of center-outward ranks and signs. We thus can state the following corollary.

%{\it center-outward rank-based} central sequence  %(see \cite{HLL2019} for motivation of the construction of $\utilde\bDelta^{(n)}_{J_1, J_2} (\bth_0)$).   

%\color{red} Isn't an essential piece missing here? Before calling $\utilde\bDelta^{(n)}_{J_1, J_2} (\bth_0)$ the {\it center-outward rank-based} central sequence, we should show that indeed, for adequate choices of $J_1$ and $J_2$, it si a version of the central sequence $\bDelta^{(n)}_{f} (\bth_0)$! Therefore we need an asymptotic representation result of the form $\utilde\bDelta^{(n)}_{J_1, J_2} (\bth_0)= \bDelta^{(n)}_{J_1, J_2} (\bth_0) + o_{\rm P}(1)$, then the scores $J_1, J_2$ such that~$\bDelta^{(n)}_{J_1, J_2} (\bth_0)= \bDelta^{(n)}_{f} (\bth_0) $. \color{black}

\begin{corollary}\label{lem.bar.til.Delta.n}
Let Assumption~\ref{ass.vp} hold. Then, \vspace{-2mm}
\begin{equation}\label{bar.til.Delta.n}
{\rm E} \left\Vert \utilde\bDelta^{(n)}_{\J_1 \J_2}(\bth_0) - %\bar
\underline{\bDelta}^{(n)}_{\J_1 \J_2}(\bth_0) \right\Vert^2 = o(1)
\end{equation}
under  ${\rm P}^{(n)}_{\bth_0 ;f}$ (any~$f\in{\cal F}_d$) and~${\rm P}^{(n)}_{\bth_0 + n^{-1/2} \btau ;f}$ (any~$f\in{\cal F}^*_d$), as $n \rightarrow \infty$.
\end{corollary}

\subsection{Asymptotic normality %of %$\utilde{\bDelta}^{(n)}_{\J_1 \J_2}(\bth_0)$
}\label{Sec42}

The  joint  asymptotic normality of the statistics $(n-i)^{1/2}\text{\rm vec} (\tenq{\bGamma}_{i, \J_1, \J_2}^{(n)}(\bth_0) - \m_{\rm a}^{(n)})$ and their asymptotic  linearity, which entail the asymptotic normality and asymptotic  linearity of~$\utilde{\bDelta}^{(n)}_{\J_1 \J_2}(\bth_0)$,  are essential in the construction of our rank-based tests. In view of the H\' ajek asymptotic representation result of Section~\ref{SecAsRepr}, it is sufficient to derive these asymptotic results for $(n-i)^{1/2}\text{\rm vec} (\underline{\bGamma}_{i, \J_1, \J_2}^{(n)}(\bth_0) - \m)$. 

 Let \vspace{-3mm} 
\begin{align*}
{\C}_{\J_1 \J_2}:= &  \int_{\mathbb{S}_d} (\I_d \otimes \J_1(\u_1)) \left[ \int_{\mathbb{S}_d} \J_2(\u_2) \J_2^\prime(\u_2) {\rm dU}_d(\u_2) \right] (\I_d \otimes \J_1^\prime(\u_1)) {\rm dU}_d(\u_1)  \\
&- \Big( \I_d \otimes  \int_{\mathbb{S}_d} \J_1(\u) {\rm dU}_d \Big)   \int_{\mathbb{S}_d} \J_2(\u) {\rm dU}_d   \int_{\mathbb{S}_d} \J_2\pr(\u) {\rm dU}_d \Big( \I_d \otimes  \int_{\mathbb{S}_d} \J_1\pr(\u) {\rm dU}_d \Big) 
 \end{align*}
 \vspace{-8mm}
 
\noindent and\vspace{-2mm}
\begin{align*}
\K_{\J_1, \J_2, f} &:=   \int_{\mathbb{S}_d} (\I_d \otimes \J_1(\u_1)) \left[ \int_{\mathbb{S}_d} \J_2(\u_2) \F_{\pms}^{-1\prime}(\u_2) {\rm dU}_d(\u_2) \right] (\I_d \otimes \bvp^\prime_{f}(\F_{\pms}^{-1}(\u_1)) {\rm dU}_d(\u_1):
\end{align*}
 under Assumptions~\ref{ass.den} and \ref{ass.vp},  $\K_{\J_1, \J_2, f}$ exists and is finite by the Cauchy--Schwarz inequality. Decompose the matrix $\Q^{(n)}_{\bth}$ defined in \eqref{defPQ}  into $d^2 \times d^2 p_1$   blocks\vspace{-2mm}  
$$\Q^{(n)}_{\bth_0} =\big(\Q_{1, \bth_0}^\prime, \ldots , \Q_{i, \bth_0}^\prime, \ldots, \Q_{n-1, \bth_0}^\prime \big)^\prime\vspace{-2mm}$$
 (note that these blocks depend on $i$, not on $n$).  Lemma~\ref{asy.Gami.bar} states the asymptotic normality of $(n-i)^{1/2}\text{\rm vec} (\underline{\bGamma}_{i, {\J}_1, {\J}_2}^{(n)}(\bth_0) - \m)$; see Appendix~C for the proof.

 \begin{lemma}\label{asy.Gami.bar}
Let Assumption~\ref{ass.vp} hold. Then, 
 for any positive integers $i\neq j$, the vectors
 $$(n-i)^{1/2}\text{\rm vec} (\underline{\bGamma}_{i, {\J}_1, {\J}_2}^{(n)}(\bth_0) - \m)\quad\text{and}\quad (n-j)^{1/2}\text{\rm vec} (\underline{\bGamma}_{j, {\J}_1, {\J}_2}^{(n)}(\bth_0) - \m)$$
  are jointly asymptotically normal, 
with mean $(\0\pr, \0\pr)\pr$ under  ${\rm P}^{(n)}_{\bth_0 ;f}$ (any~$f\in{\cal F}_d$),  mean\vspace{-2mm}
$$\left((\K_{\J_1, \J_2, f} \Q_{i, \bth_0} \P_{\bth_0} \M_{\bth_0} \btau)\pr, ((\K_{\J_1, \J_2, f} \Q_{j, \bth_0} \P_{\bth_0} \M_{\bth_0} \btau)\pr\right)\pr
$$
under ${\rm P}^{(n)}_{\bth_0 + n^{-1/2}\btau ;f}$ (any~$f\in{\cal F}^*_d$), and   covariance  $\left(
\begin{array}{cc} {\C}_{\J_1 \J_2}& \0 \\ \0& {\C}_{\J_1 \J_2} \end{array}
\right)$  under both.
\end{lemma}

%\color{blue} Note that, by assuming ${\rm E}_{{\bth_0};f} (\ZZ_t(\bth_0)) = \0$, we have $\m = {\rm E}_{{\bth_0};f}  [\underline{\bGamma}_{\J_1 \J_2}^{(n)}(\bth_0)] = \0$ \color{olive} is that correct?needed?\color{black}. 
%Denote by $\underline{\bDelta}^{(n)}_{\J_1 \J_2} (\bth_0)$ the  asymptotic representation of $\utilde{\bDelta}^{(n)}_{\J_1 \J_2} (\bth_0)$ under ${\rm P}^{(n)}_{\bth_0 ;f}$: clearly, $\underline{\bDelta}^{(n)}_{\J_1 \J_2} (\bth_0)$ is obtained by 
%
%
%sequence resulting from substituting $\underline{\bGamma}_{\J_1 \J_2}^{(n)}(\bth_0)$  for $\tenq{\bGamma}_{\J_1 \J_2}^{(n)}(\bth_0)- {\rm E}_{{\bth_0};f} [\tenq{\bGamma}_{\J_1 \J_2}^{(n)}(\bth_0)]$  in $\utilde{\bDelta}^{(n)}_{\J_1 \J_2} (\bth_0)$. Lemma~\ref{lem.bar.til.Delta.n} extends the H\' ajek asymptotic representation result in Proposition~\ref{Prop.HajekRep} to 
% $\utilde{\bDelta}^{(n)}_{\J_1 \J_2} (\bth_0)$.
%
%
%\begin{proof}
%The result follows by using Proposition~\ref{Prop.HajekRep} and by moving along the same lines as the proof of Lemma~B.4 in \cite{HLL2019}.
%\end{proof}
Piecing together Lemma~\ref{asy.Gami.bar}, the definition of ${\bDelta}^{(n)}_{\J_1 \J_2} (\bth_0)$, and Corollary~\ref{lem.bar.til.Delta.n}, we then readily obtain the asymptotic normality of  $\utilde{\bDelta}^{(n)}_{\J_1 \J_2} (\bth_0)$ under ${\rm P}^{(n)}_{\bth_0 ;f}$ (any~$f\in{\cal F}_d$) and contiguous alternatives (any~$f\in{\cal F}^*_d$). Let  $\T^{(n)}_{\bth_0} := \M_{\bth_0}^\prime  \P_{\bth_0}^\prime \Q_{\bth_0}^{(n)\prime},$\vspace{-2mm} 
\begin{equation}\label{defLam}
\bUpsilon^{(n)}_{\J_1, \J_2, f}(\bth_0) :=   \T^{(n)}_{\bth_0} (\I_{n-1} \otimes \K_{\J_1, \J_2, f}) \T^{(n)\prime}_{\bth_0}, \quad   \bLam^{(n)}_{\bth_0} :=  \T^{(n)}_{\bth_0} (\I_{n-1} \otimes \C_{\J_1, \J_2}) \T^{(n)\prime}_{\bth_0},
\end{equation}
$$\bUpsilon_{\J_1, \J_2, f}(\bth_0) := \underset{n\rightarrow \infty}{\lim} \bUpsilon^{(n)}_{\J_1, \J_2, f}(\bth_0), \quad \text{and} \quad  \bLam_{\bth_0} :=   \underset{n\rightarrow \infty}{\lim} \bLam^{(n)}_{\bth_0}.$$
where the limits exist due to the exponential decrease of the Green matrices; see Appendix~A for details.  Note that when~$\J_1=\bvp\circ {\F}_\pm^{-1}$ and $\J_2={\F}_\pm^{-1}$, ${\C}_{\J_1 \J_2}$ and $\K_{\J_1 \J_2, f}$ coincide, so that~$\bUpsilon_{\J_1, \J_2, f}(\bth_0)= \bLam_{\bth_0}$  is the information matrix $\bLam_{f} (\bth_0)$ in Proposition~\ref{Prop.LAN1}.

%Then the asymptotic normality of $\utilde\bDelta^{(n)}_{\J_1 \J_2}(\bth_0)$ follows directly from Lemma~\ref{asy.Gami.bar} and \eqref{lem.bar.til.Delta.n}.

\begin{proposition}\label{til.Delta.n}
Let Assumption~\ref{ass.vp} hold. Then 
$\utilde\bDelta^{(n)}_{\J_1 \J_2}(\bth_0)$ is asymptotically normal 
with mean $\0$ under  ${\rm P}^{(n)}_{\bth_0 ;f}$ (any~$f\in{\cal F}_d$),  mean $\bUpsilon_{\J_1, \J_2, f}(\bth_0) \btau$
under ${\rm P}^{(n)}_{\bth_0 + n^{-1/2}\btau ;f}$ (any~$f\in{\cal F}^*_d$), and  covariance~$\bLam_{\bth_0}$  
 under both.\vspace{-1mm}
\end{proposition}

Finally, in order to construct our test statistics, we still need the asymptotic linearity of $\utilde\bDelta^{(n)}_{\J_1 \J_2}(\bth_0)$; the latter is an immediate consequence of the asymptotic linearity for all~$i$ of~$(n-i)^{1/2}\text{\rm vec} (\tenq{\bGamma}_{i, \J_1, \J_2}^{(n)}(\bth_0))$, which we now assume; the form of the linear term in the right-hand side of \eqref{assA4} follows from the form of the asymptotic shift in Lemma~\ref{asy.Gami.bar}. \medskip\vspace{-2mm}

\begin{assumption}\label{asylin1}
For any positive integer $i$, as $n\to\infty$, \vspace{-2mm}
\begin{equation}\label{assA4}
(n-i)^{1/2}\left[\text{\rm vec} (\tenq{\bGamma}^{(n)}_{i, \J_1, \J_2}(\bth_0 + n^{-1/2}\btau)) - \text{\rm vec} (\tenq{\bGamma}^{(n)}_{i, \J_1, \J_2}(\bth_0)) \right] = - \K_{\J_1, \J_2, f} \Q_{i, \bth_0} \P_{\bth_0} \M_{\bth_0} \btau + o_{\rm P}(1)
\end{equation}
\vspace{-8mm}

\noindent under  ${\rm P}^{(n)}_{\bth_0 ;f}$ (any~$f\in{\cal F}^*_d$) (hence also under ${\rm P}^{(n)}_{\bth_0 + n^{-1/2}\btau ;f}$).
\end{assumption}

The asymptotic linearity of 
%$(n-i)^{1/2}\text{\rm vec} (\tenq{\bGamma}_{i, \J_1, \J_2}^{(n)}(\bth_0)$ entails that of
  $\utilde{\bDelta}^{(n)}_{\J_1 \J_2}(\bth_0))$ readily follows.

\begin{proposition}\label{asy0}
Let Assumptions \ref{ass.den},  \ref{ass.vp}, and \ref{asylin1} hold. Then,  as  $n \to \infty $, \vspace{-1mm} 
\begin{equation}\label{asy.linear2}
\utilde{\bDelta}^{(n)}_{\J_1 \J_2}(\bth_0 + n^{-1/2}\btau) - \utilde{\bDelta}^{(n)}_{\J_1 \J_2}(\bth_0) =  -  \bUpsilon^{(n)}_{\J_1, \J_2, f}(\bth_0)  \btau + o_{\rm P}(1)
\end{equation}
\vspace{-9mm}

\noindent under  ${\rm P}^{(n)}_{\bth_0 ;f}$ (any~$f\in{\cal F}^*_d$),  hence also under ${\rm P}^{(n)}_{\bth_0 + n^{-1/2}\btau ;f}$.\vspace{-1mm} 
\end{proposition}

\section{Center-outward rank-based tests}\label{Sec:Rtests}

%The null hypothesis we consider imposes some linear constraints on $\bth$, that is $\bth \in \bth_0 + \mathcal{M}(\bvartheta)$, where $\mathcal{M}(\bvartheta)$ denotes the vector subspace of $\mathbb{R}^{d^2p_1}$ spanned by the columns of some $d^2 p_1 \times d^2 p_0$ matrix $\bvartheta$. Here we discuss two special cases separately.

In this section, based on central sequences of the form   $\utilde\bDelta^{(n)}_{\J_1 \J_2}(\bth_0)$, we propose center-outward rank-based tests  for VAR models and, depending on the scores $\J_1$ and $\J_2$,   derive their asymptotic properties. 
 We separately consider two cases:  \vspace{-2mm}
\begin{enumerate}
\item[(a)] testing  the null hypothesis $\bth = \bth_0$, with specified $\bth_0$ (VAR order $p_0$), against the alternative of the form $\bth \neq \bth_0$ with VAR order $p$ ($p$ potentially larger than $p_0$); since $\bth_0$ is specified, rank-based statistics are fully distribution-free under the null (any~$f\in{\cal F}_d$);  \vspace{-7mm}
\item[(b)] testing  the null hypothesis of a VAR of order $p_0$ (with unspecified parameter $\bth_0$) against the alternative of a VAR of order $p > p_0$ (with $p\leq p_1$ and unspecified parameter $\bth$). For $p_0>0$, the null value of $\bth_0$ will be estimated;   rank-based statistics then are only {\it strongly} asymptotically distribution-free, that is, asymptotically equivalent, under the null (any~$f\in{\cal F}_d^*$),  to a strictly distribution-free statistic.\vspace{-1mm}
\end{enumerate}
Case (a) is essentially of theoretical interest, and a preparation for case~(b). The latter has an obvious application in  the sequential  identification  of the order of VAR models via  a sequence of tests of the null hypotheses of VAR($p_0$) against VAR($p_0 +1$) dependence,~$p_0=0, 1, \ldots$
%. More precisely, the procedure starts from $p_0 = 0$, where we test white noise against VAR($1$), i.e., existence of autocorrelation with order one. When the test statistic is significant, we move to test of VAR($1$) against VAR($2$), and so on.  

\subsection{Testing $\bth = \bth_0$ ($\bth_0$ specified)}\label{Sec.theta0}

We first consider the case of  the null hypothesis of $\bth = \bth_0$, with $\bth_0$ being specified, against the alternative of VAR($p$), $ p_0\leq p \leq p_1$ with parameter $\bth \neq \bth_0$.   More precisely, the null hypothesis is $\mathcal{H}^{(n)}_{\bth_0}:=\{{\rm P}\n_{{\bth_0};f}\vert f\in{\mathcal F}_d\}$. 
The center-outward rank-based test statistics we are proposing are of form\vspace{-2mm}\footnote{See \cite{HP04} for a similar test statistic based on Mahalanobis ranks and signs.}
\begin{equation}\label{Stilde}
\utilde{S}_{\J_1 \J_2}^{(n)}(\bth_0) := \utilde{\H}_{\J_1 \J_2}^{(n)\prime}(\bth_0) \left(  \Q^{(n)\prime}_{\bth_0} (\I_{n-1} \otimes \C_{\J_1, \J_2})  \Q^{(n)}_{\bth_0} \right)^{-1} \utilde{\H}_{\J_1 \J_2}^{(n)}(\bth_0),
\end{equation}
\vspace{-8mm}

\noindent where 
%we set  $p = p_1$ \color{blue} (I think I made a mistake here, and we should delete $p = p_1$) \color{red} but what should we write here, then? \color{black} in the definition~\eqref{defPQ} of~$\Q^{(n)}_{\bth_0}$ and 
\begin{equation}\label{Htilde}
\utilde{\H}_{\J_1 \J_2}^{(n)}(\bth_0) :=  n^{1/2}  \Q^{(n)\prime}_{\bth_0} \left(\utilde{\bGamma}_{\J_1 \J_2}^{(n)}(\bth_0) - {\rm E}_{{\bth_0};f}  [\tenq{\bGamma}_{\J_1 \J_2}^{(n)}(\bth_0)]\right)
\end{equation}
with ${\rm E}_{{\bth_0};f}  [\tenq{\bGamma}_{\J_1 \J_2}^{(n)}(\bth_0)]$ (which does not depend on $f$)  given in \eqref{EGam} and $\Q^{(n)}_{\bth_0}$ in~\eqref{defPQ}.

The following proposition establishes the asymptotic distribution of $\utilde{S}^{(n)}_{\J_1\J_2}(\bth_0)$ under~$\mathcal{H}^{(n)}_{\bth_0}$ and % $\mathcal{H}^{(n)}_{\bth_0 + n^{-1/2} \btau}$
contiguous alternatives of the form ${\rm P}^{(n)}_{\bth_0 + n^{-1/2} \btau; f}$ (any $f\in{\cal F}_d^*$).  Also, in order to discuss local optimality properties, we consider scores $\J_1 = \J_{1;f}  := \bvp_{f} \circ \Q_{\pm}$ and $\J_2=\J_{2;f}  :=  \Q_{\pm}$ and denote by~$\utilde{S}^{(n)}_f(\bth_0)$ ($f\in{\cal F}_d^*$) the corresponding test statistic. See Appendix~C for~a~proof.

\begin{proposition}\label{Prop.Stilde}
Let Assumption~\ref{ass.vp} hold. Then,
\begin{enumerate}[label=(\roman*)]
\item under $\mathcal{H}^{(n)}_{\bth_0}$, $\utilde{S}_{\J_1 \J_2}^{(n)}(\bth_0)$  is fully distribution-free and asymptotically chi-square with $d^2 p_1$    degrees of freedom;
\item under  ${\rm P}^{(n)}_{\bth_0 + n^{-1/2} \btau; f}$ (any~$f\in{\cal F}_d^*$), 
%$\mathcal{H}^{(n)}_{\bth_0 + n^{-1/2} \btau}$
 $\utilde{S}_{\J_1 \J_2}^{(n)}(\bth_0)$ is asymptotically non-central chi-square with $d^2 p_1$   degrees of freedom and non-centrality parameter  
\begin{align*}
&\btau^\prime   \underset{n\rightarrow \infty}{\lim} \left[ \T^{(n)}_{\bth_0} (\I_{n-1} \otimes \K_{\J_1, \J_2, f}) \Q^{(n)}_{\bth_0}  \left(  \Q^{(n)\prime}_{\bth_0} (\I_{n-1} \otimes \C_{\J_1, \J_2}) \Q^{(n)}_{\bth_0} \right)^{-1}\!\! \right.\\
& \qquad\qquad \qquad\qquad \qquad \qquad\qquad \times \left. \Q^{(n)\prime}_{\bth_0} (\I_{n-1} \otimes \K_{\J_1, \J_2, f}^{\prime}) \T^{(n)\prime}_{\bth_0} \right]  \btau; 
\end{align*}
\vspace{-10mm}
\item the sequence of tests rejecting $\mathcal{H}^{(n)}_{\bth_0}$ whenever   $\utilde{S}^{(n)}_f(\bth_0)$ (any~$f\in{\cal F}_d^*$) exceeds the~$(1-\alpha)$ chi-square quantile with $d^2 p_1$   degrees of freedom  is locally asymptotically maximin for~$\mathcal{H}^{(n)}_{\bth_0}$ versus~$\bigcup_{\bth \neq \bth_0} {\rm P}^{(n)}_{{\bth};f}$ at asymptotic  level $\alpha$. 
\end{enumerate}
\end{proposition}

\subsection{Testing the order of a VAR (unspecified $\bth_0$)}\label{Sec.testorder}

% 
%To formalize this, we write $\mathcal{H}^{(n)}(\bth_0;f)$ for the hypothesis under which an observation $\X^{(n)}$ is generated by 
%the VAR($p_0$) model with ${\rm P}^{(n)}_{\bth_0;f}$, namely by a model with parameter value $\bth_0$ satisfying Assumption~\ref{assVAR} and innovation process having density 
%$f\in \mathcal{F}_d$ satisfying Assumption~\ref{ass.den}. Our inference problem is to
%test 
%$$\mathcal{H}^{(n)}(\bth_0):= \cup_{f} \mathcal{H}^{(n)}(\bth_0;f) \quad \text{vs} \quad \cup_{\bth \neq \bth_0} \mathcal{H}^{(n)}(\bth_0).$$ 
%Therefore,  in our problem the innovation density $f$ is a nuisance parameter and the union in the definition
%of $\mathcal{H}^{(n)}_{\bth_0}$ is taken over all possible $f\in \mathcal{F}_d$.

Most null hypotheses  of practical interest involve  incompletely  specified values of the parameter, though. In the problem of testing a VAR of order $p_0$ against a VAR of order $p$ with~$0<p_0<p\leq p_1$, for instance, the value $\bth_0$ of the parameter under the null remains unspecified. A natural idea then consists in replacing the unspecified $\bth_0$ with a root-$n$ consistent (constrained) estimator~$\hat{\bth}\n$ in  \eqref{Stilde}, yielding the {\it aligned rank} test statistic $\utilde{S}_{\J_1 \J_2}^{(n)}(\hat{\bth}\n)$. Such plug-in, however, has an impact on the asymptotic distribution of $\utilde{S}_{\J_1 \J_2}^{(n)}$; that impact is neutralized via the following classical
%\footnote{See \citet[Chapter~??]{LeCam86}, where the same construction is applied to a general central sequence~$\bDelta\n (\bth)$.} 
construction  inspired from the traditional theory of Neyman  $C(\alpha)$ tests; see Section~11.9 of~\cite{LeCam86} for details.

%For a real data, $\bth_0$ is usually unknown. Therefore, an estimator of $\bth_0$ is needed before applying a test statistic. Moreover, asymptotically, the plug-in of this estimator should have no impact on the test statistic. Since $\utilde{S}^{(n)}(\bth_0)$ in Section~\ref{Sec.theta0} does not meet this requirement, we consider constructing our test statistic as follows.

Considering the null hypothesis $\mathcal{H}^{(n)}_{0}:=\bigcup_{f\in{\cal F}_d^*}\left\{ {\rm P}\n_{{\bth_0};f} \big\vert \bth_0 \text{ of the form } \eqref{bth0}\right\}$ of a VAR of order~$p_0$,  split the $d^2p_1$-dimensional\footnote{The rank-based  central sequence associated with  the reduced-rank  VAR($p_1$) process characterized by $\bth_0$ of the form \eqref{bth0}} $\utilde{\bDelta}^{(n)}_{\J_1 \J_2}(\bth_0)$ into $\left(\utilde{\bDelta}^{(n)\prime}_{I; \J_1, \J_2}(\bth_0),\ \utilde{\bDelta}^{(n)\prime}_{II; \J_1, \J_2}(\bth_0)
\right)\pr$ where~$\bth_0$ is an arbitrary parameter value of the form \eqref{bth0}, and 
$\utilde{\bDelta}^{(n)}_{I; \J_1, \J_2}(\bth_0)$ and~$\utilde{\bDelta}^{(n)}_{II; \J_1, \J_2}(\bth_0)$ 
 are~$d^2 p_0$-  and $d^2 (p_1 - p_0)$-dimensional, respectively. Conformably partition the $d^2p_1\times d^2p_1$ matrix $\bUpsilon^{(n)}_{\J_1, \J_2, f}(\bth_0)$ defined in \eqref{defLam} into\vspace{-2mm}
\begin{equation}\label{Psin}
\bUpsilon^{(n)}_{\J_1, \J_2, f}(\bth_0)  =
\begin{bmatrix}
\bUpsilon^{(n)}_{11; \J_1, \J_2, f}(\bth_0) & \bUpsilon^{(n)}_{12; \J_1, \J_2, f}(\bth_0) \\
\bUpsilon^{(n)}_{21; \J_1, \J_2, f}(\bth_0)  & \bUpsilon^{(n)}_{22; \J_1, \J_2, f}(\bth_0) 
\end{bmatrix}
,
\end{equation}
with square diagonal blocks $\bUpsilon^{(n)}_{11; \J_1, \J_2, f}(\bth_0)$ and $\bUpsilon^{(n)}_{22; \J_1, \J_2, f}(\bth_0)$  of size $d^2 p_0$ and $d^2 (p_1~\!-~\!p_0)$, respectively. 
%, for~$\bUpsilon^{(n)}_{\J_1, \J_2, f}(\bth_0)$, $\bUpsilon_{\J_1, \J_2, f}(\bth_0)$,  $\bUpsilon^{(n)}_{ij; \J_1, \J_2, f}(\bth_0)$ and $\bUpsilon_{ij; \J_1, \J_2, f}(\bth_0)$, respectively.
 Similarly partition $\bUpsilon_{\J_1, \J_2, f}(\bth_0)$ and, for the sake of simplicity,  use the notation~$\bUpsilon^{(n)}$, $\bUpsilon^{(n)}_{ij}$, $\bUpsilon$, and~$\bUpsilon_{ij}, \, i, j = 1, 2$  in an obvious way. 

Next,   consider the residual $\utilde{\bDelta}^{(n)*}_{II; \J_1, \J_2}(\bth_0)$ of the regression of~$\utilde{\bDelta}^{(n)}_{II; \J_1, \J_2}(\bth_0)$ on~$\utilde{\bDelta}^{(n)}_{I; \J_1, \J_2}(\bth_0)$ in the shift matrix $\bUpsilon^{(n)}$, namely,\vspace{-2mm}  
 \begin{equation}\label{defDeltastar}
\utilde{\bDelta}^{(n)*}_{II; \J_1, \J_2}(\bth_0) := \utilde{\bDelta}^{(n)}_{II; \J_1, \J_2}(\bth_0) - {\bUpsilon}^{(n)}_{21} ({\bUpsilon}^{(n)}_{11})^{-1} \utilde{\bDelta}^{(n)}_{I; \J_1, \J_2}(\bth_0).
\end{equation}
Partitioning $\bLam^{(n)}_{\bth_0}$ in \eqref{defLam} into the same  block form\vspace{-2mm}
 \[
\bLam^{(n)}_{\bth_0}  =
\begin{bmatrix}
\bLam^{(n)}_{11; \bth_0} & \bLam^{(n)}_{12; \bth_0} \\
\bLam^{(n)}_{21; \bth_0}  & \bLam^{(n)}_{22; \bth_0}
\end{bmatrix}
=:
\begin{bmatrix}
\bLam^{(n)}_{11} & \bLam^{(n)}_{12} \\
\bLam^{(n)}_{21}  & \bLam^{(n)}_{22}
\end{bmatrix}
\]
\vspace{-8mm}

\noindent as $\bUpsilon^{(n)}%=\bUpsilon^{(n)}_{\J_1, \J_2, f}(\bth_0)
$, let\vspace{-5mm} 
\begin{align}
\bLam^{(n)*}_{II; \bth_0} 
& := \bLam^{(n)}_{22} + {\bUpsilon}^{(n)}_{21} ({\bUpsilon}^{(n)}_{11})^{-1} \bLam_{11} ({\bUpsilon}^{(n)\prime}_{11})^{-1} {\bUpsilon}^{(n)\prime}_{21}  \nonumber \\
&\qquad - \bLam^{(n)}_{21} ({\bUpsilon}^{(n)\prime}_{11})^{-1} {\bUpsilon}^{(n)\prime}_{21} 
- {\bUpsilon}^{(n)}_{21} ({\bUpsilon}^{(n)}_{11})^{-1} \bLam^{(n)}_{12} \label{defLamstar}
\end{align}
\vspace{-10mm}

\noindent and define 
\begin{equation}\label{defW}
\utilde{W}^{(n)}_{\J_1 \J_2}(\bth_0) := \utilde{\bDelta}^{(n)*^\prime}_{II; \J_1, \J_2}(\bth_0)  (\bLam^{(n)*}_{II; \bth_0})^{-1} \utilde{\bDelta}^{(n)*}_{II; \J_1, \J_2}(\bth_0).
\end{equation}

 It is easy to see that $\utilde{\bDelta}^{(n)*}_{II; \J_1, \J_2}(\bth_0)$ is the projection of $\utilde{\bDelta}^{(n)}_{II; \J_1, \J_2}(\bth_0)$ onto the space which, in the metric characterized by \eqref{Psin}, is  orthogonal to $\utilde{\bDelta}^{(n)}_{I; \J_1, \J_2}(\bth_0)$ hence, by Le Cam's third lemma, is insensitive to local perturbations of $\bth_0$; the matrix   $\bLam^{(n)*}_{II; \bth_0} $   then is~$\utilde{\bDelta}^{(n)*}_{II; \J_1, \J_2}(\bth_0)$'s (residual) covariance. Unfortunately, $\utilde{W}^{(n)}_{\J_1 \J_2}(\bth_0)$, because it depends on the unspecified $\bth_0$ and involves $\bUpsilon_{\J_1, \J_2, f}(\bth_0)\vspace{1mm}$, cannot be used as a test statistic.

Therefore, let 
 $\hat{\bth}^{(n)} := ((\text{vec}{\hat{\A}^{(n)}_1})^\prime, \ldots , (\text{vec}{\hat{\A}^{(n)}_{p_0}})^\prime, \0_{d^2(p_1-p_0)\times 1}^\prime)^\prime$ 
 denote a  constrained---that is, of the form \eqref{bth0}---root-$n$ consistent (under~${\rm P}^{(n)}_{\bth_0; f}$)\footnote{In view of contiguity, this also holds under ${\rm P}^{(n)}_{\bth_0 + n^{-1/2} \btau; f}$ with $\btau$ of the form \eqref{bth0}.} estimator of~$\bth_0$. For technical reasons,   assume that  $\hat{\bth}^{(n)}$ moreover is {\it asymptotically discrete}.\footnote{An estimator $\hat{\bth}^{(n)}$ of $\bth_0$ is called {\it asymptotically discrete} if, for $n$ large enough,  it takes   at most $K(c)$ distinct values in balls of the form 
 $%\mathcal{Q}_n :=
  \{\bth \in \mathbb{R}^{p_1 d^2}: n^{-1/2} \left\Vert \bth - \bth_0 \right\Vert \leq c\}$, $c > 0$ fixed, $K(c)$ independent of $n$. Any root-$n$ consistent estimator is easily discretized into an asymptotically discrete one. Asymptotic discreteness, however, is only  a theoretical requirement since, in practice, $\hat{\bth}^{(n)}$ anyway only has a finite number of digits. While discretization is necessary in asymptotic statements, it is not needed  in finite-$n$ implementation of tests; see \citet[Chapter 6]{LeCam2000} and \citet[Section 5.7]{Vaart1998} for   details.
}
  %Now we show the impact from substituting $\hat{\bth}^{(n)}$ for $\bth_0$ in $\utilde{W}^{(n)}_{\J_1 \J_2}$ is asymptotically negligible.
%
%
%
%More precisely,  
%see \citet[Section 4]{Kreiss87}  for details.
%
%
%
%
%\begin{assumption}\label{assRootn}
%Under ${\rm P}^{(n)}_{\bth_0; f}$, as $n \rightarrow \infty$,
%$\hat{\bth}^{(n)} - \bth_0 = O_{\rm p}(n^{-1/2}).$
%\end{assumption}
%
%Note that LAN implies 
%that ${\rm P}^{(n)}_{\bth_0 + n^{-1/2} \btau ;f}$ and  ${\rm P}^{(n)}_{\bth_0 ;f}$ are mutually contiguous; Assumption~\ref{assRootn} hence also holds under ${\rm P}^{(n)}_{\bth_0 + n^{-1/2} \btau; f}$.
Also,  denote by $\hat{\bUpsilon}^{(n)}$ a consistent (under~${\rm P}^{(n)}_{\bth_0; f}$) estimator of $\bUpsilon_{\J_1, \J_2, f}(\bth_0)$. Such an estimator can be obtained, e.g.,  from the asymptotic linearity in \eqref{asy.linear2}: see Section~\ref{sec.comp} for details. 
 Our test statistic is obtained by plugging these estimators into \eqref{defW}, yielding\vspace{-2mm} 
\begin{equation}\label{testeq}
\widehat{\utilde{W}}^{(n)}_{\J_1 \J_2}(\hat{\bth}^{(n)}) :=  \utilde{\bDelta}^{(n)*^\prime}_{II; \J_1, \J_2}(\bth_0)  (\bLam^{(n)*}_{II; \bth_0})^{-1} \utilde{\bDelta}^{(n)*}_{II; \J_1, \J_2}(\bth_0).\vspace{-3mm}
\end{equation}

For scores $\J_1 = \J_{1;f}  := \bvp_{f} \circ \Q_{\pm}$ and $\J_2=\J_{2;f}  :=  \Q_{\pm}$ write~$\widehat{\utilde{W}}^{(n)}_{f}$ instead of $\widehat{\utilde{W}}^{(n)}_{\J_1 \J_2}(\hat{\bth}^{(n)})$. We then have the following asymptotic results. % of $\widehat{\utilde{W}}^{(n)}_{\J_1 \J_2}(\hat{\bth}^{(n)})$.

\begin{proposition}\label{Prop.Wtilde}
Let Assumptions %\ref{ass.den}, 
\ref{ass.vp} and \ref{asylin1} hold. Then,\vspace{-2mm}
\begin{enumerate}[label=(\roman*)]
\item   under any ${\rm P}^{(n)}_{\bth_0}$ in  $\mathcal{H}^{(n)}_{0}$  and contiguous alternatives of the form 
$\mathcal{H}^{(n)}_{\bth_0+ n^{-1/2} \btau; f}$  ($f\in{\cal F}_d^*$), \vspace{-2mm}
$$\widehat{\utilde{W}}^{(n)}_{\J_1 \J_2}(\hat{\bth}^{(n)}) - \utilde{W}^{(n)}_{\J_1 \J_2}(\bth_0) = o_{\rm P}(1)\quad \text{ as } n\to\infty;$$
\item under any ${\rm P}^{(n)}_{{\bth_0};f}$  (with~$f\in{\cal F}_d^*$) in  $\mathcal{H}^{(n)}_{0}$, %${\rm P}^{(n)}_{\bth_0; f}$, 
$\widehat{\utilde{W}}^{(n)}_{\J_1 \J_2}(\hat{\bth}^{(n)})$ is asymptotically chi-square with $d^2 (p_1 - p_0)$ degrees of freedom;
\item under %${\rm P}^{(n)}_{\bth_0 + n^{-1/2} \btau; f}$, 
${\rm P}^{(n)}_{\bth_0+ n^{-1/2} \btau; f}$ with  $f\in{\cal F}_d^*$, $\btau\pr=(\btau_{I}\pr ,\btau_{II}\pr)$, $\btau_{I}\in\mathbb{R}^{d^2p_0}$ and $\btau_{II}\in\mathbb{R}^{d^2(p_1-p_0)}$, 
$\widehat{\utilde{W}}^{(n)}_{\J_1 \J_2}(\hat{\bth}^{(n)})$ is asymptotically non-central chi-square with $d^2 (p_1 - p_0)$ degrees of freedom and non-centrality parameter\vspace{-2mm} 
\begin{align*}
& \btau_{II}^\prime 
 \left( \bUpsilon\pr_{22,\J_1, \J_2, f}(\bth_0)  %\right. 
% \left. 
 - \bUpsilon\pr_{12,\J_1, \J_2, f}(\bth_0)(\bUpsilon\pr_{11,\J_1, \J_2, f}(\bth_0))^{-1} \bUpsilon\pr_{21,\J_1, \J_2, f}(\bth_0) \right) 
  \\ 
&\quad 
\times\! (\bLam^{*}_{II; \bth_0})^{-1}\! \left( \bUpsilon_{22,\J_1, \J_2, f}(\bth_0) - \bUpsilon_{21,\J_1, \J_2, f}(\bth_0) \left(\bUpsilon_{11,\J_1, \J_2, f}(\bth_0)\right)^{-1}\! \bUpsilon_{12,\J_1, \J_2, f}(\bth_0)  \right) \btau_{II},
\end{align*}
\vspace{-10mm}

\noindent where $\bLam^{*}_{II; \bth_0} := \underset{n \rightarrow \infty}{\lim} \bLam^{(n)*}_{II; \bth_0}$;\vspace{-3mm}
\item the sequence of tests  rejecting $\mathcal{H}^{(n)}_{0}$ whenever $\widehat{\utilde{W}}^{(n)}_{f}$ exceeds   the chi-square  $(1-\alpha)$ quantile with $d^2 (p_1 - p_0)$ degrees of freedom  is locally asymptotically  most stringent for $\mathcal{H}^{(n)}_{0}$ versus $\bigcup_{\bth \neq \bth_0}{\rm P}^{(n)}_{\bth ;f}$  ($f\in{\cal F}_d^*$)   at asymptotic  level $\alpha$.
\end{enumerate}
\end{proposition}

Propositions \ref{Prop.Stilde} {\it (i)} and  \ref{Prop.Wtilde} {\it (ii)} provide the asymptotic null distributions of the test statistics $\utilde{S}
^{(n)}_{\J_1 \J_2}$ and $\widehat{\utilde{W}}
^{(n)}_{\J_1 \J_2}$ %under the null hypothesis $\mathcal{H}^{(n)}_{\bth_0}$
 that can be used to construct uniformly valid (in the sense of \eqref{UpseudoG}) asymptotic critical values.  
%. The implication of the result is that test procedures based on $\utilde{S}
%^{(n)}_{\J_1 \J_2}$ and $\widehat{\utilde{W}}
%^{(n)}_{\J_1 \J_2}$ remain valid (have the right asymptotic size) for any distribution $f$  satisfying Assumption~\ref{ass.den}; see \cite{Hallinetal2020}
%for a discussion on the connection of this result with the distribution-freeness of the center-outward ranks and signs. We flag that, 
%differently from the test procedures based on Mahalanobis signs and ranks (\cite{HP04, HP05}), our test 
%statistics can be applied for VAR adequacy testing problems under a broad class of underlying innovation densities, which includes both elliptical  and non-elliptical densities, and thus it is larger than the one considered by \cite{HP04, HP05}.
%
Propositions \ref{Prop.Stilde} {\it (ii)} and \ref{Prop.Wtilde} {\it (iii)} provides  asymptotic distributions under contiguous alternatives, allowing for local power evaluations.   Propositions~\ref{Prop.Stilde}~{\it (iii)} and~\ref{Prop.Wtilde} {\it (iv)} establish the optimality properties of the same tests.  

For $p_0=0$, that is, when   testing white noise against   VAR($p_1$) dependence,  
$\utilde{\bDelta}^{(n)}_{\J_1 \J_2}(\bth_0)$ reduces to~$n^{1/2}\{\tenq{\bGamma}_{\J_1, \J_2}^{(p_1, n)}(\bth_0) - {\rm E}_{{\bth_0};f} (\tenq{\bGamma}_{\J_1, \J_2}^{(p_1, n)}(\bth_0))\}$ with $\bth_0 = \0$ and\vspace{-2.5mm} 
$$\tenq{\bGamma}_{\J_1, \J_2}^{(p_1, n)}(\bth_0) := n^{-1/2} \left( (n-1)^{1/2} (\text{vec}\,  \tenq{\bGamma}_{1, \J_1, \J_2}^{(n)}(\bth_0))^\prime, \ldots ,    (n-p_1)^{1/2} (\text{vec}\,  \tenq{\bGamma}_{p_1, \J_1, \J_2}^{(n)}(\bth_0))^\prime \right)^\prime\vspace{-2mm}$$
being the $p_1$-block  truncated version of $\tenq{\bGamma}_{\J_1, \J_2}^{(n)}(\bth_0)$. In this case, the test statistics $\utilde{S}
^{(n)}_{\J_1 \J_2}(\0)$ and $\widehat{\utilde{W}}
^{(n)}_{\J_1 \J_2}(\0)$ coincide and both take the  form\vspace{-2mm}
\begin{align}
%\utilde{W}^{(n)}_{\J_1 \J_2} &:=
 n \left\{ \tenq{\bGamma}_{\J_1, \J_2}^{(p_1, n)}(\0)  - {\rm E}_{{\bth_0};f} (\tenq{\bGamma}_{\J_1, \J_2}^{(p_1, n)}(\0)) \right\}^\prime (\I_{p_1} \otimes \C_{\J_1, \J_2})^{-1}
\left\{ \tenq{\bGamma}_{\J_1, \J_2}^{(p_1, n)}(\0) - {\rm E}_{{\bth_0};f} (\tenq{\bGamma}_{\J_1, \J_2}^{(p_1, n)}(\0)) \right\}. 
\label{tildeW.white}
\end{align}
\vspace{-10mm}

\noindent If $p_1=1$,   Proposition~\ref{Prop.Wtilde} implies that the null hypothesis of white noise is rejected in favor of VAR($1$) dependence whenever $\utilde{S}
^{(n)}_{\J_1 \J_2}\!\!=\widehat{\utilde{W}}
^{(n)}_{\J_1 \J_2}$ exceeds the $(1-\alpha)$ quantile $\chi^2_{d^2;1-\alpha}$ of the chi-square distribution with $d^2$  degrees of freedom. \vspace{-2mm}

\subsection{Some standard score functions}\label{Sec: examples}

The rank-based cross-covariance matrices \vspace{-1mm}
$\tenq{\bGamma}^{(n)}_{\J_1 \J_2}$, hence also the  test statistics~$\utilde{S}_{\J_1 \J_2}^{(n)}(\bth_0)$ and $\widehat{\utilde{W}}^{(n)}_{\J_1 \J_2}(\hat{\bth}^{(n)})$,  depend on the choice of score functions $\J_1$ and $\J_2$, to be selected by the practitioner. In principle, in view of the maximin and stringency properties in Propositions~\ref{Prop.Stilde}~{\it (iii)} and  \ref{Prop.Wtilde} {\it (iv)}, an optimal choice is $\J_1 = \bvp_{f} \circ \Q_{\pm}$ and $\J_2 =  \Q_{\pm}$, where $f$ is the actual innovation density and $\Q_{\pm}$ the corresponding 
 center-outward quantile function. Such a choice, unfortunately, is unfeasible   since $f$ is unspecified. Alternatives can privilege either simplicity, robustness, or efficiency at some chosen reference density. Here are three examples 
% 
% 
% 
%one should choose . To implement these optimal scores, however, one needs an estimator of $\bvp_{f}$ and an explicitly form or estimator of $\Q_{\pm}$ for a specified reference density. The estimation of $\bvp_{f}$ involves estimation of some parameters if the density is pre-specified to belong to a particular parametric distribution family (e.g. skew-$t$, mixtures of Gaussians or elliptical distributions). This procedure is feasible (see \cite{Hallinetal2020} for more details), but computationally challenging and costly. For $\Q_{\pm}$, except for a few cases of  spherical distributions, no explicit form is available in the literature. In principle, $\Q_{\pm}$ can be estimated arbitrarily precise via simulation once the reference density is fully specified. This procedure, however, also involves large computational cost\footnote{Fortunately, increasingly efficient algorithms are available in the very active domain of numerical measure transportation (see, e.g., \cite{Merigot} and \cite{PC2019})}.
%
%
%Considering the challenges of choosing optimal scores for any $f \in \mathcal{F}_d$, one may choose scores that are easy-to-implement as long as Assumption~\ref{ass.vp} is satisfied. Here we provide three examples
 of sensible choices   extending scores that are widely applied in the univariate  (see e.g. \cite{HL2017}) and  the elliptical multivariate setting  (see \cite{HP04}). Note that for these scores, the centering term $\m_{\rm a}\n$ for~$\tenq{\bGamma}_{i, \J_1, \J_2}^{(n)}(\bth_0)$ disappears  as soon as the grid ${\mathfrak G}\n$ is symmetric with respect to the origin while $\m = \0$ since~$\J_{\ell}(-{\bf u}) = -\J_{\ell}({\bf u})$,  $\ell = 1, 2$.\medskip

\textbf{Example 1} ({\it Sign test} scores). Setting %\vspace{-2mm} 
$\J_{\ell}\left(\frac{R^{(n)}_{\pm, t}}{n_R + 1}\S^{(n)}_{\pm, t}\right)  = \S^{(n)}_{\pm, t}, \quad \ell = 1, 2$  
yields the center-outward sign-based  cross-covariance matrices\vspace{-2mm} %as follows
\begin{align}
\tenq{\bGamma}_{i, \text{sign}}^{(n)}(\bth) = (n-i)^{-1} \sum_{t=i+1}^n   \S^{(n)}_{{\pms}, t}(\bth) \S^{(n)\prime}_{{\pms}, t-i}(\bth), \quad i = 1, \ldots , n-1. \label{Eq. Signs}
\end{align}
\vspace{-8mm}

\noindent The resulting   $\utilde{\bDelta}^{(n)}_{\text{sign}}(\bth)$  entirely relies  on the center-outward  signs $\S^{(n)}_{{\pms}, t}(\bth)$, which thus explains the terminology  {\it sign test}  scores.   Since the ranks do not enter the score function, there is no need for  the grid ${\mathfrak G}\n$ to fill  the unit ball as $n\to\infty$, and we can set $n_S=n$, producing a grid over the unit hypersphere rather than the unit ball. \vspace{2mm}% (see Section~\ref{Sec: NumStud}). \color{black}
 %\\

\textbf{Example 2} ({\it Spearman} scores).   Another simple choice is $\J_1(\u)  = \u =  \J_2(\u)$.  
The corresponding rank-based cross-covariance matrices   are\vspace{-2mm}
\begin{align}
\tenq{\bGamma}_{i, \text{Sp}}^{(n)}(\bth) = (n-i)^{-1} \sum_{t=i+1}^n   \F^{(n)}_{{\pms}, t} \F^{(n)\prime}_{{\pms}, t-i}, \quad i = 1, \ldots , n-1
\end{align}
\vspace{-8mm}

\noindent reducing, for $d=1$,   to Spearman autocorrelations, whence the terminology  {\it Spearman}~scores. \vspace{1mm}%\\

%\vspace{-2mm}

\textbf{Example 3} ({\it Spherical van der Waerden} or {\it normal} scores). Let
$$\J_{\ell} \left(\frac{R^{(n)}_{\pm, t}}{n_R + 1}\S^{(n)}_{\pm, t}\right) = J_{\ell} \left(\frac{R^{(n)}_{\pm, t}}{n_R + 1}\right) \S^{(n)}_{\pm, t}, \quad  \ell  = 1, 2,\vspace{-2mm}$$
with $J_{\ell}(u) = \big((F^{\chi^2}_d)^{-1} (u)\big)^{1/2}\!$, where $F^{\chi^2}_d$ denotes the chi-square distribution function with~$d$ degrees of freedom. This yields   the {\it spherical van der Waerden} (vdW) {\it rank scores}, with 
% in Hallin and Paindaveine (2014, p. 2658), and we thus label the resulting scores {the vdW rank scores}.  In this case, 
 cross-covariance matrices\vspace{-2mm}  
\begin{align}
\tenq{\bGamma}_{i, \text{vdW}}^{(n)}(\bth) &= (n-i)^{-1}\!\! \sum_{t=i+1}^n \left[\big(F^{\chi^2}_d\big)^{-1}\! \left(\frac{R^{(n)}_{{\pms}, t}(\bth)}{n_R + 1}\right)\right]^{1/2}\! \left[\big(F^{\chi^2}_d\big)^{-1}\! \left(\frac{R^{(n)}_{{\pms}, t-i}(\bth)}{n_R + 1}\right)\right]^{1/2} 
\!\!  \S^{(n)}_{{\pms}, t} (\bth)\S^{(n)\prime}_{{\pms}, t-i}(\bth),\nonumber \\
& \qquad \qquad \qquad\quad\qquad \qquad \qquad\quad\qquad \qquad\qquad \qquad\qquad
\quad i = 1, \ldots , n-1. \label{vdWGam}
\end{align}
\vspace{-10mm}

\noindent In view of Propositions~\ref{Prop.Stilde} and~\ref{Prop.Wtilde},  the resulting tests 
%$\utilde{S}_{\J_1 \J_2}^{(n)}(\bth_0)$ and $\widehat{\utilde{W}}^{(n)}_{\J_1 \J_2}(\hat{\bth}^{(n)})$
 are asymptotically optimal under spherical Gaussian innovations.  
 
One also may wish to consider more sophisticated reference densities, such as skew-normal or skew-$t$ ones. The problem then is the absence, for non-spherical densities, of a closed-form expression for $\Q_\pm$.\footnote{Even for nonspherical elliptical densities,  no closed forms of ${\bf F}_\pm$ and ${\bf Q}_\pm$  are  available in the measure transportation literature.} That problem, in principle,  can be overcome by means of a numerical evaluation of ${\Q}_\pm$: increasingly powerful algorithms indeed exist for the numerical computation of optimal transports.   Assume that the (unspecified) innovation covariance matrix $\bSigma$ exists and is finite and consider, for instance, the scores associated with a  Gaussian~${\cal N}(\0,\bSigma)$ reference density (for which $\bvp_{\bSigma} ({\x})= \bSigma^{-1} \x$): \vspace{-2mm} 
 $$\J_1:= \bSigma^{-1} \circ \Q_{{\cal N}(\scriptstyle{{\bf 0},{\boldsymbol \Sigma}});\pm} , \qquad \J_2:=\Q_{{\cal N}(\scriptstyle{{\bf 0},{\boldsymbol \Sigma}});\pm}
\vspace{-2mm} $$
 where $\Q_{{\cal N}(\scriptstyle{{\bf 0},{\boldsymbol \Sigma}});\pm}$ denotes the center-outward quantile function associated with the ${\cal N}(\0,\bSigma)$ density. That quantile function is analytically unknown, but can be evaluated (at the $\mathfrak{G}^{(n)}$ gridpoints) with arbitrary precision by (a) estimating $\bSigma$ with some consistent estimator $\widehat\bSigma\n$ measurable with respect to the order statistic $\ZZ\n_{(.)}(\bth_0)$, then (b) simulating, independently of the data under study,  a  large sample of $N>>n$ i.i.d.\  ${\cal N}(\0,\widehat\bSigma\n)$ artificial observations, and (c) transporting these $N$ observations to    a regular grid $\mathfrak{G}^{(N)}$ refining $\mathfrak{G}^{(n)}$; that simulation,  for~$N$ large enough, provides an arbitrarily precise evaluation $\widehat{\Q}_\pm^{(N)}$ of  the values, over the~$\mathfrak{G}^{(n)}$ grid, of the center-outward quantile function   $\Q_{{\cal N}(\scriptstyle{{\bf 0},\widehat{\boldsymbol \Sigma}\n});\pm}$.  Based on the resulting scores, the   tests  of Propositions~\ref{Prop.Stilde}~{\it (iii)} and  \ref{Prop.Wtilde} {\it (iv)}  can be performed at conditional (on the order statistic) level $\alpha$ since $\widehat\bSigma\n$, as a function of the order statistic, is conditionally a constant, the value of which, in view of the independence between ${\F}\n_\pm$ and the order statistic, does not affect the distribution of the ranks and the signs. As $n$ and $N$ tend to infinity, these tests are locally asymptotically optimal under ${\cal N}(\0,\bSigma)$ innovation density. 
  
 The same idea can be implemented to construct tests based on skew-normal or  skew-$t$ scores.  Exploring\footnote{A very preliminary investigation of the Gaussian case described above shows that the efficiency gains over spherical-Gaussian-score-based methods is quite small.} the feasibility and the benefits (relative efficiencies only can be obtained numerically) of such an approach, however, is beyond the scope of this paper, and is left for further research.

\section{Gaussian tests}\label{GaussSec}

In this section, we briefly introduce the routinely-applied Gaussian test procedure (the hypothesis-testing counterpart of quasi-maximum-likelihood estimation (QMLE)), the performance of which will serve as a benchmark for our center-outward rank-based test procedure in the numerical analysis of Section~\ref{Sec:NumStud}.~These tests are developed in  Hallin and Painda\-veine~(2004a). We are focusing on the two cases  described in Sections~\ref{Sec.theta0} and \ref{Sec.testorder}.

\subsection{Gaussian test for $\bth = \bth_0$ ($\bth_0$ specified)}\label{Sec.Gauss.theta0}

We first consider the Gaussian solution to the problem considered in Section~\ref{Sec.theta0}, i.e., the   test of  the null hypothesis $\mathcal{H}^{(n)}_{\bth_0 ;{\mathcal N}}:=\{{\rm P}\n_{{\bth_0};f}\vert f~\text{Gaussian}\}$   ($\bth_0$   specified), against the alternative of VAR($p$) ($p_0\leq p \leq p_1$) with parameter $\bth \neq \bth_0$. 
Write  $\ZZ _t$ for $\ZZ_t^{(n)}(\bth_0)$ and consider the {\it residual cross-covariance matrix}\vspace{-2mm}
\begin{equation}\label{Gam}
{\bGamma}_{i; \mathcal{N}}^{(n)}(\bth_0) := (n-i)^{-1} \sum_{t=i+1}^n   \ZZ_{t} \ZZ^{\prime}_{t-i}, \quad i = 1, \ldots , n - 1,
\end{equation}
the  Gaussian counterpart to the rank-based $\tenq{\bGamma}_{i; \J_1, \J_2}^{(n)}(\bth_0)$ in \eqref{tildeGam1}.  Similarly denote by~${\bGamma}_{\mathcal{N}}^{(n)}(\bth_0)$ the Gaussian counterpart to  $\tenq{\bGamma}_{\J_1, \J_2}^{(n)}(\bth_0)$ (which results from substituting $\tenq{\bGamma}_{i; \J_1, \J_2}^{(n)}(\bth_0)$\linebreak with~${\bGamma}_{i; \mathcal{N}}^{(n)}(\bth_0)$). 
%As in Section~\ref{Sec.theta0},  we set $p_1 = p \geq p_0$ \color{blue}(similar to the mistake in Section~\ref{Sec.theta0}, we should delete $p_1 = p \geq p_0$ here) \color{red} what should we write here, then? \color{black}  for $ \Q^{(n)}_{\bth_0}$ defined in \eqref{defPQ}. 
Then, 
the Gaussian test statistic for $\mathcal{H}^{(n)}_{\bth_0 ;{\mathcal N}}$ is\vspace{-3mm} 
\begin{equation}\label{SGauss}
S_{\mathcal{N}}^{(n)}(\bth_0)  :=  \H_{\mathcal{N}}^{(n)\prime}(\bth_0) \left(  \Q^{(n)\prime}_{\bth_0}  \Q^{(n)}_{\bth_0} \right)^{-1}\H_{\mathcal{N}}^{(n)}(\bth_0),
\end{equation}
where $\H_{\mathcal{N}}^{(n)}(\bth_0)$ is the Gaussian counterpart of $\utilde{\H}_{\J_1 \J_2}^{(n)}(\bth_0)$,   obtained by substituting ${\bGamma}_{\mathcal{N}}^{(n)}(\bth_0)$ for  $\tenq{\bGamma}_{\J_1, \J_2}^{(n)}(\bth_0) - {\rm E}_{{\bth_0};f}  [\tenq{\bGamma}_{\J_1 \J_2}^{(n)}(\bth_0)]$ in~\eqref{Stilde}. The Gaussian test rejects $\mathcal{H}^{(n)}_{\bth_0 ;{\mathcal N}}$ whenever~$S_{\mathcal{N}}^{(n)}(\bth_0)$ exceeds the $(1-\alpha )$ quantile  of the chi-square distribution with~$d^2 p_1$ degree of freedom. That test is   locally asymptotically maximin at asymptotic level $\alpha$ against Gaussian alternatives; see Section~6.3 of \cite{HP04} for  details and a proof.

\subsection{Gaussian test for VAR order selection  (unspecified $\bth_0$)}\label{GIdSec}

Next let us consider  the Gaussian tests required for the sequential VAR order identification problem. The null  hypothesis of interest is thus $\mathcal{H}^{(n)}_{0;{\cal N}}:=\bigcup\big\{ \mathcal{H}^{(n)}_{\bth_0 ;{\cal N}} \big\vert\, \bth_0 \text{ of the form } \eqref{bth0}\big\}$; the general ideas leading to the test statistic  are very similar to those developed  in Section~\ref{Sec.testorder}. 

Some further notation is needed. Denote by $\bDelta_{\mathcal{N}}^{(n)}(\bth_0)$ the Gaussian central sequence, which is similar to~$\utilde{\bDelta}_{\J_1 \J_2}^{(n)}(\bth_0)$, with ${\bGamma}_{\mathcal{N}}^{(n)}(\bth_0)$  instead of $\tenq{\bGamma}_{\J_1, \J_2}^{(n)}(\bth_0) - {\rm E}_{{\bth_0};f}  [\tenq{\bGamma}_{\J_1 \J_2}^{(n)}(\bth_0)]$. Split~$\bDelta_{\mathcal{N}}^{(n)}(\bth_0)$ into $(\bDelta_{I; \mathcal{N}}^{(n)\prime}(\bth_0), \bDelta_{II; \mathcal{N}}^{(n)\prime}(\bth_0))^\prime$ with 
%$\bDelta_{I; \mathcal{N}}^{(n)\prime}(\bth_0)$ and $\bDelta_{II; \mathcal{N}}^{(n)\prime}(\bth_0)$ being
 $d^2 p_0$- and $d^2(p_1 - p_0)$-dimensional subvectors, respectively. 
Letting  
 $\bLam^{(n)}_{\mathcal{N}; \bth_0} := \T^{(n)}_{\bth_0} \left( \I_{n-1} \otimes \L_{\bth_0}\n \right) \T^{(n)\prime}_{\bth_0},$
where\vspace{-2mm} 
$$\L_{\bth_0}\n := (n-1)^{-1} \sum_{t=2}^n \text{vec} \left( \ZZ_t\n \ZZ_{t-1}^{(n)\prime} \right) \left( \text{vec} \left( \ZZ_t\n \ZZ_{t-1}^{(n)\prime} \right) \right)^\prime,$$ 
 partition  it into a block matrix\vspace{-3mm}
\[
\bLam^{(n)}_{\mathcal{N}; \bth_0}  =
\begin{bmatrix}
\bLam^{(n)}_{11; \mathcal{N}; \bth_0} & \bLam^{(n)}_{12; \mathcal{N}; \bth_0} \\
\bLam^{(n)}_{21; \mathcal{N}; \bth_0}  & \bLam^{(n)}_{22; \mathcal{N}; \bth_0}
\end{bmatrix}
\]
with square diagonal blocks $\bLam^{(n)}_{11; \mathcal{N}; \bth_0}$ and $\bLam^{(n)}_{22; \mathcal{N}; \bth_0}$ of size $d^2 p_0$ and $d^2 (p_1 - p_0)$, respectively.
With these notations, define\vspace{-2mm} %a counterpart of $\utilde{\bDelta}^{(n)*}_{II; \J_1, \J_2}(\bth_0)$ in \eqref{defDeltastar} as
\begin{equation}\label{defDeltastarGauss}
{\bDelta}^{(n)*}_{II; \mathcal{N}}(\bth_0) := {\bDelta}^{(n)}_{II; \mathcal{N}}(\bth_0) - \bLam^{(n)}_{21; \mathcal{N}; \bth_0}  (\bLam^{(n)}_{11; \mathcal{N}; \bth_0})^{-1} \utilde{\bDelta}^{(n)}_{I; \mathcal{N}}(\bth_0),
\end{equation}
\vspace{-12mm}

%and let
$$
\bLam^{(n)*}_{II;\mathcal{N}; \bth_0} 
 := \bLam^{(n)}_{22; \mathcal{N}; \bth_0} -
\bLam^{(n)}_{21; \mathcal{N}; \bth_0} (\bLam^{(n)}_{11; \mathcal{N}; \bth_0})^{-1} \bLam^{(n)}_{12; \mathcal{N}; \bth_0},
\vspace{-2mm}$$
and%Then, the Gaussian test statistic is defined as
\begin{equation}\label{Gaussiantest}
W\n_{\mathcal{N}}(\bth_0) :=  {\bDelta}^{(n)*^\prime}_{II; \mathcal{N}}(\bth_0) \left(  \bLam^{(n)*}_{II;\mathcal{N}; \bth_0} \right)^{-1} {\bDelta}^{(n)*}_{II; \mathcal{N}}(\bth_0).
\end{equation}
Finally, let $\hat{\bth}^{(n)}$ be a constrained estimator of $\bth_0$ satisfying the same assumptions as   in Section~\ref{Sec.testorder}. 
The Gaussian  test rejects the null hypothesis $\mathcal{H}^{(n)}_{0;{\cal N}}$ whenever $W\n_{\mathcal{N}}(\hat{\bth}^{(n)})$ exceeds the~$(1-\alpha )$ quantile  of the chi-square distribution with $d^2 (p_1 - p_0)$  degree of freedom and is locally asymptotically most stringent against Gaussian alternatives; see \citet[Section~5.2]{HP05}. 

When   testing the null hypothesis of white noise  (VAR($0$)) against VAR($p_1$) dependence, we have $p_0 = 0$ and $\bth_0 = \0$ (no need for an estimator $\hat{\bth}^{(n)}$  of $\bth_0$). Then the Gaussian test statistics  
 $S\n_{\mathcal{N}}(\0)$ and $W\n_{\mathcal{N}}(\0)$ coincide, taking the form\vspace{-2mm}  
\begin{equation}
n {\bGamma}_{\mathcal{N}}^{(p_1, n)\prime}(\0)
( \I_{p_1} \otimes  \L_{\0}\n )^{-1} {\bGamma}_{\mathcal{N}}^{(p_1, n)}(\0), \label{Eq: Wn}
\end{equation}
\vspace{-11mm}
 
\noindent where\vspace{-2mm} 
$${\bGamma}_{\mathcal{N}}^{(p_1, n)}(\bth_0) := n^{-1/2} \left( (n-1)^{1/2} (\text{vec}\,  {\bGamma}_{1, \mathcal{N}}^{(n)}(\bth_0))^\prime, \ldots ,    (n-p_1)^{1/2} (\text{vec}\,  {\bGamma}_{p_1, \mathcal{N}}^{(n)}(\bth_0))^\prime \right)^\prime\vspace{-2mm}$$
is the truncated version of ${\bGamma}_{\mathcal{N}}^{(n)}(\bth_0)$; if $p_1 =1$,  the   critical value is the $(1-\alpha )$ quantile $\chi^2_{d^2;1-\alpha}$ of the chi-square distribution with $d^2$  
degrees of freedom.

%\color{red}
%I propose that Davide goes through this section (maybe he already did?) My concern is the very severe bias (null rejection frequency $<<\alpha$) of rank-based tests in Tables 1 and 2. This is disturbing. Why does it occur? Any idea on how to correct it? 
%
%As for the real data case, it is mildly convincing, as all tests yield the same conclusions. Ideally, we should come up with an example where the VAR order detected differ (e.g.,   the rank-based method detecting a VAR dependence where the Gaussian test concludes for white noise).
%\color{black}
%

\section{Examples and numerical results} \label{Sec:NumStud}

We illustrate numerically the performance of our test statistics in two benchmark inferential problems:
testing for serial dependency in a multivariate time series and selecting the order of autoregression in a VAR model.
As discussed in \cite{HP04}, these testing problems are common to many econometric applications. 

%\color{blue} (Not sure if we need so much space here and at the beginning of  Section~\ref{sec.NumWhite} and Section~\ref{VARorder} describing model~\eqref{6themodel}? I think maybe it is not much relevant since we are not dealing with any estimation of $\betab$). \color{olive} not a big problem, I think. But the terminology ``rank-based multivariate  Durbin-Watson test" (with reference) should be used as a marketing tool. \color{black}

To elaborate further,   consider the $d$-variate general regression model with VAR error terms. Under this model, the observation is an $n\times d$ array\vspace{-1mm}
$$  \Y:=
             \left(
             \begin{array}{cccc}
             Y_{1,1} & Y_{1,2} & \ldots & Y_{1,d} \\
             \vdots & \vdots & & \vdots \\
             Y_{n,1} & Y_{n,2} & \ldots & Y_{n,d} \\
             \end{array}
             \right):= \left(
\begin{array}{c}{\Y}_1\pr  \\  \vdots \\ {\Y}_n\pr 
            \end{array}
 \right)
$$
of $d$-variate random vectors  $\Y_t$ satisfying \vspace{-2mm} %, where at each time point $t$ the model 
\begin{equation} \label{6themodel}
             \Y_t=g(\vbf_t,\betab)+\X_t.  \vspace{-2mm}
\end{equation}
where the regression  function $g : \mathbb{R}^q \times \mathbb{R}^k \to \mathbb{R}^d$ has known functional form and depends on an unknown parameter $\betab \in \mathbb{R}^k$ and on observed covariates $\vbf_t\in\mathbb{R}^q$. % (they can also be a matrix), which are assumed to be measurable at time $t$. 
%
%$$
%             \V\!\! :=\!\!
%             \left(\! 
%             \begin{array}{cccc}
%             v_{1,1} & v_{1,2} & \ldots & v_{1,m} \\
%             \vdots & \vdots & & \vdots \\
%             v_{n,1} & v_{n,2} & \ldots & v_{n,m} \\
%             \end{array}
%            \! \right)\!\!
%:= \!\!\left(\!
%\begin{array}{c}{\bf v}_1\pr  \\  \vdots \\ {\bf v}_n\pr 
%            \end{array}
%\! \right)\!\!
%\quad\mbox {\rm and}\quad
%         \!\!    \betab\!\! :=\!\! 
%             \left(\!
%             \begin{array}{cccc}
%             \beta_{1,1} & \beta_{1,2} & \ldots & \beta_{1,d} \\
%             \vdots & \vdots & & \vdots \\
%             \beta_{m,1} & \beta_{m,2} & \ldots & \beta_{m,d} \\
%             \end{array}
%            \! \right) \!\!
%:=\!\! \left(\!
%\begin{array}{c}\betab_1\pr  \\  \vdots \\ \betab_m\pr 
%            \end{array}
%\! \right)
%$$
%denote an $n\times m$ matrix of constants (the design 
%matrix), and the $m\times d$ parameter, respectively. 
Instead of the
traditional
             assumption that the error term
%             $$
%             \X:=
%             \left(
%             \begin{array}{cccc}
%             X_{1,1} & X_{1,2} & \ldots & X_{1,d} \\
%             \vdots & \vdots & & \vdots \\
%             X_{n,1} & X_{n,2} & \ldots & X_{n,d} \\
%             \end{array}
%             \right):= \left(
%\begin{array}{c}{\X}_1\pr  \\  \vdots \\ {\X}_n\pr 
%            \end{array}
% \right)
%             $$
             is white noise, we rather assume $\{{\X}_t ; t = 1,\ldots, n\}$ to be the finite realization (of length $n$) of some 
VAR$(p)$ process generated by \vspace{-2mm}
\[%\begin{equation*} 
\X_t = \sum _{i=1}^{p}\A_i\X _{t-i} +\bepsilon_t , \quad \quad
t\in\Z ,\vspace{-2mm}
\]%\end{equation*}
  where    $\{\bepsilon_t;  t\in\Z \}$ is  $d$-dimensional white noise  with density $f$, satisfying Assumption~\ref{ass.den}.  %Under (\ref{6themodel}) and (\ref{6ARMA}),  
%             \begin{equation}\label{Wold}
%      \Yb _t = \betab\pr \vb_t    + \sum_{u=0}^{t-1}\Gb _{u}\bepsilon _{t-u} + {\bf r}_t    ,\quad t=1,\ldots , n      ;    \end{equation}
%where  $\Gb _{u}$ are Green's matrices of the VAR operator. The remainder term ${\bf r}_t$ is related to the influence of the unobserved initial values $\Vb _0 ,\ldots , \Vb _{-p+1}$;  it is easy to see that, under Assumption 1, we have  that $\lim _{t\rightarrow\infty}\Lambda ^t{\bf r}_t$ is bounded in probability, where $1<\Lambda $ is the modulus of the smallest   root of the characteristic polynomial associated with (\ref{6ARMA}). 

A similar framework is considered in \cite{HP04b,HP05}, where the function $g$ is linear in   $\betab$ (hence, without any loss of generality, in the covariates).   Hallin and Paindaveine moreover assume that the innovation density $f$ is elliptical and  consider test procedures %which are
 based on Mahalanobis ranks and signs. This very restrictive assumption of elliptical symmetry is precisely what we are dropping here as we only assume that~$f\in{\cal F}_d^*$.

%Unfortunately, the validity of those procedure is limited to the family of elliptical innovation densities, under which Mahalanobis ranks and signs are distribution-free. In the numerical exercises of this section, we illustrate the performance of our test statistics based on ceter-outward quantities in a framework that extends the one discussed in Hallin and Paindaveine papers: we allow for a nonlinear specification of the conditional mean (the $g$ function) and for a generic innovation density belonging to ${\cal F}_d$---the considered $f$ is not necessarily an elliptic one, thus we do not need an estimator of the scatter matrix. 
%

\subsection{A multiple-output rank-based  Durbin-Watson test}\label{sec.NumWhite}

\subsubsection{Testing for serial correlation in (nonlinear) multiple-output regression}
%\textit{Inference problem.}
%\color{olive} 
A classical problem in single-output linear regression is the so-called Durbin-Watson problem of testing the null hypothesis under which $\X _t$, in the single-output  linear version of \eqref{6themodel} with unspecified $\betab$, is second-order white noise against the alternative of VAR(1) dependence: see \cite{DWI,DWII}. Here, we extend that problem to nonlinear multiple-output regression and propose a rank-based solution. For the sake of simplicity, we limit ourselves to VAR(1) alternatives (writing $\A$  instead of $\A_1$), but extensions to higher-order VAR($p$) are straightforward. 

It is easy to see that, under Assumption~\ref{ass.den} and appropriate conditions on $g$ and the asymptotic behavior of the covariates, model \eqref{6themodel} is LAN with block-diagonal information matrix (a nonserial $\betab$-block and a serial one related with $\bth=\text{vec}\A$).  LAN then implies that the unknown $\betab$ safely can be replaced by any root-$n$ consistent (and, in principle, locally discrete) estimator $\hat\betab\n$ without having any local impact on the $\bth$-part of the central sequence (see  \cite{HP04b} for the case of a linear $g$; the only difference is that the~$\betab$ part of the central sequence here has an additional $g^\prime$ factor which does not affect that fact that its covariance with the $\bth$ part vanishes). As a consequence, one can construct a test as of $\betab$ were known, then safely  replace $\betab$ with $\hat\betab\n$ in the resulting test statistic. 

The traditional implementation of  Durbin-Watson tests involves least-squares estimators of~$\betab$. We rather  suggest  robust root-$n$ consistent M-estimators  $\hat\betab\n$ such as the  bounded-influence M-estimator for dynamic location models proposed in~\cite{MRT}.

\subsubsection{Numerical results: size, power,  and bias correction}\label{sec.d2}
%\textit{Numerical results: size and power.} 
% Let us consider the model (\ref{6themodel}) and write $\A$ instead of $\A_1$. We consider t
 The Durbin-Watson   problem thus is about testing 
 %where the null hypothesis is that $\{\X _t\}$ is a white noise, that is,
  $\bth ={\bf 0}$ against   $\bth \neq {\bf 0}$, on the basis of  estimated residuals $Y_t-g(\vbf_t,\hat\betab\n)$ as if they were the exact residuals $\X_t=\Y_t-g(\vbf_t, \betab)$. Under the null hypothesis, the $\X_t$'s are white noise; under the alternative, they are VAR(1). The Gaussian solution (which is  routinely applied) is based on the Gaussian test statistic~\eqref{Eq: Wn}, its rank-based competitors on   test statistics of the form~\eqref{tildeW.white}, for which we consider the spherical van der Waerden,  Spearman, and sign scores as described in Section~\ref{Sec: examples}. 
  
To simplify the numerical implementation and avoid  specifying any particular $g$,  our simulations directly proceed with the residuals $\X_t$ in \eqref{6themodel}. 
% we set $ g(\vbf_t,\betab) \equiv \boldsymbol{0}$, for every $t$---so, there is no need for the estimation step of the parameter $\betab$ and the numerical results are not affected by the properties (e.g. efficiency and finite sample bias) of the estimator of $\betab$. To compare the size and power of the Gaussian and center-outward rank-based tests of white noise against VAR($1$), 
We thus generated $N = 1000$ replications of size $n = 800$ ($n_R = 20$, $n_S = 40$, and $n_0 = 0$)  from the bivariate models \vspace{-2mm}
\begin{equation}\label{bimod} \X_t - \ell \A \X_{t-1}= \bepsilon_t, \quad \ell = 0, 1, 2 \vspace{-2mm}\end{equation}
with $\text{vec}(\A) = (0.05, -0.01, 0.02, 0.05)^\prime$;  
%and $N = 1000$ replications of size $n = 1400$  ($n_R = 21$, $n_S = 66$, and $n_0 = 14$) from 
%the trivariate VAR models\vspace{-2mm}
%$$ \X_t - \ell \A \X_{t-1}= \bepsilon_t, \quad \ell = 0, 1, 1.5,\vspace{-2mm}$$
%with $\text{vec}(\A) = (0.05, 0.01, 0.011, 0.01, 0.02, 0.01, -0.01, 0.013, 0.033)^\prime$.
 Innovation densities   are   spherical normal and Student   with $3$ degrees of freedom (denoted as $t_3$), mixtures of  normal,   and skew-$t_3$ densities. The mixtures are of the form\vspace{-2mm}
$$%\begin{equation}
\frac{3}{8}   {\cal N}(\bmu_1, \bSigma_1) + \frac{3}{8}   {\cal N}(\bmu_2, \bSigma_2) + \frac{1}{4} {\cal N}(\bmu_3, \bSigma_3),
$$%\label{Eq. Mixture}
%\end{equation}
with $\bmu_1 = (-5, 0)^\prime,\ \bmu_2 = (5, 0)^\prime, \ \bmu_3 = (0, 0)^\prime$ and\vspace{-2mm}  
$$\bSigma_1 = 
\begin{bmatrix}
7 & 5\vspace{-1mm} \\
5 & 5\vspace{-1mm}
\end{bmatrix}, \  
\bSigma_2 = 
\begin{bmatrix}
7 &  -6\vspace{-1mm} \\
6 &  6\vspace{-1mm}
\end{bmatrix}, \ 
\bSigma_3 = 
\begin{bmatrix}
4 & 0\vspace{-1mm} \\
0 & 3\vspace{-1mm}
\end{bmatrix}
\vspace{-2mm}$$%\]
%for $d=2$, with
% $\bmu_1 = (-5, -5, 0)^\prime,\ \bmu_2 = (5, 5, 2)^\prime, \ \bmu_3 = (0, 0, -3)^\prime$ 
%and\vspace{-2mm} 
% $$\bSigma_1 = 
%\begin{bmatrix}
%7 & 3 & 5\vspace{-1mm} \\
%3 & 6 & 1\vspace{-1mm} \\
%5 & 1 & 7\vspace{-1mm}
%\end{bmatrix}, \  
%\bSigma_2 = 
%\begin{bmatrix}
%7 & -5 & -3\vspace{-1mm} \\
%-5 & 7 & 4\vspace{-1mm} \\
%-3 & 4 & 5\vspace{-1mm} 
%\end{bmatrix}, \ 
%\bSigma_3 = 
%\begin{bmatrix}
%4 & 0 & 0\vspace{-1mm}\\
%0 & 3  & 0\vspace{-1mm} \\
%0 & 0 & 1\vspace{-1mm}
%\end{bmatrix}
%\vspace{-2mm}$$
% for $d=3$. 
 The  $d$-dimensional skew-$t_\nu$  distribution  has  density\vspace{-4mm} 
\begin{align}\label{St.density}
%&  \\
f (\z; \bxi, \bSigma, \al, \nu) &\vspace{-2mm}\\ 
:=&2 {\rm det}({\mbf w)}^{-1} 
 t_d(\x; \bSigma, \nu) T\left(\al^\prime \x  \big({{(\nu + d)}/{(\nu + \x^\prime \bSigma^{-1} \x)}}\big)^{1/2}; \nu + d \right),\quad\!\!\! \z \in \R^d, %\vspace{-2mm}
\nonumber\end{align}
\vspace{-9mm}

\noindent   (see \citealt{AC03}) where  $\x = {\mbf w}^{-1}(\z - \bxi)$;  ${\mbf w} = \left({\rm diag} (\bSigma)\right)^{1/2}$, $\bxi \in \R^d$, and~$\al \in \R^d$      are scale,  location, and shape   parameters respectively,     $T(y; \nu)$ denotes the univariate $t_\nu$ distribution function, and\vspace{-2mm}
%$t_d(\x; \bSigma, \nu)$   the  multivariate $t_\nu$ density, viz. 
%$$t_d(\x; \bSigma, \nu) := \frac{\Gamma((\nu+d)/2)}{(\nu \pi)^{d/2} \Gamma(\nu/2) {\rm det}(\bSigma)^{1/2}} \left( 1+ \frac{\x^\prime \bSigma^{-1} \x}{\nu} \right)^{-(\nu+d)/2},\quad \x \in \R^d.\vspace{-2mm}$$ 
$$t_d(\x; \bSigma, \nu) := \frac{\Gamma((\nu+d)/2)}{(\nu \pi)^{d/2} \Gamma(\nu/2) {\rm det}(\bSigma)^{1/2}} \left( 1+ \frac{\x^\prime {\mbf w} \bSigma^{-1} {\mbf w} \x}{\nu} \right)^{-(\nu+d)/2},\quad \x \in \R^d.\vspace{-2mm}$$
Here ($d=2$), the values $\bxi = \0$, $\al = (5, 2)^\prime$, and 
$\bSigma = 
\begin{bmatrix}
7 & 4\vspace{-1mm} \\
4 & 5\vspace{-1mm}
\end{bmatrix}
$  were selected. \vspace{1mm}

 %%%%%%%%%%%%%%%%%%%%%%%%%%%%%%
\begin{table}[!ht]

\caption{Rejection frequencies (out of $N=1000$ replications),  under   values~$\ell \A$, $\ell=0, 1, 2$  of the VAR(1) autoregression matrix, for the bivariate model \eqref{bimod}, and various innovation densities, of the Gaussian, vdW, bias-corrected vdW, Spearman, bias-corrected Spearman, and sign  tests of white noise against VAR($1$); the sample size is $n = 800$   ($n_R = 20$, $n_S = 40$, and~$n_0 = 0$); the nominal level is $\alpha=5\%$;  rejection frequencies of the sign tests with~$n_R=1$ are also included. Permutational critical values are based on $M=5000$ random permutations. }\label{Tab.white}\vspace{2mm}

\centering

\footnotesize

\begin{tabular}{lcccc | lcccc}

\hline\hline

$\quad\ f$ &\hspace{-31mm} Test     &\hspace{-2mm}  $\0$  & $\A$  & $2\A$ &   $\quad\ f$      &\hspace{-31mm} Test     &\hspace{-2mm} $\0$  & $\A$  & $2\A$ \\ \hline\hline

\underline{Normal}\vspace{-1mm}&  &  &  & & \underline{Mixture}    & &  &   &   \\

 &\hspace{-31mm}  Gaussian &\hspace{-2mm}  0.056 & 0.353 & 0.933 &     &\hspace{-31mm} Gaussian &\hspace{-2mm}  0.047 & 0.391 & 0.946 \\

       &\hspace{-31mm} vdW      &\hspace{-2mm}  0.029 & 0.231 & 0.880 &            &\hspace{-31mm} vdW      &\hspace{-2mm}  0.024 & 0.562 & 0.998 \\

     {} &{\hspace{-31mm} bias-corrected vdW} &\hspace{-2mm} {0.055} &{0.331} &{0.920} & {} &{\hspace{-31mm} bias-corrected vdW} &\hspace{-2mm} {0.056} &{0.668} &{0.998} \\

% &              (bias-corrected)                        &                        &                        &                        &  &                          (bias-corrected)             &                        &                        &                      \\

       &\hspace{-31mm} Spearman &\hspace{-2mm}  0.052 & 0.281 & 0.898 &            &\hspace{-31mm} Spearman &\hspace{-2mm}  0.048 & 0.654 & 0.998  \\

      {} &\hspace{-31mm} bias-corrected Spearman &\hspace{-2mm} {0.060} &{0.319} &{0.909} & {} &\hspace{-31mm} bias-corrected Spearman &\hspace{-2mm} {0.060} &{0.687} &{0.999} \\

%                  &         (bias-corrected)                               &                        &                        &                        &                   &                                       (bias-corrected) &                        &                        &               \\

       &\hspace{-31mm}  Sign   ($n_S =40$)   &\hspace{-2mm}  0.052 & 0.222 & 0.726 &            &\hspace{-31mm}  Sign  ($n_S =40$)    &\hspace{-2mm}  0.053 & 0.284 & 0.863 \\

%   &                                        &                        &                        &                        &                   &                                         &                        &                        &               \\ 

%   &\hspace{-31mm} {bias-corrected Sign  ($n_R =20$)} &\hspace{-2mm} {0.049} &{0.213} &{0.710} &  &\hspace{-31mm} {bias-corrected Sign  ($n_R =20$)} &\hspace{-2mm} {0.050} &{0.276} &{0.853} \\

        &\hspace{-31mm}  {Sign ($n_S =800$)} &\hspace{-2mm} {0.053} &{0.218} & {0.734} &  &\hspace{-31mm}  {Sign ($n_S =800$)} &\hspace{-2mm}  {0.048} &{0.247} &{0.788} \\

%             &         (bias-corrected I)                               &                        &                        &                        &                   &                                       (bias-corrected I) &                        &                        &               \\

%             & {\hspace{-31mm} bias-corrected Sign ($n_R =1$)} &\hspace{-2mm} {0.052} &{0.208} &{0.727} &  &{\hspace{-31mm} bias-corrected Sign ($n_R =1$)} &\hspace{-2mm} {0.045} &{0.241} &{0.783} \\

%             &         (bias-corrected II)                               &                        &                        &                        &                   &                                       (bias-corrected II) &                        &                        &               \\

   \hline

$ \underline{\, t_3\,}$  &  &  &  &  & \underline{Skew-$t_3\,$} & &  &   &   \\

 &\hspace{-31mm} Gaussian &\hspace{-2mm} 0.043 & 0.326 & 0.918 &  &\hspace{-31mm} Gaussian &\hspace{-2mm} 0.051 & 0.306 & 0.898 \\

       &\hspace{-31mm} vdW      &\hspace{-2mm} 0.026 & 0.331 & 0.968 &            &\hspace{-31mm} vdW      &\hspace{-2mm} 0.025 & 0.476 & 0.997 \\

      {} &{\hspace{-31mm}bias-corrected vdW} &\hspace{-2mm}{0.055} &{0.417} &{0.985} & {} &{\hspace{-31mm} bias-corrected vdW} &\hspace{-2mm}{0.043} &{0.590} &{0.999} \\

%                  &         (bias-corrected)                               &                        &                        &                        &                   &                                       (bias-corrected) &                        &                        &               \\  

       &\hspace{-31mm} Spearman &\hspace{-2mm} 0.041 & 0.383 & 0.975 &            &\hspace{-31mm} Spearman &\hspace{-2mm} 0.030 & 0.543 & 0.999 \\

{} & {\hspace{-31mm} bias-corrected Spearman} &\hspace{-2mm} {0.053} & {0.398} & {0.979} &  {} & {\hspace{-31mm} bias-corrected Spearman} &\hspace{-2mm} {0.036} & {0.573} & {0.999} \\

%                  &         (bias-corrected)                               &                        &                        &                        &                   &                                       (bias-corrected) &                        &                        &               \\

       &\hspace{-31mm}  Sign   ($n_S =40$)   &\hspace{-2mm} 0.055 & 0.325 & 0.929 &            &\hspace{-31mm}  Sign   ($n_S =40$)   &\hspace{-2mm} 0.050 & 0.367 & 0.945 \\

% &          ($n_R =1$)                              &                        &                        &                        &                   &                                       ($n_R =1$) &                        &                        &               \\ 

%   & \hspace{-31mm}  bias-corrected Sign ($n_R =20$)&\hspace{-2mm} {0.051} & {0.308} & {0.922} &  &  \hspace{-31mm}  bias-corrected Sign ($n_R =20$) &\hspace{-2mm} {0.042} & {0.350} & {0.935} \\

        & \hspace{-31mm}  {Sign ($n_S =800$)} &\hspace{-2mm} {0.056} & {0.335} & {0.929} &  &\hspace{-31mm}  {Sign ($n_S =800$)} &\hspace{-2mm} {0.050} & {0.374} & {0.945} \\

%       &         (})                               &                        &                        &                        &                   &                                       (bias-corrected I) &                        &                        &               \\

%       & \hspace{-31mm}  bias-corrected Sign ($n_R =1$) &\hspace{-2mm} {0.054} & {0.326} & {0.925} &  &  \hspace{-31mm}  bias-corrected Sign ($n_R =1$)&\hspace{-2mm} {0.047} & {0.363} & {0.942} \\

%         &         (bias-corrected II)                               &                        &                        &                        &                   &                                       (bias-corrected II) &                        &                        &               \\

 \hline\hline

\underline{AOs}  (${\mbf s} = (6, 6)^\prime$) &  &  &  &  &  \underline{AOs}  (${\mbf s} = (9, 9)^\prime$) &  &  &     &   \\

 &\hspace{-31mm}  Gaussian &\hspace{-2mm} 0.179 & 0.221 & 0.573 &    &\hspace{-31mm}  Gaussian &\hspace{-2mm} 0.417 & 0.400   & 0.658 \\

                              &\hspace{-31mm}  vdW      &\hspace{-2mm} 0.035 & 0.131 & 0.711 &                            &\hspace{-31mm}  vdW      &\hspace{-2mm} 0.023 & 0.154 & 0.678 \\

   {} &{\hspace{-31mm}bias-corrected vdW}&\hspace{-2mm} {0.059} & {0.188} & {0.790} &  {} & {\hspace{-31mm}bias-corrected vdW} &\hspace{-2mm} {0.060} & {0.217} & {0.769} \\

%                     &         (bias-corrected)                               &                        &                        &                        &                   &                                       (bias-corrected) &                        &                        &               \\  

                            &\hspace{-31mm} Spearman &\hspace{-2mm} 0.049 & 0.166 & 0.756 &                             &\hspace{-31mm} Spearman &\hspace{-2mm} 0.041 & 0.197 & 0.727 \\

                            {} &\hspace{-31mm}{bias-corrected Spearman} &\hspace{-2mm} {0.060} & {0.187} & {0.777} &  {} &\hspace{-31mm}{bias-corrected Spearman} &\hspace{-2mm} {0.059} & {0.218} & {0.761} \\

%                     &         (bias-corrected)                               &                        &                        &                        &                   &                                       (bias-corrected) &                        &                        &               \\  

                            &\hspace{-31mm} Sign    ($n_S =40$)  &\hspace{-2mm} 0.069 & 0.160 & 0.626 &                             &\hspace{-31mm} Sign   ($n_S =40$)   &\hspace{-2mm} 0.051 & 0.191 & 0.605         \\

%&          ($n_R =1$)                              &                        &                        &                        &                   &                                       ($n_R =1$) &                        &                        &               \\ 

%&\hspace{-31mm}{bias-corrected Sign ($n_R =20$)} & \hspace{-2mm}{0.065} & {0.149} & {0.601} &  &\hspace{-31mm}{bias-corrected Sign ($n_R =20$)} &\hspace{-2mm}{0.044} & {0.177} & {0.589} \\

                            &\hspace{-31mm}{Sign ($n_S =800$)} &\hspace{-2mm} {0.060} & {0.145} & {0.601} &  &\hspace{-31mm}{Sign ($n_S =800$)} &\hspace{-2mm} {0.055} & {0.184} & {0.593} \\

%&         (bias-corrected I)                               &                        &                        &                        &                   &                                       (bias-corrected I) &                        &                        &               \\

%&\hspace{-31mm}{bias-corrected Sign ($n_R =1$)} &\hspace{-2mm}{0.060} & {0.141} & {0.593} &  &\hspace{-31mm}{bias-corrected Sign ($n_R =1$)}  &\hspace{-2mm}{0.051} & {0.179} & {0.585} \\

%&         (bias-corrected II)                               &                        &                        &                        &                   &                                       (bias-corrected II) &                        &                        &               \\

\hline\hline

\end{tabular}

\end{table}

%%%%%%%%%%%%%%%%%%%%%%%%%%%%

To investigate the robustness of the center-outward rank-based tests,  %compared with the Gaussian one, 
we also considered the case of spherical Gaussian   $\X_t$'s  contaminated by additive outliers (AOs). More precisely, we generated observations of the form 
 $\{\X^*_t=  \X_t + \1_{\{t = h\}} {\mbf s}\}$, where   $h$ and ${\mbf s}$ denote the location and size of the AOs, respectively. We set $h$ in order to obtain  
$5\%$ of equally spaced~AOs and put ${\mbf s} = (6, 6)^\prime$ and $(9, 9)^\prime$. All contaminated 
observations were demeaned prior to the implementation of the testing procedures.

\begin{figure}[!t]
     \centering
     \begin{subfigure}[b]{0.4\textwidth}
         \centering
         \includegraphics[width=\textwidth]{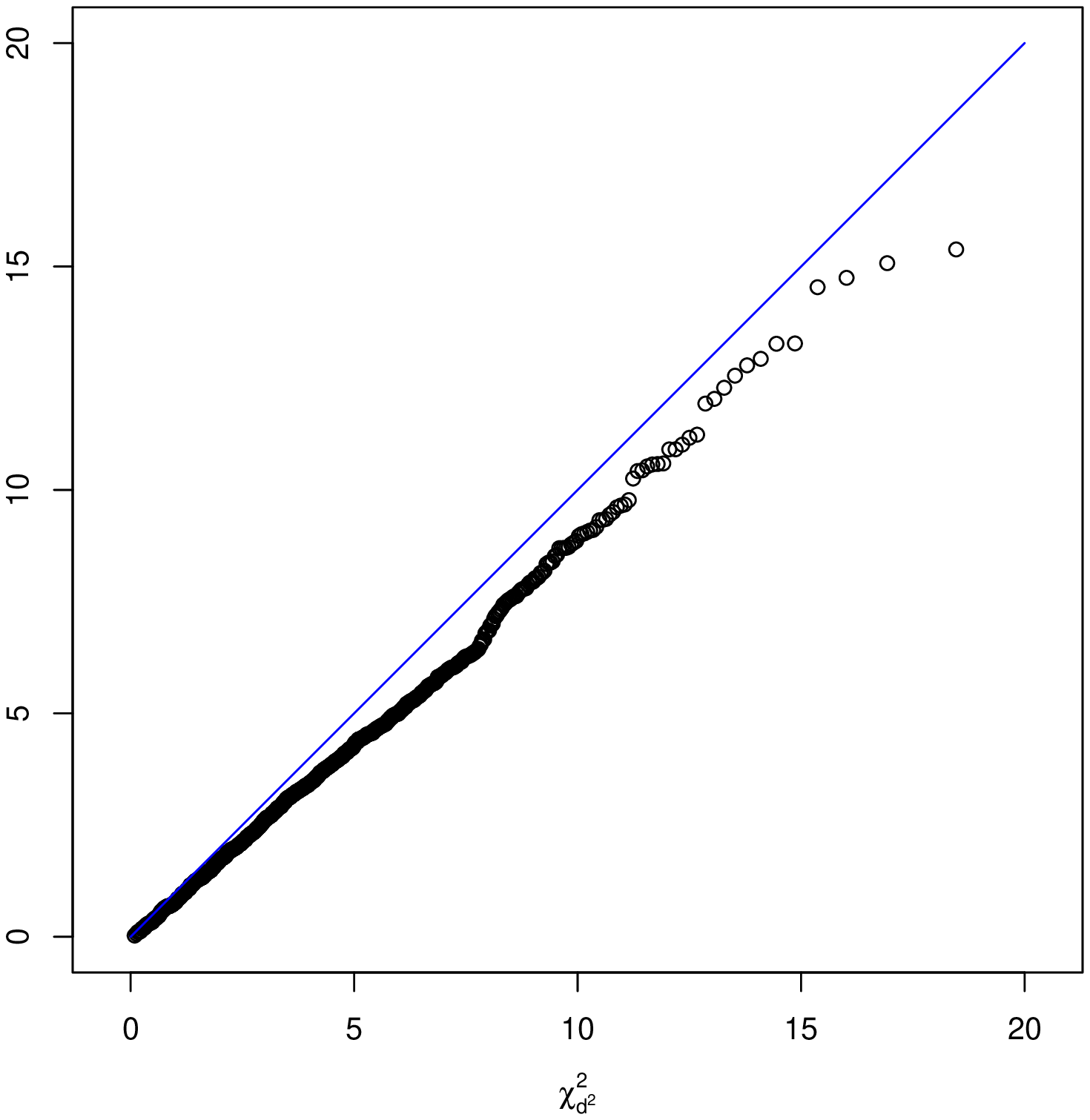}
        % \caption{$y=x$}
  %       \label{fig:y equals x}
     \end{subfigure}
%     \hfill
     \begin{subfigure}[b]{0.4\textwidth}
         \centering
         \includegraphics[width=\textwidth]{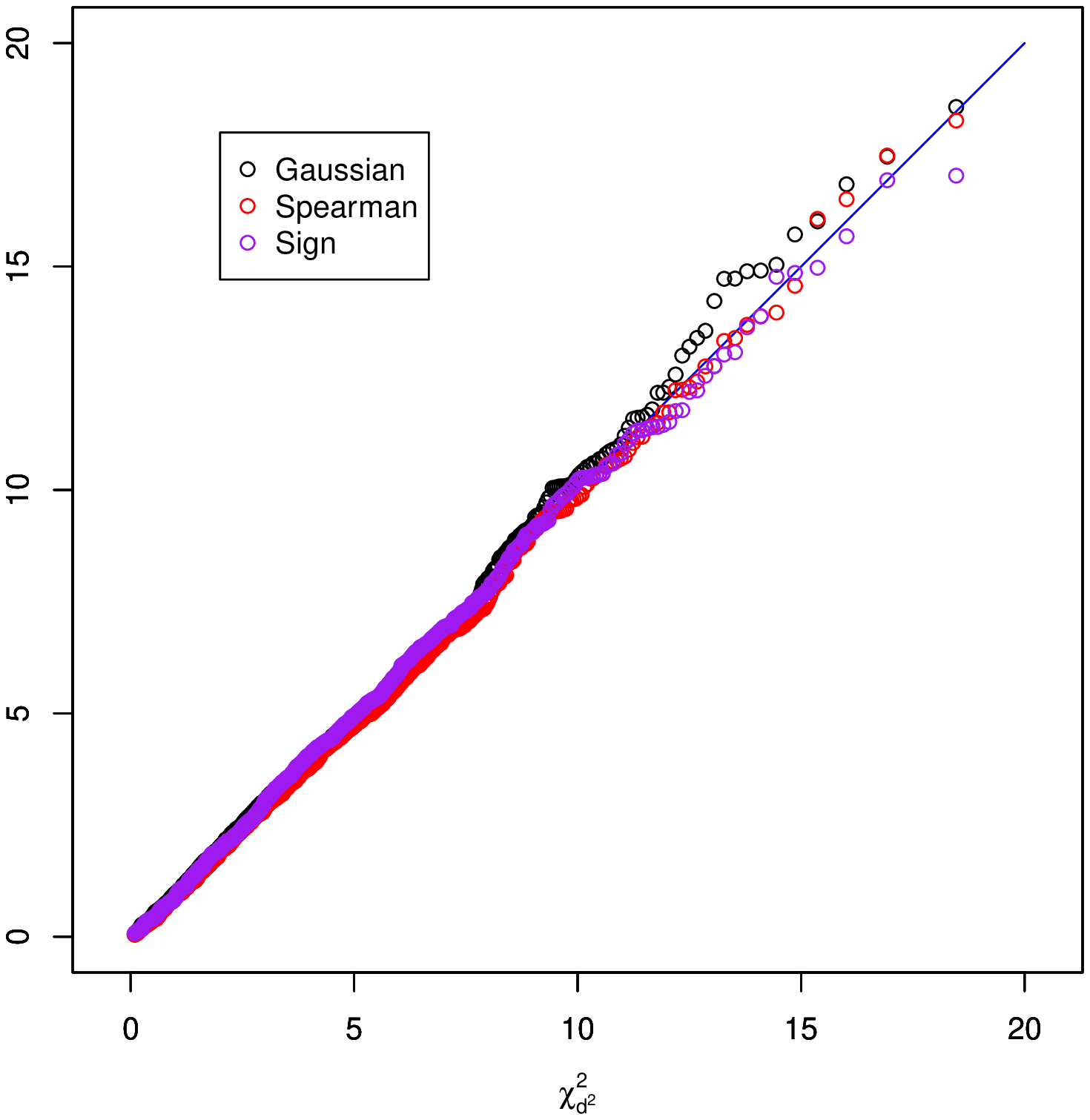}
       %  \caption{$y=3sinx$}
      %   \label{fig:three sin x}
     \end{subfigure}
%     \\
%      \begin{subfigure}[b]{0.4\textwidth}
%         \centering
%         \includegraphics[width=\textwidth]{d3vdWQQplot}
%        % \caption{$y=x$}
%  %       \label{fig:y equals x}
%    \end{subfigure}
%%     \hfill
%     \begin{subfigure}[b]{0.4\textwidth}
%         \centering
%         \includegraphics[width=\textwidth]{d3GauSpSignQQplot}
%       %  \caption{$y=3sinx$}
%      %   \label{fig:three sin x}
%  \end{subfigure}   %\vspace{-1mm}  
%     \hfill
%     \begin{subfigure}[b]{0.3\textwidth}
%         \centering
%         \includegraphics[width=\textwidth]{graph3}
%         \caption{$y=5/x$}
%         \label{fig:five over x}
%     \end{subfigure}
        \caption{QQ plots of the empirical distributions (across $N=1000$ replications) of the test statistics in  Table~\ref{Tab.white}  against their asymptotically chi-square null distributions: vdW   (left), Gaussian, Spearman, and sign test scores   (right).}
       \label{QQFig}
\end{figure}

Rejection frequencies at $5\%$ nominal level are reported  in Table~\ref{Tab.white}  for  the Gaussian, vdW, Spearman and sign tests based on the asymptotic critical values provided in Sections~\ref{Sec.theta0} and~\ref{Sec.Gauss.theta0} but also for the vdW and Spearman tests based on bias-corrected critical values as described below. %Leaving aside the latter, %(postponed to Section~\ref{biasSec}),
  Inspection of  Table~\ref{Tab.white} 
   reveals that,  under  normal innovations, the rejection frequencies under the null hypothesis of the Gaussian    test are very close to the  nominal  $5\%$ size. But the vdW and Spearman   tests are quite below that  nominal size, indicating  a severe bias. This bias is confirmed by   Figure~\ref{QQFig}, which shows the QQ plots  of the values of the vdW test statistic across the~$N=1000$ replications with Gaussian $f$ (but $f$ has no impact here, as the ranks are distribution-free under the null).  The tails of the vdW statistic, for instance, very clearly do not match those of the chi-square they are converging to. The reason for this is    the relatively slow convergence   of the distribution of center-outward rank statistics:~$n_R=20$ for $n=800$ is still a rather small $n_R$ value, which explains the phenomenon. Note that the same QQ plot  for the sign test score statistic, which does not depend  on the ranks, are more satisfactory.\vspace{-1mm}

Now, that bias of rank-based tests  is  easily   corrected by considering {\it permutational} critical values instead of the asymptotic ones.  The latter are obtained by generating $M$   permutations of the gridpoints ${\mathfrak G}\n$, and taking   the 0.95 quantile of the resulting $M$ values of the test statistic   as a critical value instead of the chi-square quantile of order 0.95.  The resulting tests have asymptotic size 5 \%  under the null.  Their rejection frequencies  are reported in Table 1  under the label ``bias-corrected,'' for the vdW and Spearman tests. %Such corrections are not possible with parametric test statistics. 

\begin{table}[!ht]
\caption{Numbers of under-, correct, and over-identification for the VAR($1$) model under various innovation densities and additive outlier schemes, of the Gaussian, vdW, bias-corrected vdW, Spearman, bias-corrected Spearman, \color{black}and sign tests out of~$N=1000$ replications; the sample size is $n = 800$ ($n_R = 20$, $n_S = 40$, and~$n_0 = 0$);   results of the sign tests with~$n_S=n$ are also included; permutational crirical values are based on $M=1000$ random permutations.\color{black}}  \label{Tab.VAR1}\vspace{1mm}
\centering
\footnotesize
\begin{tabular}{lcccc | lcccc}
\hline\hline

             &          & \multicolumn{3}{c|}{Order identification} &                   &          & \multicolumn{3}{c}{Order identification} \\ \cline{3-5} \cline{8-10}

$\quad\ f$ &\hspace{-20mm} Test         & 0         & 1          & $\geq 2$        & $\quad\ f$      &\hspace{-20mm} Test     & 0         & 1          & $\geq 2$        \\ \hline

\underline{Normal}&    &  &  & & \underline{Mixture}    & &  &   &   \\

      &\hspace{-20mm}  Gaussian &\hspace{-2mm}  0 & 953  & 47 &                    & \hspace{-20mm}  Gaussian & 0 & 950  & 50  \\

             &\hspace{-20mm}  vdW &\hspace{-2mm} 0 & 990  & 10 &                             & \hspace{-20mm}  vdW      & 0 & 991  & 9   \\

              &\hspace{-20mm}{bias-corrected vdW} & {0} & {979} & {21} &  & \hspace{-20mm}{bias-corrected vdW} & {0} & {974} & {26} \\

%             &         (bias-corrected)                               &                        &                        &                        &                   &                                       (bias-corrected) &                        &                        &               \\

                            & \hspace{-20mm} Spearman & 0 & 984  & 16 &                             & \hspace{-20mm} Spearman & 0 & 980  & 20  \\

                            & \hspace{-20mm}{bias-corrected Spearman} & {0} & {982} & {18} &  & \hspace{-20mm}{bias-corrected Spearman} & {0} & {971} & {29} \\

%                                   &         (bias-corrected)                               &                        &                        &                        &                   &                                       (bias-corrected) &                        &                        &               \\

                             & \hspace{-20mm} Sign  ($n_S=40$)   & 0 & 977  & 23 &                              & \hspace{-20mm} Sign  ($n_S=40$)    & 0 & 972  & 28  \\

                             & \hspace{-20mm} Sign  ($n_S =800$)  & {0} & {979} & {21} &   & \hspace{-20mm} Sign  ($n_S=800$)  & {0} & {979} & {21} \\ \hline

  %                           &          ($n_R =1$)                              &                        &                        &                        &                   &                                       ($n_R =1$) &                        &                        &               \\   \hline

\underline{$t_3$}          &    &  &  & & \underline{Skew-$t_3$}    & &  &   &   \\

&\hspace{-20mm}  Gaussian & 0 & 908  & 92 &                   & \hspace{-20mm}  Gaussian & 0 & 855  & 145 \\

                            & \hspace{-20mm} vdW      & 0 & 985  & 15 &                             & \hspace{-20mm} vdW      & 0 & 987  & 13  \\

                          & \hspace{-20mm}{bias-corrected vdW} & {0} & {961} & {39} &  & \hspace{-20mm}{bias-corrected vdW} & {0} & {970} & {30} \\

 %           &         (bias-corrected)                               &                        &                        &                        &                   &                                       (bias-corrected) &                        &                        &               \\

                            &  \hspace{-20mm} Spearman & 0 & 973  & 27 &                             &  \hspace{-20mm} Spearman & 0 & 970  & 30  \\

                            & \hspace{-20mm} {bias-corrected Spearman} & {0} & {961} & {39} &  &  \hspace{-20mm} {bias-corrected Spearman} & {0} & {968} & {32} \\

 %                             &         (bias-corrected)                               &                        &                        &                        &                   &                                       (bias-corrected) &                        &                        &               \\

                             & \hspace{-20mm} Sign   ($_S=40$)   & 0 & 969  & 31 &                              & \hspace{-20mm} Sign   ($n_S=40$)   & 0 & 981  & 19  \\

                             & \hspace{-20mm} Sign  ($n_S =800$)  & {0} & {968} & {32} &   & \hspace{-20mm} Sign  ($n_S =800$)  & {0} & {978} & {22} \\

%

%                              &          ($n_R =1$)                              &                        &                        &                        &                   &                                       ($n_R =1$) &                        &                        &               \\  

                            \hline

\underline{AOs}  (${\mbf s} = (6, 6)^\prime$) &  &  &  &  &  \underline{AOs}  (${\mbf s} = (9, 9)^\prime$) &  &  &     &   \\

 &  \hspace{-20mm}  Gaussian & 0 & 933  & 67 &   &  \hspace{-20mm}  Gaussian & 0 & 878  & 122 \\

                      &  \hspace{-20mm}  vdW      & 0 & 1000 & 0  &                           &  \hspace{-20mm}  vdW      & 0 & 1000 & 0   \\

 &  \hspace{-20mm}  {bias-corrected vdW} & {0} & {1000} & {0} &  &  \hspace{-20mm}  {bias-corrected vdW} & {0} & {1000} & {0} \\

%                  &         (bias-corrected)                               &                        &                        &                        &                   &                                       (bias-corrected) &                        &                        &               \\

                            & \hspace{-20mm}   Spearman & 0 & 1000 & 0  &                             &  \hspace{-20mm}  Spearman & 0 & 1000 & 0   \\

                            &  \hspace{-20mm}  {bias-corrected Spearman} & {0} & {1000} & {0} &  &  \hspace{-20mm}  {bias-corrected Spearman} & {0} & {1000} & {0} \\

%                                         &         (bias-corrected)                               &                        &                        &                        &                   &                                       (bias-corrected) &                        &                        &               \\

                             & \hspace{-20mm} Sign   ($n_S=40$)   & 0 & 998  & 2  &                              & \hspace{-20mm} Sign  ($n_S=40$)    & 0 & 987  & 13           \\

                              & \hspace{-20mm} Sign  ($n_S =800$)  & {0} & {998} & {2} &   & \hspace{-20mm} Sign  ($n_S =800$)  & {0} & {990} & {10} \\

  %                           &          ($n_R =1$)                              &                        &                        &                        &                   &                                       ($n_R =1$) &                        &                        &               \\  

                             \hline \hline

\end{tabular}\vspace{-5mm}
\label{Tab 3}
\end{table}

The vdW and Spearman tests have greater power than the sign test and slightly less power than the Gaussian one under Gaussian  innovations. Under  $t_3$ innovations, however, the vdW  and Spearman tests, despite their  bias, both outperform  the Gaussian test.  The   bias-corrected vdW and Spearman tests have  correct size under the null and even higher powers;  under the mixture and skew-$t_3$ innovations, they outperform the Gaussian test by a landslide. 
% Under the mixture and skew-$t_3$ distributions, the superiority of the vdW and Spearman tests over the Gaussian one is overwhelming for   $\ell = 1$ (for $\ell =2$ ($d=2$) and  $\ell = 1.5$  ($d=3$), all tests yield good powers).  
The tests based on signs have been implemented with the same factorization (viz., $n_R = 20$, $n_S = 40$, and~$n_0 = 0$) as vdW and Spearman; the factorization ($n_R = 1$, $n_S = 800$, and~$n_0 =~\!0$) also has been considered, but provides little improvement. 
%---it nevertheless reduces the bias under most densities.
 Turning to robustness issues,   the resistance to additive outliers of the rank-based tests over the Gaussian one is extremely good, while the Gaussian test literally crashes, with exploding size under the null---the corresponding rejection frequencies under the alternative thus are meaningless. Despite their simplicity, the sign tests do extremely well (confirming univariate robustness results by \cite{Boldin12, Boldin13}) but remain   less powerful than the bias-corrected vdW tests. Bias-corrected tests uniformly outperform their uncorrected asymptotic counterparts. \vspace{-2mm}

%%%%%%%%%%%%%%%%%%%%%%%%%%%%%%%%%%%%%

\subsection{VAR order identification\vspace{-1mm}} \label{VARorder}
Turning to  VAR order identification, let us consider  the  bivariate VAR(1) model with autoregression matrix coefficient $\text{vec}(\A) = 
(0.30, -0.06, 0.12, 0.24)^\prime$, the same innovation densities and the same contamination schemes as in Section~\ref{sec.NumWhite}, and sample size $n = 800$.   The sequential method of  Sections~\ref{Sec.testorder}  and \ref{GIdSec} was applied to $N=1000$ replications thereof. 

Each step consists in testing (at 5\% nominal level) a VAR($p$) against a VAR($p+1$), which requires the estimation of the null VAR($p$) model;  throughout, this estimation is performed via the center-outward R-estimator based on  vdW scores as  proposed by \cite{HLL2019}, where it is shown to be root-$n$ consistent. 
% and was shown in \cite{HLL2019} to have good performance under all distributions we consider here. 
 The numbers of under-, correct, and over-identification for the Gaussian,   (corrected) vdW,  (corrected) Spearman, and  (corrected) sign tests  are reported in Table~\ref{Tab.VAR1}.

All tests are rejecting the hypothesis of white noise, irrespective of the underlying innovation density or the presence of outliers: no under-identifications, thus. Also, the rank-based procedures yield uniformly lower over-identification rates than the Gaussian one. Under the $t_3$ distribution, the Gaussian procedure yields an almost 10\% over-identification rate, indicating that its performance is badly  affected by heavy tails. The rank-based procedures, on the other hand, do much better under heavy tails, skewness, and additive outliers.

\subsection{A real data example}\label{Sec:emp}

To illustrate the applicability of our VAR order selection procedure, we consider the joint behaviour of 
two macroeconomic time series in the U.S. economy. Specifically, we downloaded  the M1SL Money Stock 
 and  the Current Price Index All Urban Customers for All Items (labeled as M1SL and CPIAUCSL, respectively) time series  from the Federal Reserve Bank of Saint Louis economic data server (see \url{https://
fred.stlouisfed.org}). The M1SL includes funds that are 
readily accessible for spending and  represent the most liquid portions of the money supply quickly convertible into cash. 
%For instance, 
%it contains physical currency and coins, demand deposits, traveler's checks of nonbank issuers, demand deposits, and other 
%checkable deposits. 
The CPIAUCSL is a measure of the average monthly change in the price for goods and services paid by urban consumers between any 
two time periods; we refer to U.S. Bureau of Labor Statistics for further technical details. In our statistical analysis, we consider monthly 
records for the period 01-Jan-1960 to 01-Sept-2020.  For additional info, we refer to the Board of Governors of the Federal Reserve System. 
The joint semiparametric modeling of these
%time
 series can be of help for monetary policy decisions.\vspace{-3mm}
 
 \begin{center}
\begin{figure}[!ht]
\caption{Top  panels: M1SL series in levels (left panel, in Billions of Dollars);  CPIAUCSL series in levels (right panel, the Index $1982-1984=100$). Bottom panels: same   series, differentiated (first-order).\vspace{-2mm}}
\begin{center}
\includegraphics[width=0.95\textwidth, height=0.65\textwidth]{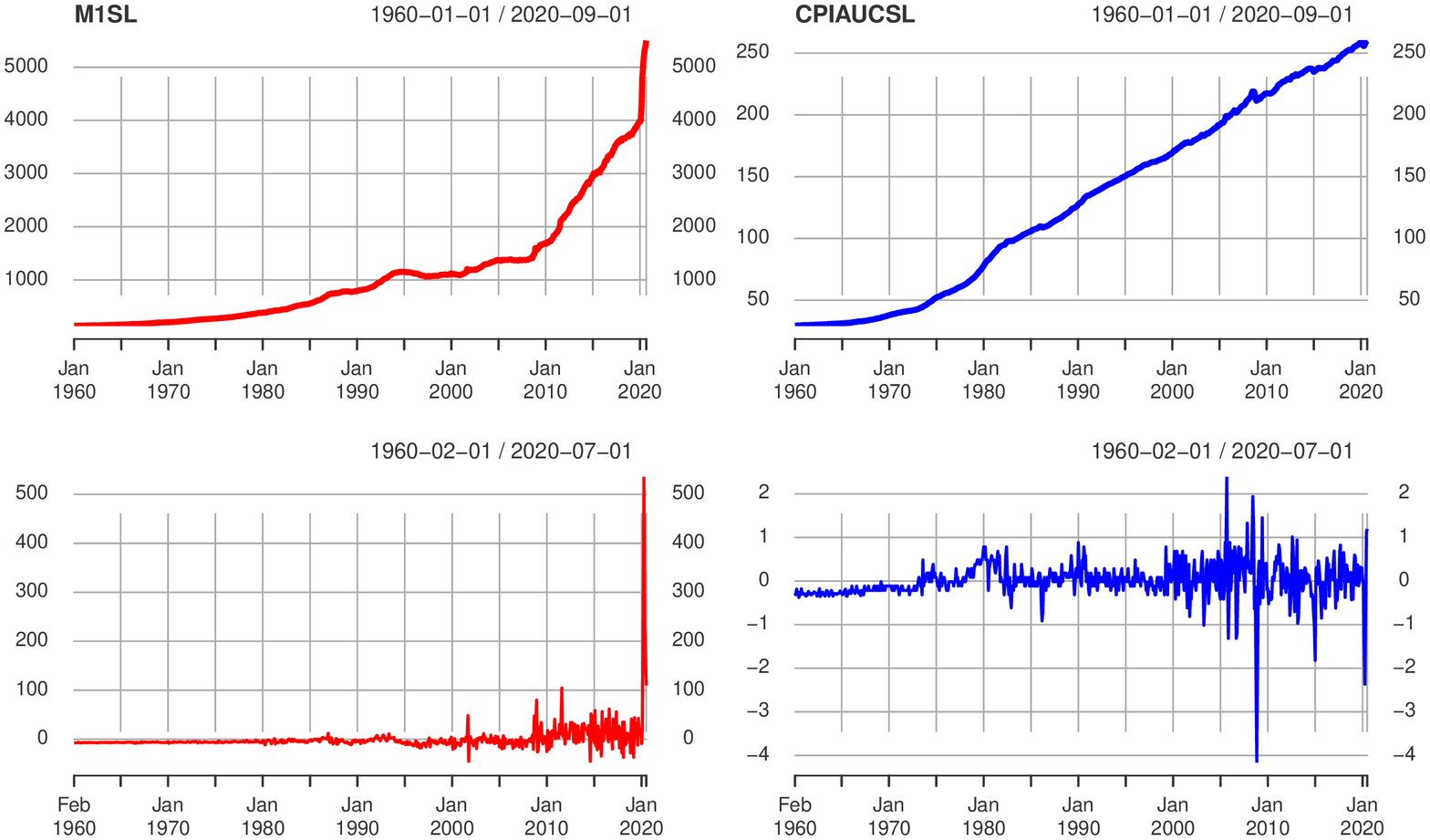}
\end{center}
\label{Fig. all}
\end{figure}
\end{center}
\begin{table}[!hb]
\caption{Values of test statistics and $p$-values (in parenthesis) of the Gaussian, vdW, Spearman, sign ($n_S = 22$), and sign ($n_S = 726$) tests, along with their $p$-values  (in  brackets) and (for vdW and Spearman) bias-corrected $p$-values  (in square brackets)      under the null hypotheses of  a VAR($p_0$) model ($p_0=0, ..., 7$).}\label{testReal}\vspace{-3mm} 
\centering
\scriptsize
\begin{tabular}{cccccc}
\hline
$p_0$ & Gaussian         & vdW              & Spearman            & Sign ($n_S = 22$) & Sign ($n_S = 726$) \\ \hline
$0$ & 25.05 (4.91$\times 10^{-5}$) & 224.67 (0) [0]     & 238.71 (0)  [0]       & 356.90 (0)      & 342.58 (0)         \\ 
$1$ & 25.83 (3.42$\times 10^{-5}$) & 16.77 (0.002) [0.001]   & 15.42 (0.004) [0.001]     & 26.16 (2.94$\times 10^{-5}$)  & 59.14 (4.40$\times 10^{-12}$)   \\
$2$ & 2.30 (0.681)    & 62.35 (9.30$\times 10^{-13}$) [0] & 71.57 (1.07$\times 10^{-14}$) [0]   & 16.56 (0.002)    & 56.74 (1.40$\times 10^{-11}$)   \\
$3$ &                  & 48.59 (7.18$\times 10^{-10}$) [0] & 56.90 (1.30$\times 10^{-11}$)  [0]  & 47.99 (9.49$\times 10^{-10}$)  & 55.64 (2.38$\times 10^{-11}$)   \\
$4$ &                  & 10.74 (0.030) [0.016]  & 13.14 (0.011)     [0.010] & 26.58 (2.41$\times 10^{-5}$)  & 12.94 (0.012)     \\
$5$ &                  & 19.02 (7.81$\times 10^{-5}$) [0] & 21.37 (2.67$\times 10^{-5}$) [0] & 16.57 (0.002)    & 40.57 (3.29$\times 10^{-8}$)   \\
$6$ &                  & 7.58 (0.108)  [0.068]  & 24.14 (7.49$\times 10^{-5}$)  [0.005]  & 32.58 (1.45$\times 10^{-6}$)  & 26.22 (2.85$\times 10^{-5}$)   \\
$7$ &                  &                  & 4.27 (0.371)  [0.322]     & 8.60 (0.072)     & 8.24 (0.083)    
  \\  \hline
\end{tabular}\vspace{-0mm}
\label{Tab: real}
\end{table}

$\,$\vspace{-8mm}

To start with, in the top panels of Figure~\ref{Fig. all}, we display the seasonally adjusted time series. 
Even a visual inspection reveals some interesting characteristics. A first clearly visible aspect is that both time series display a trend in 
time. To deal with this, we differentiate the series and plot the resulting outputs in the bottom panels of Figure~\ref{Fig. all}. We notice 
that the differentiated  series seem to be 
correlated and display common movements (mainly in opposite directions). A second noticeable aspect is that, starting from 2005, the 
trajectories of the differentiated series look increasingly 
asymmetric and spiky. %with respect to the zero line,
% display  a number of spikes (both below and above zero).
 The period 
 March-September-2020 reveals a cluster of outlying values, most likely due to  COVID-19-related policy decisions.
\begin{figure}[!t]
\caption{Plots of the fitted residuals of the M1SL (upper left panel) and CPIAUCSL (upper right  panel) series. Scatterplots of the M1SL residuals (bottom left panel); the same scatterplot (bottom right) with center-outward quantile contour of order 0.94 (rank 32) (see \cite{Hallinetal2020}) for details.\vspace{-9mm}}
\center
\begin{tabular}{cc}
\includegraphics[width=0.34\textwidth, height=0.35\textwidth]{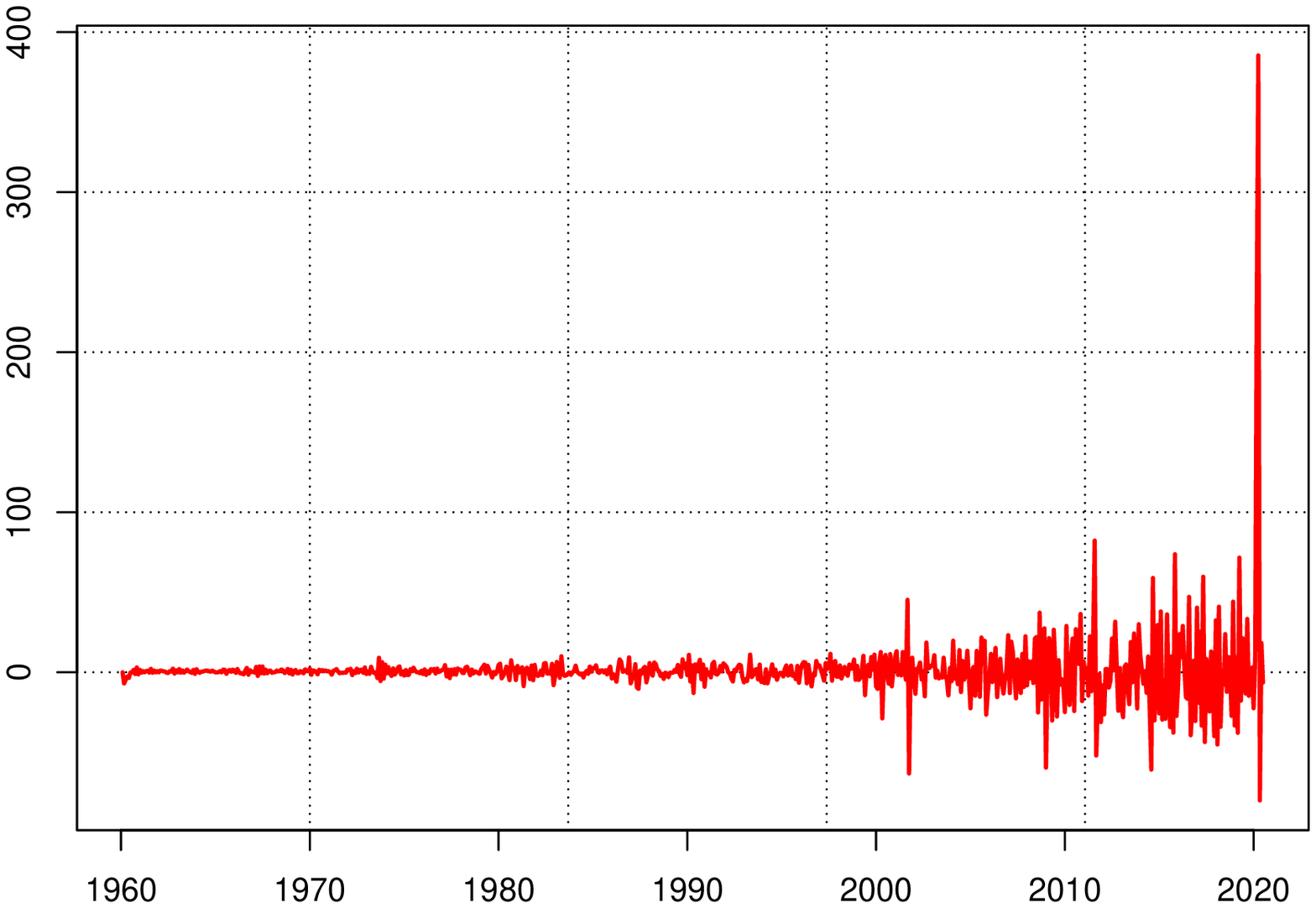}&
\includegraphics[width=0.34\textwidth, height=0.35\textwidth]{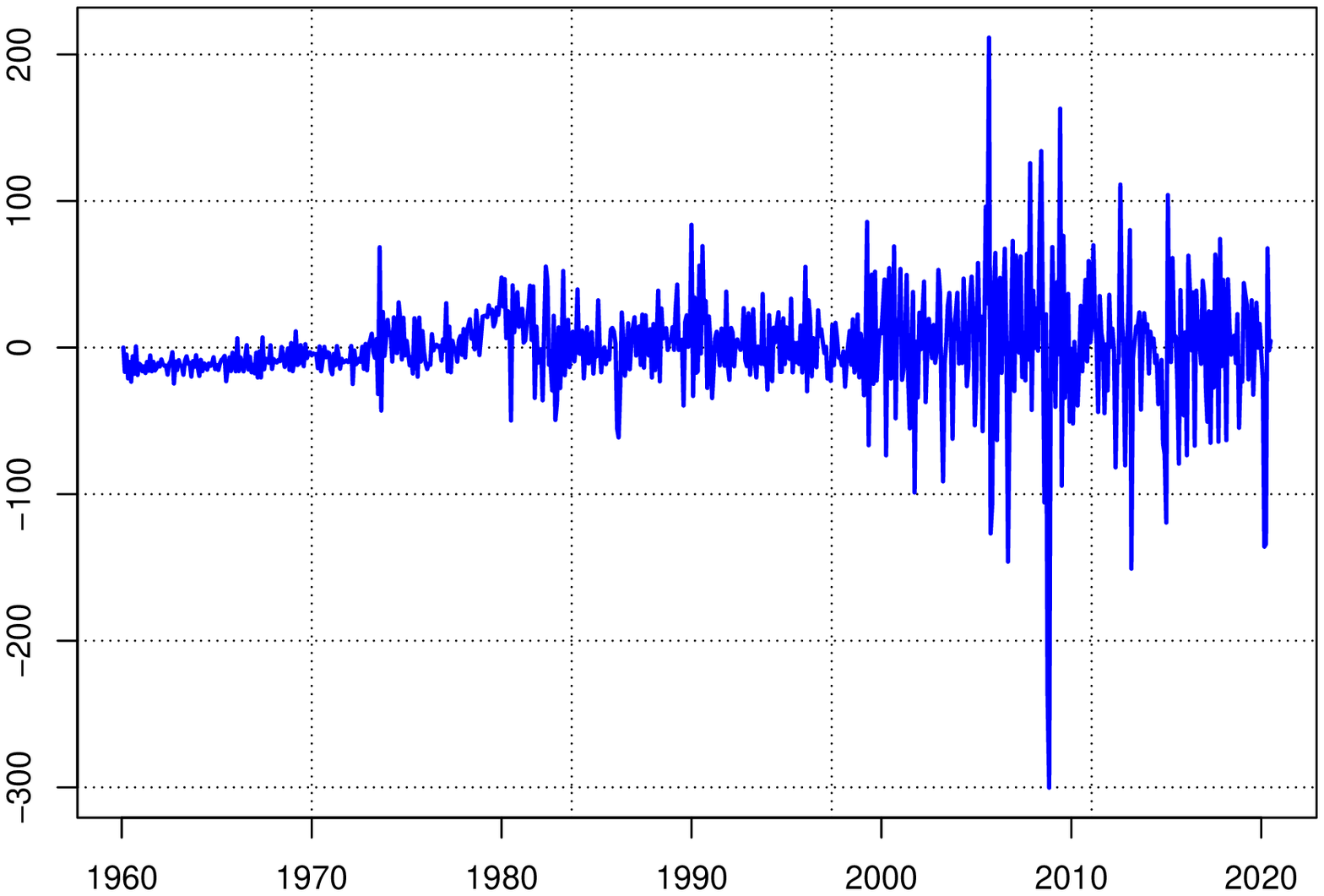}\vspace{-9mm}
\\ 
%\end{tabular}
%\begin{tabular}{c}
%\hspace{2.5cm} 
\includegraphics[width=0.34\textwidth, height=0.35\textwidth]{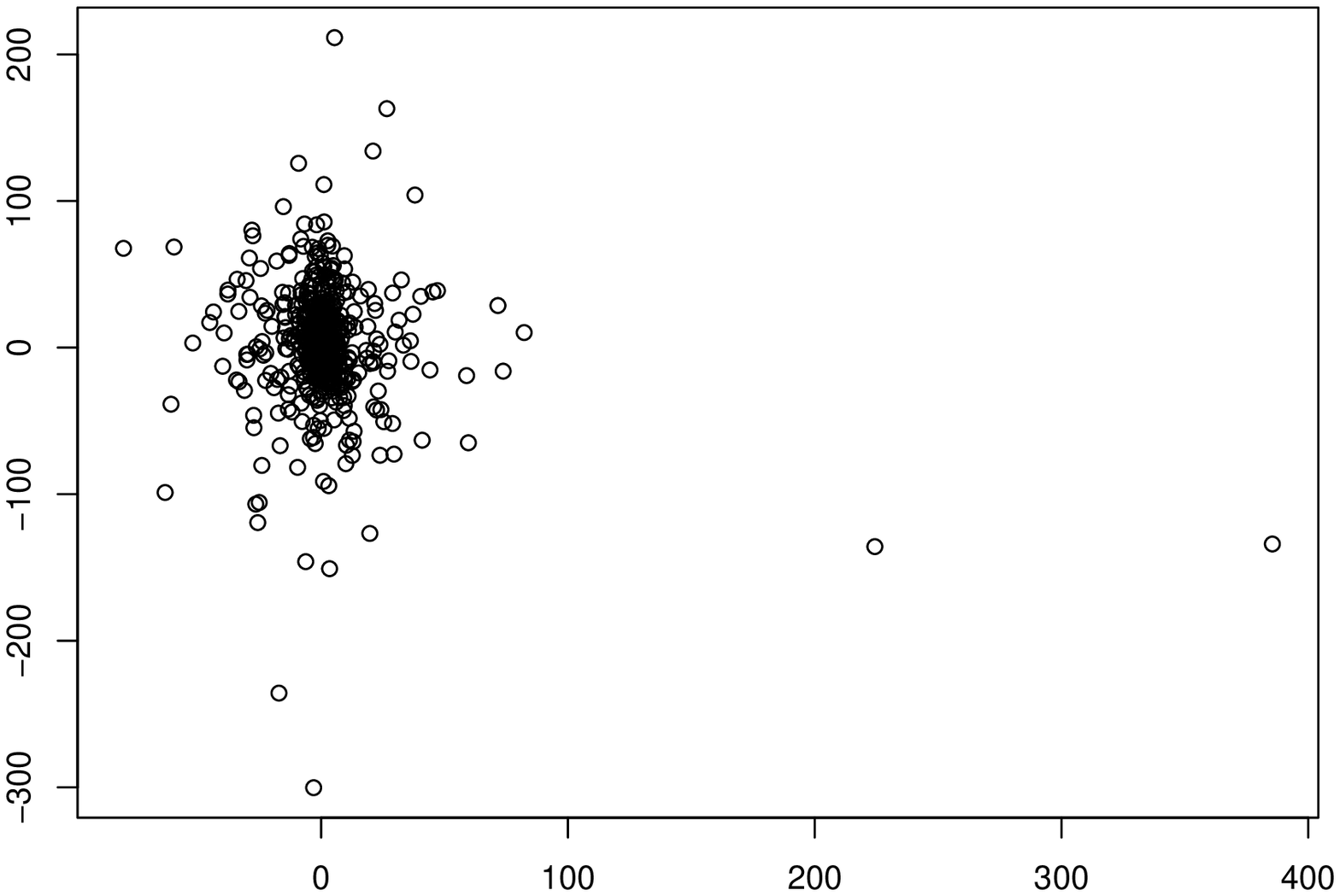}&\hspace{-3mm}\raisebox{8mm}{
\includegraphics[width=0.24\textwidth, height=0.28\textwidth]{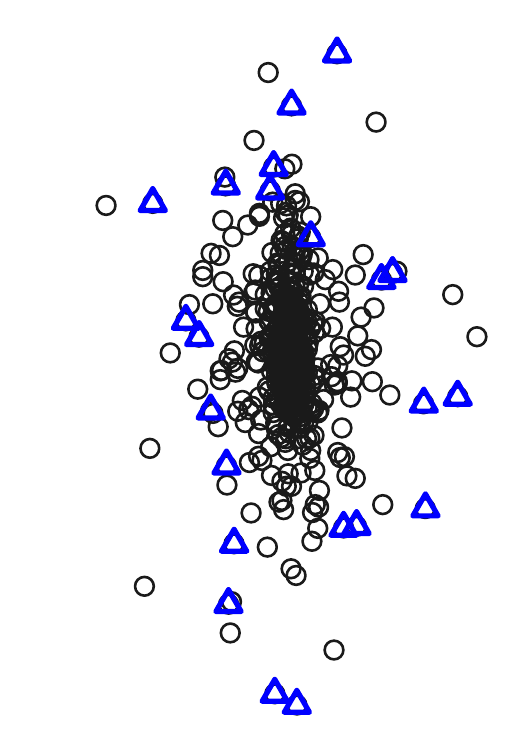}
}
\end{tabular}\vspace{-11mm}
\label{PlotResidual}
\end{figure}
%\end{center}

%Combining the above considerations, we believe that
A VAR model of appropriate order is likely to provide an adequately describe the joint behavior of the M1SL and CPIAUCSL  
series. However, the plots in Figure~\ref{Fig. all} suggest possible heavy tails and the presence, during the COVID-19 pandemic period,  of large outliers. Gaussian procedures, therefore, might be inappropriate in view of the Monte Carlo analyses of Sections~\ref{sec.NumWhite} and \ref{VARorder}. 
% that we need to resort on inference procedures which are able, by design, to 
%handle a sequence of extreme observations, such as those recorded during the COVID-19 pandemic. As a result, we have serious concerns about 
%the suitability of Gaussian tests. Indeed, our 
%experience gained via  Monte Carlo studies (see Tables 1 and 2) indicates that the selection of the VAR order via 
%Gaussian tests can entail the underestimation of the number of lags. This problem is particularly serious in the presence of heavy-tails. 
%Therefore, we believe that the use of our rank-based tests has to be 
%preferred to the routinely-applied Gaussian test. 
To investigate that point,  besides  the traditional Gaussian  one, we also ran, on the same series, the rank-based order-identification procedures described in Section~\ref{VARorder}. The results are shown in Table~\ref{testReal}. %\vspace{-2mm}

%Each step consists in testing (at 5\% nominal level) a VAR($p$) against a VAR($p+1$).

Inspection of Table~\ref{testReal} reveals a sharp contrast between the conclusions of the Gaussian and rank-nased methods. The Gaussian procedure indeed very clearly  selects  (at  nominal level~$5\%$) a VAR(2) model while all   rank-based procedures, whether bias-corrected (with~$M = 1000$ permutations)  or not, agree on a larger number of 6-7 lags. 
%Specifically, the vdW test indicates that a VAR(6) model is suitable, while the Spearman and the sign (with different values of $n_S$)  tests suggest a VAR(7).
% In Table~\ref{testReal}, we also report $p$-values of the bias-corrected vdW and Spearman (with $M = 1000$ permutations) in square brackets. In line with their non bias-corrected counterparts, the bias-corrected vdW and Spearman also select order 6 and order 7, respectively. 
%\begin{center}

Our educated guess is that the combination of skewness, kurtosis, and outliers  are blurring the conclusions of Gaussian tests;    the graphical diagnostics in Figure~\ref{PlotResidual} bring some evidence in favor of that guess.  Figure~\ref{PlotResidual} displays   the residuals of a VAR(7) fit of the M1SL  and CPIAUCSL   series. The  panels show that the fitted residuals of each time series are centered about zero but   asymmetrically scattered, apparently with heavy tails. The scatterplots (bottom panels)   confirms the presence of skewness and kurtosis; some isolated, large outlying values are clearly visible, which supports our conjecture. These features of the underlying distribution obviously have a significant impact on the Gaussian procedure whereas our rank-based approach is much less affected. \vspace{-2mm}

\section{Conclusion\vspace{-1mm}} This paper introduces   rank-based tests for VAR models with unspecified innovation density. Based on the residual center-outward ranks and signs recently proposed by \cite{Chernetal17} and \cite[where they are shown to be essentially maximal ancillary]{Hallinetal2020}, they constitute  the testing counterpart of the R-estimators proposed in \cite{HLL2019}. When testing a VAR with  specified parameter $\bth_0$, they are  fully distribution-free; when $\bth_0$ remains unspecified under the null, they are {\it strongly} asymptotically distribution-free in the sense of being asymptotically equal to a fully-distribution-free test. When based on appropriate score functions, they achieve parametric efficiency at selected reference densities. Monte Carlo experiments and an empirical example demonstrate the excellent performance of the proposed tests, which significantly outperform the traditional pseudo-Gaussian methods under skew, heavy-tailed, and contaminated innovations. %\vspace{-4mm}

\newpage

\appendix \section*{Appendix A: Algebraic preparation for LAN}\label{algebrasec}

Some linear difference equation algebra is required in order to obtain the explicit form of the central sequences and information matrices in the LAN result of Section~\ref{secLAN}---to be used in the construction of our center-outward rank-based test statistics.

 Denote by $\G_u$, $u \in \Z$ the {\it Green's matrices} associated with the difference operator~$\A_{\bth_0}(L)$, with~$\bth_0\in\bTH_0\subset\bTH_1\subset \mathbb{R}^{p_1d^2}$ for some  $p_0 <p_1$; those matrices    are defined as the solutions of  the homogeneous   linear recursions
$$\A_{\bth_0}(L)\G_u= \G_u - \sum_{i=1}^{p_0} \A_i   \G_{u-i} = {\bf 0}, \quad u\in\mathbb{Z}$$
 with initial values $\I_d, {\bf 0},\ldots ,{\bf 0}$ at $u=0,-1,\ldots ,-p_0+1$.

For $\bth_0\in\bTH$, the decrease of $\{\Vert \G_u\Vert , u\in \N\}$ as $u\to\infty$ is exponential. Specifically, there exists some $\varepsilon > 0$ (depending on the smallest  root of the %rdeterminantal
 equation~$\text{\rm det} \left(\A_{\bth} (z)
\right)~\!=~\!0$)  such that~$\Vert \G_u\Vert  (1+\varepsilon)^u$ converges to $0$ as $u \rightarrow \infty$.
Denoting  by~$\X_{-p_0 + 1}, \ldots , \X_{0}$ the initial values (which are typically unobservable) of a solution of the homogenous difference equation~$\A_{\bth_0}(L)\X_t={\bf 0}$, this exponential decrease
 ensures
 that these  initial values
  have no asymptotic  influence on   asymptotic results and, therefore,   safely can be set to zero in the sequel. This allows us to invert the autoregressive polynomial;  the Green matrices $\G_u$  then are the matrix coefficients\footnote{In the econometric terminology, the {\it impulse response matrix  coefficients}.} of the inverted operator $\A_{\bth_0} (L)^{-1}$:
$$\A_{\bth_0} (z)^{-1}=\sum_{u=0}^\infty \G_u z^u := \left(\I_d -  \sum_{i=1}^{p_0} \A_i z^i\right)^{-1}, \quad z \in \C, |z| < 1.$$

More generally, associated with an arbitrary $d$-dimensional linear difference ope-\linebreak  rator~$\CC (L) := \sum_{i=0}^\infty \CC_i L^i$ (this of course includes operators of finite order $s$), define, for any integers  $ v \leq u$, the $d^2 u \times d^2 v$ matrices (the form of which is motivated by the vec$\A_u$ form under which the VAR matrix coefficients  $\A_u$ enter $\bth_0$ and the $d^2v$-dimension nature of the perturbations $\btau$ on the LAN property )
\[
\CC_{u, v} :=
\begin{bmatrix}
\CC_0 \otimes \I_d & \0 & \ldots  & \0 \\
\CC_1 \otimes \I_d & \CC_0 \otimes \I_d & \ldots  & \0 \\
\vdots & & \ddots & \vdots \\
\CC_{v-1} \otimes \I_d & \CC_{v-2} \otimes \I_d & \ldots  & \CC_0 \otimes \I_d \\
\vdots & & & \vdots \\
\CC_{u-1} \otimes \I_d & \CC_{u-2} \otimes \I_d & \ldots  & \CC_{u-v} \otimes \I_d
\end{bmatrix}.
\]
%and write $\CC_u$ for $\CC_{u, u}$. 
 With this notation, $\G_{u, u}$ (associated with $\CC (L)=\A_{\bth_0} (L)^{-1}$) is the inverse of~$\A_{u, u}$ (associated with $\CC (L)=\A_{\bth_0} (L)$).   Denoting by $\CC_{u, v}^{\prime}$   the matrices  associated with the transposed operator $\CC^\prime (L) := \sum_{i=0}^\infty \CC^\prime_i L^i$, we also have that~$\G_{u, u}^{\prime} = (\A_{u, u}^{\prime})^{-1}$. Then,   for any $\bth_0\in\bTH_0$ and any~$p_1>p_0$,   the $d^2 p_1 \times d^2 p_1$ matrix 
\begin{equation*}\label{defM}
\M_{\bth_0} := \G^{\prime}_{p_1, p_1};\tag{A.1}
\end{equation*}
 is of full rank.

Still associated with $\A_{\bth_0}(L)$, consider the operator 
%\footnote{This operator   and most quantities defined below depend on $\bth_0$; for the sake of  simplicity, however, we are dropping this reference to $\bth_0$ in the notation.}
 $\D_{\bth_0} (L) := \I_d + \sum_{i=1}^{p_0} \D_i L^i$ where
\[
\begin{bmatrix}
\D_1^\prime \\
\vdots \\
\D_{p_0}^\prime
\end{bmatrix}
:= -
\begin{bmatrix}
\G_0 & \G_{-1} & \ldots  & \G_{-p_0+1} \\
\G_{1} & \G_{0} & \ldots  & \G_{-p_0+2} \\
\vdots & & \ddots & \vdots \\
\G_{p_0-1} & \G_{p_0-2} & \ldots  & \G_{0} \\
\end{bmatrix}
^{-1}
\begin{bmatrix}
\G_{1} \\
\vdots \\
\G_{p_0}
\end{bmatrix}
\]
(recall that $\G_{-1} = \G_{-2} = \cdots = \G_{-p_0+1} = \0$).

 Let $\{\bpsi_t^{(1)}, \ldots , \bpsi_t^{(p_0)} \}$ be any set of $d \times d$ matrices forming a fundamental system of solutions of the homogeneous linear difference equation associated with $\D_{\bth_0} (L)$. Such a system can be obtained, for instance,  from the Green matrices of   $\D_{\bth_0} (L)$ (see, e.g., \citealt{Hallin1986}). Defining
\[
{\boldsymbol\Psi}_m (\bth_0) :=
\begin{bmatrix}
\bpsi_{p_1 - p_0 + 1}^{(1)} & \ldots  & \bpsi_{p_1 - p_0+ 1}^{(p_0)} \\
\bpsi_{p_1 - p_0+ 2}^{(1)} & \ldots  & \bpsi_{p_1 - p_0 + 2}^{(p_0)} \\
\vdots & & \vdots \\
\bpsi_{m}^{(1)} & \ldots  & \bpsi_{m}^{(p_0)}
\end{bmatrix}
\otimes \I_d   \qquad (m > p_1 - p_0),
\]
the {\it Casorati matrix} %$\mathbf{C}_{\bpsi}$
 associated with $\D_{\bth_0} (L)$ is ${\boldsymbol\Psi}_{p_1}$, which has full rank. Finally, let 
% \color{red}What is the connection between $\D (L)$ (hence the matrices below) and our problem?? What is $p_1$ when $\D (L)$ is of order $p_0$??  \color{black} \color{blue} Also, $\D (L)$, $\P_{\bth_0}$ and $\Q^{(n)}_{\bth_0}$ are from Hallin and Paindaveine (2004, The Ann. of Stats.). $\P_{\bth_0}$ and $\Q^{(n)}_{\bth_0}$ are used in the central sequence \color{black}
\begin{equation*}
\P_{\bth_0} :=
\begin{bmatrix}
\I_{d^2(p_1 - p_0)} & \0 \\
\0 & {\boldsymbol\Psi}_{p_1}^{-1}(\bth_0)
\end{bmatrix}
\quad\text{and}\quad
\Q^{(n)}_{\bth_0} :=
\begin{bmatrix}
\I_{d^2(p_1 - p_0)} & \0 \\
\0 & {\boldsymbol\Psi}_{n-1}(\bth_0)
\end{bmatrix}
.\label{defPQ}\tag{A.2}
\end{equation*}
Note that while $\M_{\bth_0}$, $\P_{\bth_0}$, and  $\Q_{\bth_0}^{(n)}$ depend on the choice of the fundamental system $\{\bpsi_t^{(1)}, \ldots , \bpsi_t^{(p_0)} \}$,  the product $\M_{\bth_0}^\prime  \P_{\bth_0}^\prime \Q_{\bth_0}^{(n)\prime} $ appearing in the definition of the central sequence \eqref{Delta} 
 does not.  We refer to \cite{HP04} for details.

\appendix \section*{Appendix B:  
%\subsection{
Proofs of Proposition~\ref{Prop.bar.til.Gam} and Lemma~\ref{lem.xi}}\label{Sec3}

\renewcommand\theequation{B\arabic{equation}}
\setcounter{equation}{0}

Let $\a_{t, s}^{(n)} := \a (\F^{(n)}_{\pm, t}, \F^{(n)}_{\pm, s})$. The following   Lemmas 
%~\ref{a21n} and Lemma~\ref{a1234n}
 will be used to prove Proposition~\ref{Prop.bar.til.Gam}; throughout this section, ${\rm E}$ stands for ${\rm E}_{\bth_0 ;f}$ ($f\in{\cal F}_d$).
 
 \begin{customlemma}{B1}\label{a21n}
%\begin{lemma}\label{a21n}
$%\displaystyle{
{\rm E}(\a_{2, 1}^{(n)}) = \frac{(n-2)(n-3)}{n(n-1)} {\rm E}(\a_{2, 1}^{(n)}| \F^{(n)}_{\pm, 4}, \F^{(n)}_{\pm, 3}) + \frac{1}{n(n-1)} \b^{(n)}_{4, 3}
%}
$,
where
\begin{align*}
\b^{(n)}_{4, 3} :=& (n-2)\Big[  {\rm E}(\a_{4, 1}^{(n)}| \F^{(n)}_{\pm, 4}, \F^{(n)}_{\pm, 3}) + {\rm E}(\a_{3, 1}^{(n)}| \F^{(n)}_{\pm, 4}, \F^{(n)}_{\pm, 3}) + {\rm E}(\a_{2, 4}^{(n)}| \F^{(n)}_{\pm, 4}, \F^{(n)}_{\pm, 3}) \\
& + {\rm E}(\a_{2, 3}^{(n)}| \F^{(n)}_{\pm, 4}, \F^{(n)}_{\pm, 3}) \Big] +  {\rm E}(\a_{3, 4}^{(n)}| \F^{(n)}_{\pm, 4}, \F^{(n)}_{\pm, 3}) + {\rm E}(\a_{4, 3}^{(n)}| \F^{(n)}_{\pm, 4}, \F^{(n)}_{\pm, 3}).
\end{align*}
\end{customlemma}
\begin{proof}
The proof follows from an enumeration of the $n(n-1)$ terms in $\underset{1\leq i_1 \neq i_2 \leq n}{\sum \sum} \a_{i_1, i_2}^{(n)}\vspace{1mm}$, which contains $(n-2)(n-3)$ terms with indices $i_1\neq i_2$ such that $\#[\{ i_1 \neq i_2\} \cap \{3, 4\}] =0$, $4(n-2)$ terms with indices $i_1\neq i_2$ such that $\#[\{ i_1 \neq i_2\} \cap \{3, 4\}] = 1$ and $2$ terms with indices $i_1\neq i_2$ such that $\#[ \{ i_1 \neq i_2\} \cap \{3, 4\}] = 2$.
\end{proof}

 \begin{customlemma}{B2}\label{a1234n}
$\left\vert {\rm E} (\a_{2, 1}^{(n)\prime}  \a_{4, 3}^{(n)}) \right\vert \leq \frac{n(n-1)}{(n-2)(n-3)} {\rm E} (\a_{2, 1}^{(n)\prime}) {\rm E} (\a_{2, 1}^{(n)}) + \frac{4(n-2) + 2}{(n-2)(n-3)} {\rm E} \left\Vert \a_{2, 1}^{(n)} \right\Vert^2.$ 
\end{customlemma}

\begin{proof}
Lemma~\ref{a21n} implies
\begin{align*}
 \left\vert {\rm E} (\a_{2, 1}^{(n)\prime}  \a_{4, 3}^{(n)}) \right\vert &= \left\vert  {\rm E}\left[  {\rm E}\left(\a_{2, 1}^{(n)\prime} \big| \F^{(n)}_{\pm, 4}, \F^{(n)}_{\pm, 3}\right) \a_{4, 3}^{(n)}\right] \right\vert \\
 &\leq \frac{n(n-1)}{(n-2)(n-3)} {\rm E} (\a_{2, 1}^{(n)\prime}) {\rm E} (\a_{4, 3}^{(n)}) + \frac{1}{(n-2)(n-3)} \left\vert {\rm E} (\b^{(n)\prime}_{4, 3} \a_{4, 3}^{(n)}) \right\vert \\ 
&=:\frac{n(n-1)}{(n-2)(n-3)} C\n_1 +  \frac{1}{(n-2)(n-3)}C\n_2,\quad\text{say.}  \end{align*}
The desired result then follows from the fact that $C\n_1 ={\rm E} (\a_{2, 1}^{(n)}) = {\rm E} (\a_{4, 3}^{(n)})$ and  the  inequality~$C\n_2\leq  (4(n-2)+2){\rm E} \left\Vert \a_{2, 1}^{(n)} \right\Vert^2$.   More specifically, note that
\begin{align*}
\left\vert {\rm E} (\b^{(n)\prime}_{4, 3} \a_{4, 3}^{(n)}) \right\vert  
&\leq  (n-2)\Big(  \left\vert {\rm E}(\a_{4, 1}^{(n)\prime}\a_{4, 3}^{(n)})  \right\vert  + \left\vert {\rm E}(\a_{3, 1}^{(n)\prime}\a_{4, 3}^{(n)})  \right\vert  + \left\vert {\rm E}(\a_{2, 4}^{(n)\prime}\a_{4, 3}^{(n)}) \right\vert +  \left\vert {\rm E}(\a_{2, 3}^{(n)\prime}\a_{4, 3}^{(n)})  \right\vert \Big)  \\
& \quad + \left\vert {\rm E}(\a_{3, 4}^{(n)\prime}\a_{4, 3}^{(n)})  \right\vert + \left\vert {\rm E}(\a_{4, 3}^{(n)\prime}\a_{4, 3}^{(n)}) \right\vert \\
&=: (n-2) \Big( C\n_{21} + C\n_{22} + C\n_{23} + C\n_{24} \Big) + C\n_{25} + C\n_{26},\quad\text{say.\vspace{-2mm}}
\end{align*}
It follows from %\vspace{-2mm}
$$\left\vert {\rm E}(\a_{4, 1}^{(n)\prime}\a_{4, 3}^{(n)})  \right\vert \leq ({\rm E}\Vert (\a_{4, 1}^{(n)}\Vert^2)^{1/2} ({\rm E}\Vert (\a_{4, 3}^{(n)}\Vert^2)^{1/2} = {\rm E} \left\Vert \a_{2, 1}^{(n)} \right\Vert^2$$
that $C\n_{21} \leq {\rm E} \left\Vert \a_{2, 1}^{(n)} \right\Vert^2$. For the same reason, we have $C\n_{2j} \leq {\rm E} \left\Vert \a_{2, 1}^{(n)} \right\Vert^2, \   j = 2, ..., 6$. The result then follows.
\end{proof}

\textbf{Proof of Proposition~\ref{Prop.bar.til.Gam}}. Due to  LAN, it is sufficient to prove the desired result  under~${\rm P}^{(n)}_{\bth_0 ;f}$ ($f\in{\cal F}_d$). Let
$\al_{t, s}^{(n)} \!:=\! (\bxi_{t, s}^{(n)} -  {\mbf \gamma}^{(n)}) -  (\bxi_{t, s} - {\mbf \gamma}),$
 where 
 \[
 \bxi_{t, s}^{(n)} \!:=\! \text{vec} (\a (\F^{(n)}_{\pm, t},  \F^{(n)}_{\pm, s})), \ \ \bxi_{t, s} \!:=\! \text{vec} (\J_1(\F_{\pm, t}) \J_2\pr(\F_{\pm, s})), \ \ {\mbf \gamma}^{(n)} \!\!:=\! \text{vec} (\m^{(n)}),   \text{ and } 
  {\mbf \gamma}  \!:=\! \text{vec} (\m).
  \] 
We have%\vspace{-2mm}
\begin{align*}
&{\rm E}\left\Vert (n-i)^{1/2}  \text{\rm vec} ( \tenq{\bGamma}_{i, \J_1, \J_2}^{(n)}(\bth_0) - \m^{(n)}  - \underline{\bGamma}_{i, \J_1, \J_2}^{(n)}(\bth_0) + \m) \right\Vert^2 = {\rm E} \Big\Vert (n-i)^{-1/2} \sum_{t=i+1}^n \al_{t, t-i}^{(n)} \Big\Vert^2 \\
%& \\
&= (n-i)^{-1} \!\!  \sum_{t=i+1}^n\! {\rm E} \Vert \al_{t, t-i}^{(n)} \Vert^2\! +2 (n-i)^{-1}\!\!  \underset{\substack{i+1 \leq t, s \leq n \\ t = s-i}}{\sum \sum}\! {\rm E} (\al_{t, t-i}^{(n)\prime}  \al_{s, s-i}^{(n)})    + (n-i)^{-1} \!\! \underset{\substack{i+1 \leq t, s \leq n \\ t\neq s, s-i  \\ s \neq t-i}} {\sum \sum} \!{\rm E} (\al_{t, t-i}^{(n)\prime}  \al_{s, s-i}^{(n)}) \\
&=  {\rm E} \Vert \al_{2, 1}^{(n)} \Vert^2 + \frac{2(n-2i)}{n-i}  {\rm E} (\al_{2, 1}^{(n)\prime}  \al_{3, 2}^{(n)}) + \frac{(n-i)^2-(n-i)-2(n-2i)}{n-i}   {\rm E} (\al_{2, 1}^{(n)\prime}  \al_{4, 3}^{(n)})  \\
&=: D_1^{(n)} + \frac{2(n-2i)}{n-i} D_2^{(n)} + \frac{(n-i)^2-(n-i)-2(n-2i)}{n-i}  D_3^{(n)},\quad\text{say.}
\end{align*}
Noting that 
 $D_2^{(n)} \leq ({\rm E} \Vert \al_{2, 1}^{(n)} \Vert^2)^{1/2} ({\rm E} \Vert \al_{3, 2}^{(n)} \Vert^2)^{1/2} = {\rm E} \Vert \al_{2, 1}^{(n)} \Vert^2 = D_1^{(n)},$ 
 it suffices to prove that~$D_1^{(n)} = o(1)$ and $D_3^{(n)} = o(n^{-1})$.

We first deal with $D_1^{(n)}$. It follows from \eqref{ass.a} that ${\rm E} \left\Vert \bxi_{2, 1}^{(n)} -  \bxi_{2, 1} \right\Vert^2 = o(1)$. Also, by Jensen's inequality,
$\Vert {\mbf \gamma}^{(n)} - {\mbf \gamma} \Vert^2   \leq {\rm E} \left\Vert \bxi_{2, 1}^{(n)} - \bxi_{2, 1} \right\Vert^2$. 
Therefore,%\vspace{-2mm}
\begin{align*}
{\rm E}\left\Vert \al_{2, 1}^{(n)} \right\Vert^2 = {\rm E} \left\Vert \bxi_{2, 1}^{(n)} -  {\mbf \gamma}^{(n)} - \bxi_{2, 1} + {\mbf \gamma} \right\Vert^2 
 \leq  2 {\rm E} \left\Vert \bxi_{2, 1}^{(n)} - \bxi_{2, 1} \right\Vert^2 + 2 \Vert  {\mbf \gamma}^{(n)} -  {\mbf \gamma} \Vert^2 
  = o(1).
\end{align*}

 In order to show that $D_3^{(n)}\!\! = o(n^{-1})$, consider an order statistic
$\ZZ_{(.)}^{(n)}\!\! :=\! (\ZZ_{(1)}^{(n)}, \ldots, \ZZ_{(n)}^{(n)})$, where $\ZZ_{(i)}^{(n)}$ is such that its first component is the $i$th order statistic of the $n$-tuple of first components. According to \cite{Hallinetal2020}, $\ZZ_{(.)}^{(n)}$ and $(\F^{(n)}_{\pm, 1}, \ldots, \F^{(n)}_{\pm, n})$ are mutually independent. 
 Conditional on~$\ZZ_{(.)}^{(n)}= {\bf z}\n\in\mathbb{R}^{nd}$, the distribution of~$\al_{2, 1}^{(n)}$ thus is that of   a score function (the form of which depends on ${\bf z}\n$)  computed at the empirical center-outward ranks and signs. 
 Hence, in view of Lemma~\ref{a1234n},
\begin{align*}
\left\vert D_3^{(n)} \right\vert &= \left\vert {\rm E}\left[ {\rm E}(\al_{2, 1}^{(n)\prime}  \al_{4, 3}^{(n)} | \ZZ\n_{(.)}) \right] \right\vert \\
& \leq \frac{n(n-1)}{(n-2)(n-3)} {\rm E} \left[ {\rm E} (\al_{2, 1}^{(n)\prime} | \ZZ\n_{(.)})  {\rm E} (\al_{2, 1}^{(n)} | \ZZ\n_{(.)}) \right] + 
 \frac{4(n-2) + 2}{(n-2)(n-3)} {\rm E}\left[ {\rm E} (\Vert \al_{2, 1}^{(n)} \Vert^2 | \ZZ\n_{(.)}) \right] \\
 &= \frac{n(n-1)}{(n-2)(n-3)}  {\rm E} \left[ {\rm E} (\al_{2, 1}^{(n)\prime} | \ZZ\n_{(.)})  {\rm E} (\al_{2, 1}^{(n)} | \ZZ\n_{(.)}) \right] + \frac{4(n-2) + 2}{(n-2)(n-3)} {\rm E} \left\Vert \al_{2, 1}^{(n)} \right\Vert^2 \\
 &=: \frac{n(n-1)}{(n-2)(n-3)} D_{31}^{(n)} + \frac{4(n-2) + 2}{(n-2)(n-3)} D_{32}^{(n)}
\end{align*}
where the second term above is $o(n^{-1})$ since  $D_{32}^{(n)}={\rm E} \left\Vert \al_{2, 1}^{(n)} \right\Vert^2$ is $o(1)$. 

It remains to show that $D_{31}^{(n)} = o(n^{-1})$. Since ${\mbf \gamma}^{(n)} = {\rm E}\, \bxi_{2, 1}^{(n)}$ and ${\mbf \gamma}  = {\rm E}\, \bxi_{2, 1}$,
\begin{align*}
 D_{31}^{(n)} 
&= {\rm E} \left[ ({\rm E} \bxi_{2, 1}^{(n)\prime} -  {\mbf \gamma}^{(n)\prime} - {\rm E}(\bxi_{2, 1}^\prime | \ZZ\n_{(.)}) + {\mbf \gamma}^{\prime})  ({\rm E}\bxi_{2, 1}^{(n)} - {\mbf \gamma}^{(n)} - {\rm E}(\bxi_{2, 1} | \ZZ\n_{(.)})  + {\mbf \gamma}) \right]  \\
&= {\rm E} \left[  ( {\rm E}(\bxi_{2, 1}^\prime | \ZZ\n_{(.)}) - {\mbf \gamma}^{\prime})  ( {\rm E}(\bxi_{2, 1} | \ZZ\n_{(.)})  - {\mbf \gamma}) \right]   \\
&= {\rm E} \left[  {\rm E}(\bxi_{2, 1}^\prime | \ZZ\n_{(.)})  {\rm E}(\bxi_{2, 1} | \ZZ\n_{(.)}) \right] - {\rm E}(\bxi_{2, 1}^\prime)  {\rm E}(\bxi_{2, 1}).
\end{align*}
Noting that 
$${\rm E}(\bxi_{2, 1} | \ZZ\n_{(.)}) = {\rm E}(\bxi_{R\n_2\!, R\n_1})
= \frac{1}{n(n-1)} \underset{1 \leq j_1 \neq j_2 \leq n}{\sum \sum} \bxi_{j_1, j_2}$$
where $R\n_2$ and $R\n_1$ denote the ranks of the first components of $\ZZ\n_2$ and $\ZZ\n_1$, respectively, 
we have%\vspace{-2mm}  
\begin{align*}
D_{31}^{(n)} 
&= \frac{1}{n^2 (n-1)^2} \underset{1 \leq j_1 \neq j_2 \leq n, 1 \leq j_3 \neq j_4 \leq n}{\sum \sum \sum \sum} \left[ {\rm E}(\bxi_{j_1, j_2}^\prime \bxi_{j_3, j_4}) -  {\rm E}(\bxi_{2, 1}^\prime)  {\rm E}(\bxi_{2, 1}) \right] \\
&=  \frac{1}{n^2 (n-1)^2} \underset{1 \leq j_1 \neq j_2 \leq n, 1 \leq j_3 \neq j_4 \leq n}{\sum \sum \sum \sum} \big( {\rm E} \left[(\J_1\pr(\F_{\pm, j_1}) \J_1(\F_{\pm, j_3})) (\J_2\pr(\F_{\pm, j_2}) \J_2(\F_{\pm, j_4}))\right] 
  \\
&\qquad \qquad \qquad \qquad \qquad \qquad \qquad \quad \quad
- \left\Vert {\rm E}   (\J_1(\F_{\pm, 2})) \right\Vert^2 \left\Vert {\rm E}   (\J_2(\F_{\pm, 1})) \right\Vert^2 \big).
\end{align*}
%\vspace{-8mm}
%
%\noindent 
Since $\ZZ_1^{(n)}, \ldots, \ZZ_n^{(n)}$ are i.i.d., only $2n(n-1)$ terms  with $j_1 \neq j_2$ and $j_3 \neq j_4$ 
%$\#\left[ \{j_1 \neq j_2\} \cap \{j_3 \neq j_4\}\right] = 2$
 contribute to this latter summation. Therefore, in view of the square-integrability of $\J_1$ and $\J_2$, we have~$D_{31}^{(n)}  = O(n^{-2})$.  The result follows.  \qed \vspace{3mm}

\textbf{Proof of Lemma~\ref{lem.xi}}.  
In view of  LAN, it is sufficient to prove the result  under ${\rm P}^{(n)}_{\bth_0 ;f}$  (with $f\in{\cal F}_d$).
It follows from the continuity of $\J_1$ and $\J_2$ and the Glivenko-Cantelli theorem in \cite{Hallinetal2020} that $\bxi_{2, 1; {\rm a}}^{(n)}$ converges to $\bxi_{2, 1}$ a.s.. Moreover,  square-integrability  of $\J_1$ and $\J_2$ and  independence of $\ZZ_1$ and $\ZZ_2$ entail  
$${\rm E} \left\Vert \bxi_{2, 1; {\rm a}} \right\Vert^2 = {\rm E} \left(\left\Vert  \J_1(\F_{\pm, 2})  \right\Vert^2\right)  {\rm E} \left(\left\Vert  \J_2(\F_{\pm, 1})   \right\Vert^2\right) < \infty.$$ 
Hence, ${\rm E} \left\Vert \bxi_{2, 1; {\rm a}}^{(n)} \right\Vert^2 - {\rm E} \left\Vert \bxi_{2, 1; {\rm a}} \right\Vert^2=o(1)$. It follows (see, e.g., part {\it (iv)} of Theorem 5.7 in Chapter~3 of~\citet{Shorack20}) that 
${\rm E} \left\Vert \bxi_{2, 1; {\rm a}}^{(n)} - \bxi_{2, 1}\right\Vert^2 =o(1).$  

Turning to the second part of \eqref{lem21}, put %\vspace{-2mm} 
$$\bxi_{2, 1}(\omega) = \text{vec} (\J_1(\F_{\pm, 2}(\omega)) \J_2\pr(\F_{\pm, 1}(\omega))), \ \bxi_{2, 1; {\rm a}}^{(n)}(\omega) = \text{vec} (\a_{\rm a} (\F^{(n)}_{\pm, 2}(\omega), \F^{(n)}_{\pm, 1}(\omega)))$$
for  any   point $\omega$ in the sample space $\Omega$. 
Since $\Vert \bxi_{2, 1; {\rm a}}^{(n)} - \bxi_{2, 1} \Vert$ converges to zero a.s., it follows from the Egorov theorem (see, e.g., part {\it (ii)} of Exercise 5.8 in \citet[Chapter~3]{Shorack20}) that, for any $\varepsilon > 0$, there exists a subset $\Acal \subset \Omega$ such that 
\begin{equation}\label{Egorov}
{\rm P}(\Acal) > 1- \varepsilon \quad \text{and}  \quad \underset{\omega  \in  \mathcal{A}}{\sup} \left\Vert \bxi_{2, 1; {\rm a}}^{(n)}(\omega) - \bxi_{2, 1}(\omega) \right\Vert =o(1).
\end{equation}
Denoting by $\Acal^c$ the complement of $\Acal$ in $\Omega$, we have
\begin{align*}
{\rm E} \left\Vert \bxi_{2, 1} -\bxi_{2, 1; {\rm e}}^{(n)} \right\Vert^2 
&= {\rm E} \left\Vert \bxi_{2, 1}(\1_{\Acal} + \1_{\Acal^c}) - {\rm E}[\bxi_{2, 1}(\1_{\Acal} + \1_{\Acal^c}) |  \F^{(n)}_{\pm, 2}, \F^{(n)}_{\pm, 1}]\right\Vert^2 \\
&\leq 3 {\rm E} \left\Vert \bxi_{2, 1}\1_{\Acal}  - {\rm E}[\bxi_{2, 1} \1_{\Acal}  |  \F^{(n)}_{\pm, 2}, \F^{(n)}_{\pm, 1}]\right\Vert^2    + 3 {\rm E} \left\Vert  \bxi_{2, 1} \1_{\Acal^c} \right\Vert^2 \\
&\qquad \qquad  + 3 {\rm E} \Vert  {\rm E}[\bxi_{2, 1} \1_{\Acal^c} |  \F^{(n)}_{\pm, 2}, \F^{(n)}_{\pm, 1}]\Vert^2 \\
&=: 3 (I_1^{(n)} + I_2^{(n)} + I_3^{(n)}).
\end{align*}
By Jensen's inequality,
$$I_3^{(n)} \leq {\rm E} \left( {\rm E}[\Vert  \bxi_{2, 1} \1_{\Acal^c} \Vert^2 |  \F^{(n)}_{\pm, 2}, \F^{(n)}_{\pm, 1}] \right) = {\rm E}\Vert  \bxi_{2, 1} \1_{\Acal^c} \Vert^2 = I_2^{(n)}$$
where, due to the square-integrability of $\J_1$ and $\J_2$ and the independence between $\ZZ_1$ and~$\ZZ_2$,~$I_2^{(n)}$ is arbitrarily small as $\varepsilon \rightarrow 0$. 

It remains to prove that $I_1^{(n)} = o(1)$. Denoting by $\mathfrak{G}\n$ the regular grid   the empirical center-outward distribution function is mapping to,   define, for any  $\mathfrak{g}_1, \mathfrak{g}_2 \in \mathfrak{G}\n$, 
$$\Bcal_{\mathfrak{g}_1 \mathfrak{g}_2}^{(n)} := \{\omega : \F^{(n)}_{\pm, 1}(\omega) = \mathfrak{g}_1, \F^{(n)}_{\pm, 2}(\omega)= \mathfrak{g}_2\}.$$
We have
\begin{align*}
I_1^{(n)} &= {\rm E}\, \Big\Vert \sum_{\mathfrak{g}_1 \in \mathfrak{G}\n} \sum_{\mathfrak{g}_2 \in \mathfrak{G}\n} \1_ {\Bcal_{\mathfrak{g}_1 \mathfrak{g}_2}^{(n)}} \left(  \bxi_{2, 1}\1_{\Acal}  - {\rm E}[\bxi_{2, 1} \1_{\Acal}  |  \F^{(n)}_{\pm, 2}, \F^{(n)}_{\pm, 1}] \right) \Big\Vert^2 \\
&= {\rm E}  \sum_{\mathfrak{g}_1 \in \mathfrak{G}\n} \sum_{\mathfrak{g}_2 \in \mathfrak{G}\n} \left\Vert \1_ {\Bcal_{\mathfrak{g}_1 \mathfrak{g}_2}^{(n)}} \left(  \bxi_{2, 1}\1_{\Acal}  - {\rm E}[\bxi_{2, 1} \1_{\Acal}  |  \F^{(n)}_{\pm, 2}, \F^{(n)}_{\pm, 1}] \right) \right\Vert^2,
\end{align*}
where the second equality follows from the fact that  $\1_ {\Bcal_{\mathfrak{g}_1 \mathfrak{g}_2}^{(n)}} \1_ {\Bcal_{\mathfrak{h}_1 \mathfrak{h}_2}^{(n)}} = \1_ {(\mathfrak{g}_1 = \mathfrak{h}_1, \mathfrak{g}_2 = \mathfrak{h}_2)}.$ Note that
$$\1_ {\Bcal_{\mathfrak{g}_1 \mathfrak{g}_2}^{(n)}} {\rm E}\left[\bxi_{2, 1} \1_{\Acal}  |  \F^{(n)}_{\pm, 2}, \F^{(n)}_{\pm, 1}\right]  = \frac{\1_ {\Bcal_{\mathfrak{g}_1 \mathfrak{g}_2}^{(n)}}}{{\rm P}(\Bcal_{\mathfrak{g}_1 \mathfrak{g}_2}^{(n)})} \int_{\eta  \in \Bcal_{\mathfrak{g}_1 \mathfrak{g}_2}^{(n)}} \bxi_{2, 1}(\eta) \1_{\Acal}(\eta ) {\rm d} {\rm P} (\eta).$$
Therefore, 
\begin{align*}
I_1^{(n)} &=  \sum_{\mathfrak{g}_1 \in \mathfrak{G}\n} \sum_{\mathfrak{g}_2 \in \mathfrak{G}\n} \int_{\Omega} \Big\Vert \1_ {\Bcal_{\mathfrak{g}_1 \mathfrak{g}_2}^{(n)}}(\omega) 
\int_{\eta \in \Bcal_{\mathfrak{g}_1 \mathfrak{g}_2}^{(n)}} \left[ \bxi_{2, 1}(\omega) \1_{\Acal}(\omega) \right.  \\
&\qquad \qquad \qquad \qquad \qquad \qquad \qquad \qquad  \left.   - \bxi_{2, 1}(\eta) \1_{\Acal}(\eta) \right] \frac{{\rm d} {\rm P} (\eta)}{{\rm P}(\Bcal_{\mathfrak{g}_1 \mathfrak{g}_2}^{(n)})} \Big\Vert^2 {\rm d} {\rm P} (\omega)  \\
&= \sum_{\mathfrak{g}_1 \in \mathfrak{G}\n} \sum_{\mathfrak{g}_2 \in \mathfrak{G}\n} \int_{\Omega} \Big\Vert \1_ {\Bcal_{\mathfrak{g}_1 \mathfrak{g}_2}^{(n)}}(\omega) 
\int_{\eta \in \Bcal_{\mathfrak{g}_1 \mathfrak{g}_2}^{(n)}} \left[ \left( \bxi_{2, 1}(\omega) - \bxi_{2, 1; {\rm a}}^{(n)} (\omega) \right) \1_{\Acal}(\omega) \right.  \\
&\qquad \qquad \qquad \qquad    \qquad  \left.  +  \left( \bxi_{2, 1; {\rm a}}^{(n)} (\eta) - \bxi_{2, 1}(\eta) \right) \1_{\Acal}(\eta) \right] \frac{{\rm d} {\rm P} (\eta)}{{\rm P}(\Bcal_{\mathfrak{g}_1 \mathfrak{g}_2}^{(n)})} \Big\Vert^2 {\rm d} {\rm P} (\omega) 
\end{align*}
since,  on $\Acal \cap \Bcal_{\mathfrak{g}_1 \mathfrak{g}_2}^{(n)}$, 
\begin{align*}
\bxi_{2, 1; {\rm a}}^{(n)} (\omega)   \1_{\Acal}(\omega) \1_ {\Bcal_{\mathfrak{g}_1 \mathfrak{g}_2}^{(n)}}(\omega)
%&
= \text{vec}(\J_1(\mathfrak{g}_2) \J_2^\prime(\mathfrak{g}_1))% \\ 
%&
= \bxi_{2, 1; {\rm a}}^{(n)} (\eta)   \1_{\Acal}(\eta) \1_ {\Bcal_{\mathfrak{g}_1 \mathfrak{g}_2}^{(n)}}(\eta).
\end{align*}
  In view of \eqref{Egorov}, for any $\tilde{\varepsilon} > 0$, there exists $n({\tilde{\varepsilon}})$ such that $\Vert \bxi_{2, 1; {\rm a}}(\omega) - \bxi_{2, 1}^{(n)} (\omega) \Vert^2 < \tilde{\varepsilon}$ for all~$n > n({\tilde{\varepsilon}})$ and all $\omega \in  \mathcal{A}$. Then, using Jensen's inequality again, 
\begin{align*}
I_1^{(n)} &\leq \sum_{\mathfrak{g}_1 \in \mathfrak{G}\n} \sum_{\mathfrak{g}_2 \in \mathfrak{G}\n} \int_{\Omega}  \1_ {\Bcal_{\mathfrak{g}_1 \mathfrak{g}_2}^{(n)}}(\omega) 
\int_{\eta \in \Bcal_{\mathfrak{g}_1 \mathfrak{g}_2}^{(n)}} \left[ 2 \Vert \bxi_{2, 1}(\omega) - \bxi_{2, 1; {\rm a}}^{(n)} (\omega) \Vert^2 \1_{\Acal}(\omega) \right.  \\
&\qquad \qquad \qquad \qquad \qquad  \qquad  \left.  +  2 \Vert \bxi_{2, 1; {\rm a}}^{(n)} (\eta) - \bxi_{2, 1}(\eta)  \Vert^2 \1_{\Acal}(\eta) \right] \frac{{\rm d} {\rm P} (\eta)}{{\rm P}(\Bcal_{\mathfrak{g}_1 \mathfrak{g}_2}^{(n)})} {\rm d} {\rm P} (\omega) \\
&\leq   \sum_{\mathfrak{g}_1 \in \mathfrak{G}\n} \sum_{\mathfrak{g}_2 \in \mathfrak{G}\n} \int_{\Omega} \1_ {\Bcal_{\mathfrak{g}_1 \mathfrak{g}_2}^{(n)}}(\omega) 4  \tilde{\varepsilon}^2 \frac{{\rm P}(\Bcal_{\mathfrak{g}_1 \mathfrak{g}_2}^{(n)})}{{\rm P}(\Bcal_{\mathfrak{g}_1 \mathfrak{g}_2}^{(n)})} {\rm d} {\rm P} (\omega) = 4  \tilde{\varepsilon}^2 {\rm E} \sum_{\mathfrak{g}_1 \in \mathfrak{G}\n} \sum_{\mathfrak{g}_2 \in \mathfrak{G}\n} \1_ {\Bcal_{\mathfrak{g}_1 \mathfrak{g}_2}^{(n)}} = 4  \tilde{\varepsilon}^2.
\end{align*}
The result follows. \qed

\appendix\section*{Appendix C:  Proofs of Lemma~\ref{asy.Gami.bar}, Proposition~\ref{Prop.Stilde}, and Proposition~\ref{Prop.Wtilde}}\label{sec.C}

\renewcommand\theequation{C\arabic{equation}}
\setcounter{equation}{0}

\hspace{7mm}\textbf{Proof of Lemma~\ref{asy.Gami.bar}}.  %\begin{proof}
The result follows from deriving the asymptotic joint distribution, under ${\rm P}^{(n)}_{\bth_0 ;f}$,  of
$$(n-i)^{1/2}\text{\rm vec} (\underline{\bGamma}_{i, {\J}_1, {\J}_2}^{(n)}(\bth_0) - \m),\quad (n-j)^{1/2}\text{\rm vec} (\underline{\bGamma}_{j, {\J}_1, {\J}_2}^{(n)}(\bth_0) - \m),\quad\text{and}\quad \bDelta^{(n)}_{f} (\bth_0)$$ 
along the same lines as in the proof of Lemma~B.1 in \cite{HLL2019}. An application of Le Cam's third Lemma concludes. 
 Details are left to the reader.\qed\bigskip
 
\textbf{Proof of Proposition~\ref{Prop.Stilde}}
Parts  {\it (i)} and  {\it (ii)} directly follow from Lemma~\ref{asy.Gami.bar} and Proposition~\ref{Prop.HajekRep}. Turning to part~{\it (iii)}, notice that $\underline{\bGamma}_{i, \J_1, \J_2}^{(n)}(\bth_0)$ coincides with $\bGamma_{i, f}^{(n)}(\bth_0)$  for score functions $\J_{\ell} =  \J_{\ell; f}$, $\ell = 1, 2$. Moreover, we then have $\m = {\rm E}_{{\bth_0};f}(\bGamma_{i, f}^{(n)}(\bth_0)) = \0$ since~${\rm E}_{{\bth_0};f}(\ZZ_t(\bth_0)) = \0$. The result then follows from using Proposition~\ref{Prop.HajekRep} and moving along the same lines as the proof of Proposition 4 {\it (v)} in \cite{HP04}.\qed\bigskip

\textbf{Proof of Proposition~\ref{Prop.Wtilde}}.  
Proposition~\ref{asy0} implies that, for $\btau =  (\btau_I^{\prime} \ \0^{\prime})^\prime$, under any~${\rm P}^{(n)}_{\bth_0}$ in  $\mathcal{H}^{(n)}_{0}$ and contiguous alternatives,%\vspace{-3mm} 
\begin{equation}\label{asy.linear3}
\utilde{\bDelta}^{(n)}_{I; \J_1, \J_2}(\bth_0 + n^{-1/2} \btau) - \utilde{\bDelta}^{(n)}_{I; \J_1, \J_2}(\bth_0) =  -  \bUpsilon^{(n)}_{11; \J_1, \J_2, f}(\bth_0)  \btau_I + o_{\rm P}(1)
\end{equation}
%\vspace{-7mm}
%
%\noindent
 and
\begin{equation}\label{asy.linear4}
\utilde{\bDelta}^{(n)}_{II; \J_1, \J_2}(\bth_0 + n^{-1/2} \btau) - \utilde{\bDelta}^{(n)}_{II; \J_1, \J_2}(\bth_0) =  -  \bUpsilon^{(n)}_{21; \J_1, \J_2, f}(\bth_0)  \btau_I + o_{\rm P}(1).
\end{equation}
Therefore, by the definition \eqref{defDeltastar} of $\utilde{\bDelta}^{(n)*}_{II; \J_1, \J_2}(\bth_0)$, we  have%\vspace{-2mm}
\begin{equation}\label{DeltaItau.diff}
\utilde{\bDelta}^{(n)*}_{II; \J_1, \J_2}(\bth_0 + n^{-1/2} \btau) - \utilde{\bDelta}^{(n)*}_{II; \J_1, \J_2}(\bth_0) = o_{\rm P}(1).
\end{equation}
%\vspace{-8mm}
%
%\noindent 
Part {\it (i)} then follows from consistency of $\hat{\bUpsilon}^{(n)}$ and root-$n$ consistency and asymptotically discreteness of $\hat{\bth}^{(n)}$ (which allows us to replace $\bth_0 + n^{-1/2} \btau$ with  $\hat{\bth}^{(n)}$ in \eqref{DeltaItau.diff}; see Lemma~4.4 in \cite{Kreiss87}).

Now, part {\it (i)} of the proposition  implies that $\widehat{\utilde{W}}^{(n)}_{\J_1 \J_2}(\hat{\bth}^{(n)})$ has the same limit distribution as $\utilde{W}^{(n)}_{\J_1 \J_2}(\bth_0)$  under $\mathcal{H}^{(n)}_{\bth_0}$
 and under contiguous alternatives. It follows from   Proposition~\ref{til.Delta.n} that~$\utilde{\bDelta}^{(n)*}_{II; \J_1, \J_2}(\bth_0)$  is asymptotically $d^2(p_1 - p_0)$-variate normal,   with mean~$\0$ under  $\mathcal{H}^{(n)}_{\bth_0}$,  mean~$\left( \bUpsilon_{22} - \bUpsilon_{21} (\bUpsilon_{11})^{-1} \bUpsilon_{12} \right) \btau_{II}$ 
under ${\rm P}^{(n)}_{\bth_0+ n^{-1/2} \btau}$ with~$\btau\pr=(\btau_I\pr, \btau_{II}\pr)$,  and   co\-variance~$\bLam^{*}_{II; \bth_0}$ 
 under both. Parts  {\it (ii)} and {\it (iii)} follow.
 
  %%%%%%%%%%%%%%%%%%%%%%%%%%%%

 Finally, part {\it (iv)} follows from Lemma~\ref{lem.bar.til.Delta.n}  by  
 noticing that, for the scores $\J_{1;f}  := \bvp_{f} \circ~\!\Q_{\pm}$ and $\J_{2;f}  :=  \Q_{\pm}$, $\m = \0$ and  $\underline{\bGamma}_{i, \J_1, \J_2}^{(n)}(\bth_0)=\bGamma_{i, f}^{(n)}(\bth_0)$.  Local asymptotic stringency then follows from general results on asymptotically optimal tests in LAN families: see Section~11.9 of  Le Cam~(1986).
\qed

\appendix\section*{Appendix D:  Further numerical results}\label{appBias}

The same Monte Carlo experiment as in Section~\ref{sec.d2} is  conducted here in dimension $d=3$. We  generated $N = 1000$ replications of size   $n = 1400$  ($n_R = 21$, $n_S = 66$, and $n_0 = 14$) from 
the trivariate VAR model
\begin{equation} \X_t - \ell \A \X_{t-1}= \bepsilon_t, \quad \ell = 0, 1, 1.5,\label{trimod}\end{equation}
with $\text{vec}(\A) = (0.05, 0.01, 0.011, 0.01, 0.02, 0.01, -0.01, 0.013, 0.033)^\prime$ and     spherical normal and $t_3$,  mixtures of  normal,   and skew-$t_3$ innovation densities. The mixtures are of the form\vspace{-2mm}
\begin{equation}
\frac{3}{8}   {\cal N}(\bmu_1, \bSigma_1) + \frac{3}{8}   {\cal N}(\bmu_2, \bSigma_2) + \frac{1}{4} {\cal N}(\bmu_3, \bSigma_3),
\label{Eq. Mixture}
\end{equation}
with
 $\bmu_1 = (-5, -5, 0)^\prime,\ \bmu_2 = (5, 5, 2)^\prime, \ \bmu_3 = (0, 0, -3)^\prime$ 
and\vspace{-2mm} 
 $$\bSigma_1 = 
\begin{bmatrix}
7 & 3 & 5\vspace{-1mm} \\
3 & 6 & 1\vspace{-1mm} \\
5 & 1 & 7\vspace{-1mm}
\end{bmatrix}, \  
\bSigma_2 = 
\begin{bmatrix}
7 & -5 & -3\vspace{-1mm} \\
-5 & 7 & 4\vspace{-1mm} \\
-3 & 4 & 5\vspace{-1mm} 
\end{bmatrix}, \ 
\bSigma_3 = 
\begin{bmatrix}
4 & 0 & 0\vspace{-1mm}\\
0 & 3  & 0\vspace{-1mm} \\
0 & 0 & 1\vspace{-1mm}
\end{bmatrix}
.$$

The  skew-$t_3$  distribution  ($d=\nu =3$)  has  density \eqref{St.density}  with    
 $\bxi = \0$,    $\al = (7,5,3)^\prime$,  and~vec$\left(\bSigma\right) = \left(7, -5,-3,-5,7,4,-3,4,5\right)\pr
$. 

\begin{table}[!b]
\caption{Rejection frequencies (out of $N=1000$ replications),  under  values~$\ell \A$, $\ell=0, 1, 1.5$  of the VAR(1) autoregression matrix, for the trivariate model~\eqref{trimod} and various innovation densities (Gaussian, mixture~\ref{Eq. Mixture}, $t_3$, skew-$t_3$, and contaminated Gaussian), of the Gaussian, vdW, bias-corrected vdW, Spearman, bias-corrected Spearman,  and  sign  tests of white noise against VAR($1$); the sample size is $n = 1400$   ($n_R = 21$, $n_S = 66$, and~$n_0 = 14$); the nominal level is $\alpha=5\%$;  rejection frequencies of the sign tests with $n_S=n$ are also included. Permutational critical values are based on $M=5000$ random permutations. }\label{Tab.d3}%\vspace{2mm}

\centering

\footnotesize

\begin{tabular}{lcccc | lcccc}

\hline\hline

$\quad\ f$ &\hspace{-31mm} Test     &\hspace{-2mm}  $\0$  & $\A$  & $1.5\A$ &   $\quad\ f$      &\hspace{-31mm} Test     &\hspace{-2mm} $\0$  & $\A$  & $1.5\A$ \\ \hline\hline

\underline{Normal}\vspace{-1mm}&  &  &  & & \underline{Mixture}    & &  &   &   \\

 &\hspace{-31mm}  Gaussian &\hspace{-2mm}  0.044 & 0.354 & 0.751 &     &\hspace{-31mm} Gaussian &\hspace{-2mm}  0.053 & 0.490 & 0.885 \\

       &\hspace{-31mm} vdW      &\hspace{-2mm}  0.021 & 0.213 & 0.600 &            &\hspace{-31mm} vdW      &\hspace{-2mm}  0.029 & 0.628 & 0.981 \\

     {} &{\hspace{-31mm} bias-corrected vdW} &\hspace{-2mm} {0.047} &{0.325} &{0.729} & {} &{\hspace{-31mm} bias-corrected vdW} &\hspace{-2mm} {0.056} &{0.736} &{0.987} \\

% &              (bias-corrected)                        &                        &                        &                        &  &                          (bias-corrected)             &                        &                        &                      \\

       &\hspace{-31mm} Spearman &\hspace{-2mm}  0.028 & 0.242 & 0.629 &            &\hspace{-31mm} Spearman &\hspace{-2mm}  0.041 & 0.686 & 0.985  \\

      {} &\hspace{-31mm} bias-corrected Spearman &\hspace{-2mm} {0.048} &{0.285} &{0.696} & {} &\hspace{-31mm} bias-corrected Spearman &\hspace{-2mm} {0.061} &{0.735} &{0.987} \\

%                  &         (bias-corrected)                               &                        &                        &                        &                   &                                       (bias-corrected) &                        &                        &               \\

       &\hspace{-31mm}  Sign   ($n_S =66$)   &\hspace{-2mm}  0.039 & 0.229 & 0.534 &            &\hspace{-31mm}  Sign  ($n_S =66$)    &\hspace{-2mm}  0.036 & 0.457 & 0.864 \\

%   &                                        &                        &                        &                        &                   &                                         &                        &                        &               \\ 

%   &\hspace{-31mm} {bias-corrected Sign  ($n_R =20$)} &\hspace{-2mm} {0.050} &{0.265} &{0.575} &  &\hspace{-31mm} {bias-corrected Sign  ($n_R =20$)} &\hspace{-2mm} {0.050} &{0.495} &{0.882} \\

 & \hspace{-31mm}{Sign ($n_S =1400$)} & \hspace{-2mm} 0.043 & 0.238 & 0.568 &  & \hspace{-31mm}{Sign ($n_S =1400$)} & 0.043 & 0.492 & 0.898 \\

%             &         (bias-corrected I)                               &                        &                        &                        &                   &                                       (bias-corrected I) &                        &                        &               \\

%             & {\hspace{-31mm} bias-corrected Sign ($n_R =1$)} &\hspace{-2mm} {0.051} &{0.246} &{0.605} &  &{\hspace{-31mm} bias-corrected Sign ($n_R =1$)} &\hspace{-2mm} {0.045} &{0.556} &{0.926} \\

%             &         (bias-corrected II)                               &                        &                        &                        &                   &                                       (bias-corrected II) &                        &                        &               \\

   \hline

$ \underline{\, t_3\,}$  &  &  &  &  & \underline{Skew-$t_3\,$} & &  &   &   \\

 &\hspace{-31mm} Gaussian &\hspace{-2mm} 0.035 & 0.376 & 0.763 &  &\hspace{-31mm} Gaussian &\hspace{-2mm} 0.038 & 0.412 & 0.819 \\

       &\hspace{-31mm} vdW      &\hspace{-2mm} 0.024 & 0.329 & 0.804 &            &\hspace{-31mm} vdW      &\hspace{-2mm} 0.021 & 0.637 & 0.980 \\

      {} &{\hspace{-31mm}bias-corrected vdW} &\hspace{-2mm}{0.060} &{0.464} &{0.888} & {} &{\hspace{-31mm} bias-corrected vdW} &\hspace{-2mm}{0.056} &{0.759} &{0.992} \\

%                  &         (bias-corrected)                               &                        &                        &                        &                   &                                       (bias-corrected) &                        &                        &               \\  

       &\hspace{-31mm} Spearman &\hspace{-2mm} 0.031 & 0.346 & 0.787 &            &\hspace{-31mm} Spearman &\hspace{-2mm} 0.035 & 0.673 & 0.983 \\

{} & {\hspace{-31mm} bias-corrected Spearman} &\hspace{-2mm} {0.046} & {0.411} & {0.834} &  {} & {\hspace{-31mm} bias-corrected Spearman} &\hspace{-2mm} {0.058} & {0.728} & {0.990} \\

%                  &         (bias-corrected)                               &                        &                        &                        &                   &                                       (bias-corrected) &                        &                        &               \\

       &\hspace{-31mm}  Sign   ($n_S =66$)   &\hspace{-2mm} 0.049 & 0.349 & 0.789 &            &\hspace{-31mm}  Sign   ($n_S =66$)   &\hspace{-2mm} 0.035 & 0.534 & 0.941 \\

% &          ($n_R =1$)                              &                        &                        &                        &                   &                                       ($n_R =1$) &                        &                        &               \\ 

%   & \hspace{-31mm}  bias-corrected Sign ($n_R =20$)&\hspace{-2mm} {0.062} & {0.395} & {0.815} &  &  \hspace{-31mm}  bias-corrected Sign ($n_R =20$) &\hspace{-2mm} {0.051} & {0.565} & {0.951} \\

&  \hspace{-31mm}  {Sign ($n_S =1400$)}  & \hspace{-2mm} 0.044 & 0.422 & 0.832 &  & \hspace{-31mm}  {Sign ($n_S =1400$)} & \hspace{-2mm} 0.052 & 0.501 & 0.924 \\

%       &         (})                               &                        &                        &                        &                   &                                       (bias-corrected I) &                        &                        &               \\

%       & \hspace{-31mm}  bias-corrected Sign ($n_R =1$) &\hspace{-2mm} {0.058} & {0.398} & {0.845} &  &  \hspace{-31mm}  bias-corrected Sign ($n_R =1$)&\hspace{-2mm} {0.044} & {0.552} & {0.937} \\

%         &         (bias-corrected II)                               &                        &                        &                        &                   &                                       (bias-corrected II) &                        &                        &               \\

 \hline\hline

\underline{AOs}  (${\mbf s} = (6, 6, 6)^\prime$) &  &  &  &  &  \underline{AOs}  (${\mbf s} = (9, 9, 9)^\prime$) &  &  &     &   \\

 &\hspace{-31mm}  Gaussian &\hspace{-2mm} 0.337 & 0.378 & 0.525 &    &\hspace{-31mm}  Gaussian &\hspace{-2mm} 0.745 & 0.768   & 0.854 \\

                              &\hspace{-31mm}  vdW      &\hspace{-2mm} 0.027 & 0.135 & 0.360 &                            &\hspace{-31mm}  vdW      &\hspace{-2mm} 0.031 & 0.125 & 0.356 \\

   {} &{\hspace{-31mm}bias-corrected vdW}&\hspace{-2mm} {0.059} & {0.215} & {0.497} &  {} & {\hspace{-31mm}bias-corrected vdW} &\hspace{-2mm} {0.059} & {0.223} & {0.499} \\

%                     &         (bias-corrected)                               &                        &                        &                        &                   &                                       (bias-corrected) &                        &                        &               \\  

                            &\hspace{-31mm} Spearman &\hspace{-2mm} 0.036 & 0.140 & 0.372 &                             &\hspace{-31mm} Spearman &\hspace{-2mm} 0.040 & 0.153 & 0.364 \\

                            {} &\hspace{-31mm}{bias-corrected Spearman} &\hspace{-2mm} {0.053} & {0.186} & {0.448} &  {} &\hspace{-31mm}{bias-corrected Spearman} &\hspace{-2mm} {0.061} & {0.190} & {0.432} \\

%                     &         (bias-corrected)                               &                        &                        &                        &                   &                                       (bias-corrected) &                        &                        &               \\  

                            &\hspace{-31mm} Sign    ($n_S =66$)  &\hspace{-2mm} 0.045 & 0.171 & 0.389 &                             &\hspace{-31mm} Sign   ($n_S =66$)   &\hspace{-2mm} 0.056 & 0.159 & 0.396         \\

%&          ($n_R =1$)                              &                        &                        &                        &                   &                                       ($n_R =1$) &                        &                        &               \\ 

%&\hspace{-31mm}{bias-corrected Sign ($n_R =20$)} & \hspace{-2mm}{0.059} & {0.198} & {0.435} &  &\hspace{-31mm}{bias-corrected Sign ($n_R =20$)} &\hspace{-2mm}{0.075} & {0.188} & {0.432} \\

                            &\hspace{-31mm}{Sign ($n_S =1400$)} &\hspace{-2mm} 0.052 & 0.162 & 0.389 &  &\hspace{-31mm}{Sign ($n_S =1400$)} &\hspace{-2mm} 0.063 & 0.165 & 0.397 \\

%&         (bias-corrected I)                               &                        &                        &                        &                   &                                       (bias-corrected I) &                        &                        &               \\

%&\hspace{-31mm}{bias-corrected Sign ($n_R =1$)} &\hspace{-2mm}{0.066} & {0.180} & {0.441} &  &\hspace{-31mm}{bias-corrected Sign ($n_R =1$)}  &\hspace{-2mm}{0.059} & {0.165} & {0.423} \\

%&         (bias-corrected II)                               &                        &                        &                        &                   &                                       (bias-corrected II) &                        &                        &               \\

\hline\hline
\vspace{-8mm}
\end{tabular}

\end{table}%\vspace{-8mm}

To investigate  robustness, % of the center-outward rank-based test,  
we also considered   contaminated spherical Gaussian   $\X_t$'s 
 of the form 
 $\{\X^*_t=  \X_t + \1_{\{t = h\}} {\mbf s}\}$, where we set~$h$   in order to obtain  
$5\%$ of equally spaced additive outliers and put    ${\mbf s} = (6, 6,6)^\prime$ and $(9, 9,9)^\prime$, respectively. All contaminated 
observations were demeaned prior to the implementation of the testing procedures.

 %%%%%%%%%%%%%%%%%%%%%%%%%%%%

%%%%%%%%%%%%%%%%%%%%%%%%%%%%%%%%%%%%%

 \begin{figure}[!t]
     \centering
%     \begin{subfigure}[b]{0.4\textwidth}
%         \centering
%         \includegraphics[width=\textwidth]{vdWQQplot.eps}
%        % \caption{$y=x$}
%  %       \label{fig:y equals x}
%     \end{subfigure}
%%     \hfill
%     \begin{subfigure}[b]{0.4\textwidth}
%         \centering
%         \includegraphics[width=\textwidth]{GauSpSignQQplot.eps}
%       %  \caption{$y=3sinx$}
%      %   \label{fig:three sin x}
%     \end{subfigure}
%     \\
      \begin{subfigure}[b]{0.4\textwidth}
         \centering
         \includegraphics[width=\textwidth]{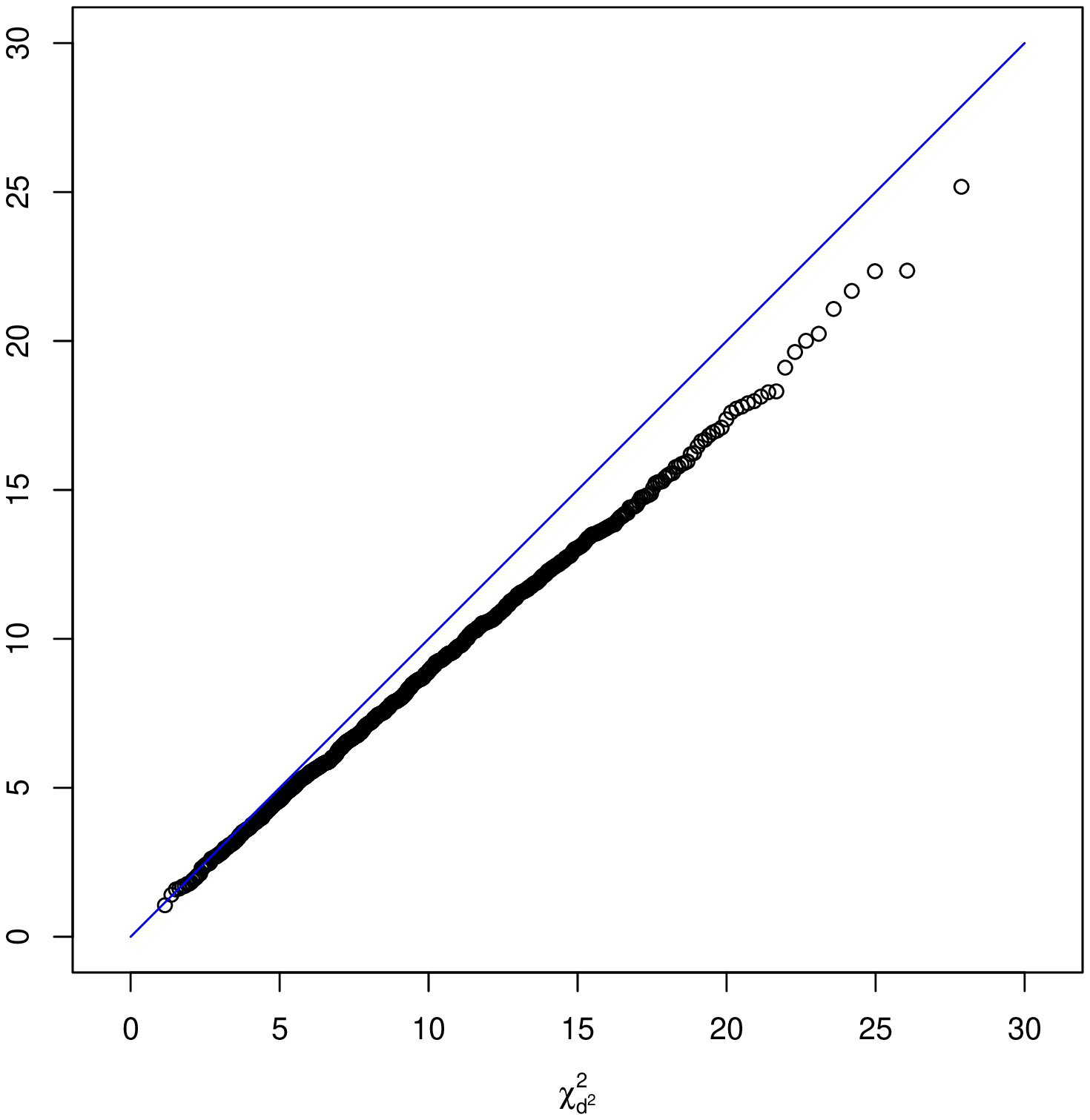}
        % \caption{$y=x$}
  %       \label{fig:y equals x}
    \end{subfigure}
%     \hfill
     \begin{subfigure}[b]{0.4\textwidth}
         \centering
         \includegraphics[width=\textwidth]{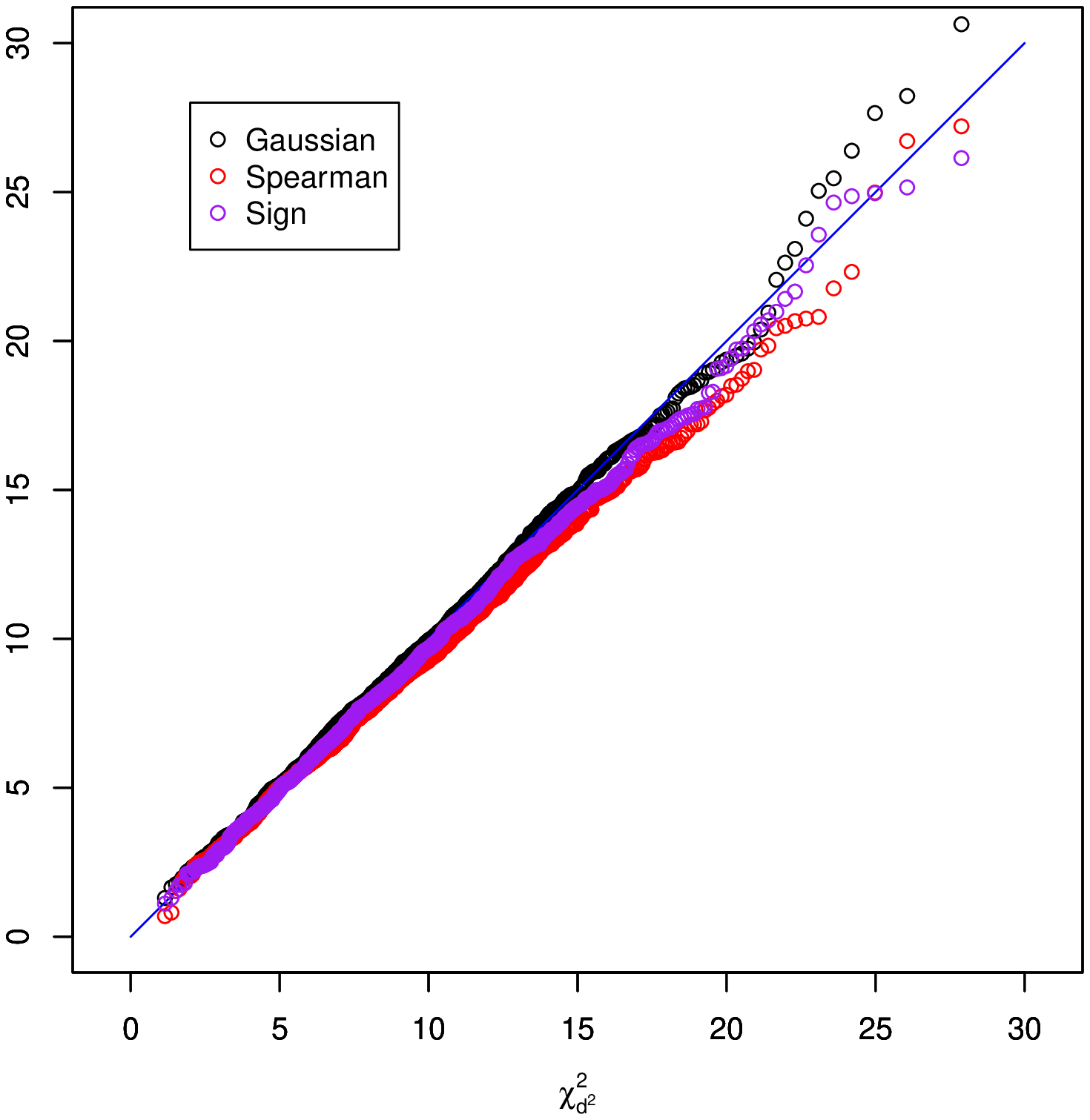}
       %  \caption{$y=3sinx$}
      %   \label{fig:three sin x}
  \end{subfigure}   %\vspace{-1mm}  
%     \hfill
%     \begin{subfigure}[b]{0.3\textwidth}
%         \centering
%         \includegraphics[width=\textwidth]{graph3}
%         \caption{$y=5/x$}
%         \label{fig:five over x}
%     \end{subfigure}
        \caption{QQ plots of the empirical distributions (across $N=1000$ replications) of the test statistics in  Tables~\ref{Tab.white} and~\ref{Tab.d3} against their asymptotically chi-square null distributions: vdW for~$d=2$ (top left)  and $d=3$ (bottom left); Gaussian, Spearman, and sign test scores  for~$d=2$ (top right)  and $d=3$ (bottom right).}
       \label{QQFig'}
\end{figure}

%In all cases, the asymptotic critical value based on the asymptotic chi-square distribution was used. 
Rejection frequencies at $5\%$ nominal level are reported  in Table~\ref{Tab.d3}  for  the Gaussian, vdW, Spearman and sign tests based on the asymptotic critical values  but also for the vdW and Spearman tests based on bias-corrected critical values as described in Section~\ref{sec.d2}. %Leaving aside the latter, %(postponed to Section~\ref{biasSec}),
  Inspection of  the   table confirm the findings of Section~\ref{sec.d2} and the same comments are in order here. %~\ref{Tab.white} and~\ref{Tab.d3}
   The vdW and Spearman   tests again suffer  a  bias  confirmed by   Figure~\ref{QQFig'}, which shows the QQ plots  of the   values of the vdW test statistic across the~$N=1000$ replications with Gaussian~$f$.  The same  relatively slow asymptotics   of center-outward ranks are to be blamed---the value of~$n_R=21$ for~$n=1400$ indeed  is quite small---and the problem  similarly is solved by resorting to permutational critical values.

\appendix\section*{Appendix E:  Computational aspects}\label{sec.comp}

In this section, we briefly discuss some computational aspects related to the implementation of our test statistics 
$\utilde{S}_{\J_1 \J_2}^{(n)}(\bth_0)$ and $\widehat{\utilde{W}}^{(n)}_{\J_1 \J_2}(\hat{\bth}^{(n)})$
 and we propose algorithms for ease of implementation.
%; see \cite{HLL2019} for similar comments in the context of center-outward rank-based estimation for VARMA models.

{\it (i)} The asymptotic distribution of our test statistic requires 
 that both $n_R $ and $n_S$  tend to infinity. In practice, we  factorize $n$  into~$n_Rn_S + n_0$ in such a way that both $n_R$ and $n_S$ are large. Typically, $n_R$  is of order~$n^{1/d}$ and~$n_S$ is of order~$n^{(d-1)/d}$, whilst $0 \leq n_0 < \min(n_S, n_R)$ has to be small as possible---its value, however, is entirely determined by  the values of $n_R$ and $n_S$.
%while $n_0$ is as small as possible. 
Generating  ``regular grids'' of $ n_S$ points over the unit sphere $\mathcal{S}_{d-1}$ as described in Section \ref{secranks} is easy  for $d=2$, where perfect regularity can be achieved by dividing the unit circle into $n_S$ arcs of equal length~$2\pi/n_S$. {For $d\geq 3$, 
``perfect regularity''  is no longer possible. 
A random array of $n_S$ independent and uniformly distributed unit vectors does satisfy (almost surely) the requirement for weak convergence to ${\mathrm U}_d$, representing the uniform distribution over the ball $\mathbb{S}_d$. More regular deterministic arrays (with faster convergence) can be constructed, though, such as the {\it low-discrepancy sequences}   (see, e.g., \cite{Niederreiter1992},  \cite{Judd1998}, or \cite{Dick2014}) considered  in numerical integration and the design of computer experiments; we suggest the use of the function {\tt UnitSphere} in R package {\tt mvmesh}. 
%See Hallin (2017, Section 4) for more details on how to generate the points on the grids.

{\it (ii)} The empirical center-outward distribution function $\F_{{\pms} }^{(n)} $ is obtained as the solution of  an optimal coupling problem. Many efficient algorithms have been proposed in the measure transportation literature (see, e.g., \cite{PC2019}). We followed  \cite{Hallinetal2020},  using a Hungarian algorithm (see the \texttt{clue} R package).  For a detailed account of the Hungarian algorithm and the complexity of different implementations, see, e.g., Chapter~4 in \cite{Betal09}). Faster algorithms are available, though, as Bertsekas’ auction algorithm or its variant, the forward/reverse auction algorithm, (Chapter 4 in \cite{Ber91}), implemented in the R package \texttt{transport}.

{\it (iii)} The computation of the test statistic $\widehat{\utilde{W}}^{(n)}_{\J_1 \J_2}(\hat{\bth}^{(n)})$  involves  two basic ingredients: a root-$n$ consistent estimator $\hat{\bth}^{(n)}$ of $\bth_0$ and a consistent estimator of the matrix~$\bUpsilon_{{\J}_1, {\J}_2, f}(\bth_0)$. For $\hat{\bth}^{(n)}$, a great number of candidates are available, e.g., the center-outward R-estimator of \cite{HLL2019}, the  reweighted multivariate least trimmed squares estimator of  \cite{Croux2008}, and  the QMLE (provided that  fourth-order moments  finite). Turning to the estimation of $\bUpsilon_{{\J}_1, {\J}_2, f}(\bth_0)$, the issue is that this matrix depends on the unknown actual density~$f$. 
A simple consistent estimator is obtained by letting $\btau= {\bf e}_i$, $i=1,\ldots, p_0 d^2$ in~\eqref{asy.linear2} where ${\bf e}_i$  denotes the $i$th vector of the canonical basis in the parameter space~${\mathbb{R}}^{p_1 d^2}$: the difference~$
\utilde{\bDelta}^{(n)}_{{\J}_1, {\J}_2}(\hat{\bth}^{(n)} + n^{-1/2}{\bf e}_i) - \utilde{\bDelta}^{(n)}_{{\J}_1, {\J}_2}(\hat{\bth}^{(n)})
$ 
then provides a consistent estimator of the $i$th column of $-\bUpsilon_{{\J}_1, {\J}_2, f}(\bth_0)$. %More sophisticated constructions also are possible: 
See \cite{HOP2006}  or \cite{Cassart2010} for more sophisticated estimation methods.

%For sake of implementation of our tests, we summarize, in  Algorithm~\ref{algorithm1} and Algorithm~\ref{algorithm2}, respectively, the algorithms of the  rank-based test procedures for the null of $\bth = \bth_0$ (with specified $\bth_0$) and the VAR order identification (with unspecified $\bth_0$). in Appendix~D.

\appendix\section*{Appendix F:  Algorithms for the center-outward rank-based tests}\label{algosec}

For the sake of implementation, we summarize, in  Algorithm~\ref{algorithm1} and Algorithm~\ref{algorithm2}, respectively, the algorithms of the  rank-based test procedures for the null of $\bth = \bth_0$ (with specified $\bth_0$) and the VAR order identification (with unspecified $\bth_0$).\bigskip

\begin{algorithm}[H]
\SetAlgoLined
\KwIn{a $d$-dimensional sample $\{\X_t; 1 \leq t \leq n\}$, the null VAR($p_0$) parameter $\bth_0$, the alternative VAR order $p$, nominal level $\alpha$.}
\KwOut{Indicator of the  rejection region (asymptotic critical value).}
%\begin{enumerate}

 Factorize $n$ into $n_Rn_S + n_0$ and  generate a ``regular grid" of $n_R n_S$ points over the unit ball $\mathbb{S}_d$.  

Specify $\J_1$ and $\J_2$, set the initial values  $\bepsilon_0$ and~$\X_{-p_0 + 1}^{(n)}, \ldots , \X_{0}^{(n)}$ all equal to zero, and compute residuals  $\ZZ_1^{(n)}(\bth_0),\ldots,\ZZ_n^{(n)}(\bth_0)$ recursively.

 Create an $n\times n$ matrix $\mathbf{D}$ with $(i,j)$ entry %on $i$-th row and $j$-th column being
 the squared Euclidean distance between~$\mathbf{Z}_i^{(n)}(\bth_0)$ and the $j$-th gridpoint. Based on that  matrix, compute $\{\F_\pm^{(n)} (\ZZ_t^{(n)}(\bth_0)); t=1,\ldots, n \}$ solving the optimal pairing problem in (\ref{Fpm0}),  using e.g.  the Hungarian algorithm.
 
 From $\F^{(n)}_\pm$, compute the center-outward ranks~(\ref{Ranks}), signs~(\ref{Signs}), and $\m_{\rm a}\n$.
 
 Compute $\Q^{(n)}_{\bth_0}$ in \eqref{defPQ},  $\utilde{\bGamma}_{i, \J_1, \J_2}^{(n)}(\bth_0)$ in \eqref{tildeGam1}, and $\utilde{\H}^{(n)}(\bth_0)$ in \eqref{Stilde}. Then combine these expressions into $\utilde{S}_{\J_1 \J_2}^{(n)}(\bth_0)$
 
Return $\1_{\{ \utilde{S}_{\J_1 \J_2}^{(n)}(\bth_0) < \chi_{d^2 p, \alpha}\}}$.

%\end{enumerate}
\caption{Center-outward rank-based tests for the null VAR with $\bth = \bth_0$}
\label{algorithm1}
\end{algorithm}

\begin{algorithm}[h!]
\setstretch{1.35}
\SetAlgoLined
\KwIn{a $d$-dimensional sample $\{\X_t; 1 \leq t \leq n\}$, nominal level $\alpha$.}
\KwOut{Order of the VAR model $p_0$.}
%\begin{enumerate}

 Factorize $n$ into $n_Rn_S + n_0$ and  generate a ``regular grid" of $n_R n_S$ points over the unit ball $\mathbb{S}_d$.

Specify $\J_1$ and $\J_2$, set $p_0 \gets 0$ and $\ZZ_t^{(n)} \gets \X_t^{(n)}$ for $t = 1, ..., n$.

 Create a $n\times n$ matrix $\mathbf{D}$ with $(i,j)$ entry %on $i$-th row and $j$-th column being
 the squared Euclidean distance between~$\mathbf{Z}_i^{(n)}$ and the $j$-th gridpoint. Based on that  matrix, compute $\{\F_\pm^{(n)} (\ZZ_t^{(n)}); t=1,\ldots, n \}$ solving the optimal pairing problem in (\ref{Fpm0}),  using e.g.  the Hungarian algorithm.
 
 From $\F^{(n)}_\pm$, compute the center-outward ranks~(\ref{Ranks}), signs~(\ref{Signs}).

 Compute $\tenq{\bGamma}_{1, \J_1, \J_2}^{(n)}$ by using \eqref{tildeGam1}, then, via \eqref{tildeW.white}, compute $\utilde{W}^{(n)}_{\J_1 \J_2}$.

 \uIf{$\utilde{W}^{(n)}_{\J_1 \J_2} < \big(F^{\chi^2}_{d^2}\big)^{-1}\!(1-\alpha)$}{
    \Return $p_0$\;
  }
  \Else{
 \Repeat{$\widehat{\utilde{W}}^{(n)}_{\J_1 \J_2}(\hat{\bth}^{(n)}) < \big(F^{\chi^2}_{d^2}\big)^{-1}\!(1-\alpha)$}{
      Set $p_0 \gets p_0 + 1$ \;
      Compute a root-$n$ consistent estimator $\hat{\bth}^{(n)}$ of $\bth_0$ using, e.g., the QMLE or center-outward R-estimator in \cite{HLL2019}\;
       Set the initial values  $\bepsilon_0$ and~$\X_{-p_0 + 1}^{(n)}, \ldots , \X_{0}^{(n)}$ all equal to zero, and compute residuals  $\ZZ_1^{(n)}(\hat{\bth}^{(n)}),\ldots,\ZZ_n^{(n)}(\hat{\bth}^{(n)})$ recursively\;
      Create a $n\times n$ matrix $\mathbf{D}$ with $(i,j)$ entry %on $i$-th row and $j$-th column being
 the squared Euclidean distance between~$\mathbf{Z}_i^{(n)}(\hat{\bth}^{(n)})$ and the $j$-th gridpoint. Based on that  matrix, compute $\{\F_\pm^{(n)} (\ZZ_t^{(n)}(\hat{\bth}^{(n)})); t=1,\ldots, n \}$\;
      From $\F^{(n)}_\pm$, compute the center-outward ranks~(\ref{Ranks}), signs~(\ref{Signs}) and $\m_{\rm a}\n$\;
      Compute $\M_{{\hat{\bth}^{(n)}}}$ in \eqref{defM}, and $\P_{{\hat{\bth}^{(n)}}}$ and $\Q^{(n)}_{{\hat{\bth}^{(n)}}}$ in \eqref{defPQ}, then   $\tenq{\bGamma}_{i, {\J}_1, {\J}_2}^{(n)}({\hat{\bth}^{(n)}})$ in \eqref{tildeGam1}. Finally, combine these expressions into $\utilde{\bDelta}^{(n)}_{{\J}_1, {\J}_2}(\hat{\bth}^{(n)})$\;
       For some chosen $\btau_1,\ldots,\btau_{p_0 d^2}$, compute $\utilde{\bDelta}^{(n)}_{{\J}_1, {\J}_2}(\hat{\bth}^{(n)} + n^{-1/2}\btau)$, then,  via~(\ref{asy.linear2}),~$\hat{\bUpsilon}_{{\J}_1, {\J}_2}^{(n)}$\;
       From $\utilde{\bDelta}^{(n)}_{{\J}_1, {\J}_2}(\hat{\bth}^{(n)})$ and $\hat{\bUpsilon}_{{\J}_1, {\J}_2}^{(n)}$, compute 
$\hat{\utilde{\bDelta}}^{(n)*}_{I; \J_1, \J_2}(\hat{\bth}^{(n)})$  and $\hat{\bLam}^{(n)*}_{I; \hat{\bth}^{(n)}}$
       using \eqref{defDeltastar} and \eqref{defLamstar} respectively, then $\widehat{\utilde{W}}^{(n)}_{\J_1 \J_2}(\hat{\bth}^{(n)})$ using \eqref{defW};
    }
  }

%\end{enumerate}
\caption{Center-outward rank-based sequential test for VAR order identification}
\label{algorithm2}
\end{algorithm}

\end{document}